\documentclass[11pt
              , a4paper
              , twoside
              , openright
              ]{report}

\usepackage{float} % lets you have non-floating floats

\usepackage{url} % for typesetting urls
\newfloat{fig}{thp}{lof}[chapter]
\floatname{fig}{Figure}

\title{The reverse mathematics of Cousin's lemma}
\author{Jordan Mitchell Barrett}

\usepackage[image,sms,bschonsmath]{vuwproject}

\supervisors{Rod Downey, Noam Greenberg}

\date{Friday 30th October 2020}

\usepackage[utf8]{inputenc}
\usepackage[T1]{fontenc}

\usepackage{pgfplots}
\usetikzlibrary{matrix}
\usetikzlibrary{arrows.meta}
\usetikzlibrary{calc}
\usetikzlibrary{cd}

\usepackage{tasks}
\usepackage{nicefrac}
\usepackage{mathtools}
\usepackage{amsmath}
\usepackage{amssymb}
\usepackage{amsthm}
\usepackage{enumitem}
\usepackage{hyperref}

\usepackage[
	backend=biber,%citestyle=alphabetic,%maxalphanames=1,
	bibstyle=jmbarrett,
	citestyle=alphabetic
]{biblatex}
\DeclareNolabel{\nolabel{\regexp{[\p{Z}\p{P}\p{S}\p{C}]+}}}
\theoremstyle{definition}
	\newtheorem{theorem}{Theorem}[section]
	\newtheorem{lemma}[theorem]{Lemma}
	\newtheorem{example}[theorem]{Example}
	\newtheorem{definition}[theorem]{Definition}
	
	\newtheorem{corollary}[theorem]{Corollary}
\newtheorem{proposition}[theorem]{Proposition}
\newtheorem{thmclaim}{Claim}[theorem]	% Claims inside proof of thm

\newtheorem{remark}[theorem]{Remark}

\newtheorem{axioms}[theorem]{Axioms}
\newtheorem{question}[theorem]{Question}

\newtheorem{qdefinition}[theorem]{``Definition''}
\newtheorem{algorithm}[theorem]{Algorithm}

\makeatletter
\def\thmhead@plain#1#2#3{%
	\thmname{#1}\thmnumber{\@ifnotempty{#1}{ }\@upn{#2}}%
	\thmnote{ {\the\thm@notefont#3}}}
\let\thmhead\thmhead@plain
\makeatother

	\setlist[enumerate,1]{label={(\roman*)}}
	\setlist[enumerate,2]{label={(\alph*)}}

\newcommand{\N}{\mathbb{N}}			% Natural numbers
\newcommand{\Z}{\mathbb{Z}}			% Integers
\newcommand{\Q}{\mathbb{Q}}			% Rational numbers
\newcommand{\R}{\mathbb{R}}			% Real numbers
\newcommand{\UI}{{[0,1]}}			% Unit interval
\newcommand{\A}{\mathcal{A}}
\newcommand{\B}{\mathcal{B}}
\newcommand{\F}{\mathcal{F}}
\newcommand{\K}{\mathcal{K}}
\newcommand{\Ll}{\mathcal{L}}		% Curly "L" for language
\newcommand{\M}{\mathcal{M}}
\newcommand{\Ss}{\mathcal{S}}
\newcommand{\proves}{\vdash}			% Turnstile for "proves"
\newcommand{\nproves}{\nvdash}			% Turnstile for "proves"
\renewcommand{\models}{\vDash}			% Double turnstile "models"
\newcommand{\ZF}{\mathsf{ZF}}			% Zermelo-Fraenkel set theory
\newcommand{\ZFC}{\mathsf{ZFC}}			% ZF + AC
\newcommand{\RCA}{\mathsf{RCA}_0}		% Recursive comprehension
\newcommand{\WKL}{\mathsf{WKL}_0}		% Weak Konig's lemma
\newcommand{\ACA}{\mathsf{ACA}_0}		% Arithmetical comprehension
\newcommand{\ATR}{\mathsf{ATR}_0}		% Arithmetic transfinite recursion
\newcommand{\PiCA}{\Pi^1_1\text{-}\mathsf{CA}_0}	% Pi-1-1 comprehension
\newcommand{\PA}{\mathsf{PA}}			% Peano arithmetic

\newcommand{\CS}{{2^\omega}}					% Cantor space 2^ω
\newcommand{\CSf}{{2^{<\omega}}}				% Cantor space 2^<ω
\newcommand{\substr}[2]{{{#1}{\upharpoonright}_{#2}}}	% Substring

\DeclareMathOperator{\dom}{dom}				% Domain
\DeclareMathOperator{\Pow}{\mathcal{P}}		% Power set

\newcommand{\abs}[1]{{\left\lvert #1 \right\rvert}}		% Absolute value
\newcommand{\ang}[1]{{\left\langle #1 \right\rangle}}	% Angle brackets
\newcommand{\ct}{{}^\frown}								% String concatenation
\newcommand{\defeq}{
	\mathrel{\mathrel{\mathop:}\mkern-1.2mu=}}			% := symbol
\newcommand{\cc}{\mathfrak{c}}			% Cardinality of the continuum
\DeclareMathOperator{\dif}{d\!}							% Differential d
\newcommand{\restr}[2]{{{#1}{\mid}_{#2}}}	% Restriction

\newcommand{\CLc}{\mathsf{CL}_\mathrm{c}}
\newcommand{\CLB}[1]{\mathsf{CL}_{\mathsf{B}#1}}

\newcommand{\REC}{\mathrm{REC}}
\newcommand{\ARITH}{\mathrm{ARITH}}

\newcommand{\chr}{\chi}
\DeclareMathOperator{\RS}{RS}
\DeclareMathOperator{\proj}{proj}
\newcommand{\Ls}{\Ll_2}%\mathrm{SOA}}

\setlist[enumerate]{itemsep=1pt}%,topsep=5pt}
\setlist[itemize]{itemsep=0pt}%,topsep=5pt}
\setlist{topsep=5pt}

\let\oldcite\cite
\let\cite\relax
\newcommand{\cite}[2][]{{\color{gray}\oldcite[#1]{#2}}}

\begin{document}

% Make the page numbering roman, until after the contents, etc.
\frontmatter

%%%%%%%%%%%%%%%%%%%%%%%%%%%%%%%%%%%%%%%%%%%%%%%%%%%%%%%

%%%%%%%%%%%%%%%%%%%%%%%%%%%%%%%%%%%%%%%%%%%%%%%%%%%%%%%

\begin{abstract}
	
Cousin’s lemma is a compactness principle that naturally arises when studying the gauge integral, a generalisation of the Lebesgue integral.	We study the axiomatic strength of Cousin's lemma for various classes of functions, using Friedman and Simpson's reverse mathematics in second-order arithmetic. We prove that, over $\RCA$:
\begin{enumerate}
	\item Cousin's lemma for continuous functions is equivalent to the system $\WKL$;
	
	\item Cousin's lemma for Baire 1 functions is at least as strong as $\ACA$;
	
	\item Cousin's lemma for Baire 2 functions is at least as strong as $\ATR$.
\end{enumerate}

%Cousin's lemma states that for every positive valued function
%there is a finite sequence

%can be viewed as a compactness principle

%In this report, we analyse

%Cousin's lemma is a compactness principle arising naturally in the study of gauge integration.

%We analyse the strength of various forms Cousin's lemma, from the point of view of reverse mathematics.

%In particular, we prove that Cousin's lemma for continuous functions is equivalent to the subsystem $\WKL$ of second-order arithmetic
%but for Baire 1 functions, it requires 

%Abstract:  Using reverse mathematics, we analyse the strength of Cousin’s lemma for various classes of functions. Joint work with Downey and Greenberg.

%Abstract: Cousin’s lemma is a compactness principle that naturally arises when studying the gauge integral, a generalisation of the Lebesgue integral. We attempt to find the optimal proof of Cousin’s lemma for various functions, using the toolkit of reverse mathematics. Joint work with Downey and Greenberg.

\end{abstract}

%%%%%%%%%%%%%%%%%%%%%%%%%%%%%%%%%%%%%%%%%%%%%%%%%%%%%%%

\maketitle

\tableofcontents

% we want a list of the figures we defined
\listof{fig}{Figures}

%%%%%%%%%%%%%%%%%%%%%%%%%%%%%%%%%%%%%%%%%%%%%%%%%%%%%%%

\mainmatter

%%%%%%%%%%%%%%%%%%%%%%%%%%%%%%%%%%%%%%%%%%%%%%%%%%%%%%%

% individual chapters included here
\chapter{Introduction}\label{C:intro}

Before the 17th century, mathematics essentially comprised arithmetic, geometry and elementary algebra, and generally only dealt with finite objects \cite{evesIntroductionHistoryMathematics1969}.
Proofs were almost always \textit{constructive}---a statement would be proved by explicitly constructing a witness.
Being largely motivated by physics, mathematics was concerned primarily with calculation, and therefore algorithms took centre stage \cite{metakidesIntroductionNonrecursiveMethods1982}.

%Before the 17th century, mathematics essentially comprised arithmetic, geometry and elementary algebra. As a result, it was largely finitary, and entirely constructive and computable. A proof of a statement essentially amounted to an algorithm to solve some problem, such as finding the roots of a polynomial.

%algorithmic focus - motivated by physics etc where calculation important

The development of calculus in the 17th century represented the first signs of departure from this.
The ideas were present in Archimedes' method of exhaustion and Cavalieri's method of indivisibles \cite{evesIntroductionHistoryMathematics1969}, but Newton and Leibniz systematised this, manipulating infinite and infinitesimal quantities as if they were numbers. These new methods proved revolutionary
%Newton and Leibniz's development of calculus in the 17th century represented a radical departure from this. For the first time, the ideas of infinity and infinitesimals were used in an essential way, and moreover, this proved extremely useful
throughout mathematics and physics. There was some concern about the rigour of such methods, and this was not fully abated until the 19th century, when Cauchy, Bolzano, and Weierstrass replaced infinitesimals with more rigorous $\varepsilon$-$\delta$ definitions \cite{cauchyCoursAnalyseEcole1821, bolzanoReinAnalytischerBeweis1817, schwarzLecturesKarlWeierstrass1861}.

%{\color{red} A systematic <> of Archimedes' method of exhaustion and Cavalieri's
%Infinities/infinitesimals were treated/manipulated as regular quantities}

As calculus flourished into real analysis, the techniques used became gradually less constructive \cite{metakidesIntroductionNonrecursiveMethods1982}. Early analytical proofs would often implicitly appeal to the infinite pigeonhole principle, and (weak forms of) the axiom of choice \cite{mendiolaSignificanceAxiomChoice2016}, thereby proving the existence of objects without actually constructing them. Despite some backlash, this trend towards nonconstructivism only continued, as analysis was later abstracted to topology and descriptive set theory.

In the 1870s, while studying a problem in topology, Georg Cantor formulated the concept of ordinals \cite{cantorUeberUnendlicheLineare1883},
leading to the creation of set theory. Cantor was the first to systematically study infinity: famously, he showed in 1874 that the real numbers $\R$ cannot be put into bijection with the natural numbers $\N$, thus demonstrating that there are different sizes of infinity \cite{cantorUeberEigenschaftInbegriffs1874}.
As set theory developed, paradoxes arose (most notably Russell's), and the need for a careful and rigorous foundation for mathematics became clear. One such foundation was provided by $\ZFC$ in the 1920s \cite{zermeloUberGrenzzahlenUnd1930a}.

Cantor's work provided new impetus to mathematical logic, a small subfield of mathematics developed by Boole, De Morgan, and Peano in the mid-to-late 1800s \cite{booleInvestigationLawsThought1854, demorganFormalLogic1847, peanoArithmeticesPrincipiaNova1889}. Around this time, the ideas of computation, mathematical truth and mathematical proof were formalised for the first time. By the 1930s, logic was a thriving area of mathematics---highlights included G\"odel's completeness \cite{godelUberVollstandigkeitLogikkalkuls1929} and incompleteness theorems \cite{goedelUeberFormalUnentscheidbare1931}, Turing's negative solution to the \textit{Entscheidungsproblem} \cite{turingComputableNumbersApplication1937}, Tarski's development of model theory \cite{vaughtAlfredTarskiWork1986}, and Hilbert's work on geometry \cite{hilbertGrundlagenGeometrie1899} and proof theory \cite{hilbertGrundlagenMathematik1934}.

A later development in logic was reverse mathematics, initiated by Harvey Friedman in the late 1960s \cite{friedmanSubsystemsSetTheory1967, friedmanBarInductionPi1969}. Reverse mathematics asks, for a given theorem of mathematics $\varphi$, ``what axioms are really necessary to prove $\varphi$?'' More broadly, it studies the logical implications between foundational principles of mathematics. An early example was the discovery of non-Euclidean geometries, thereby proving the independence of the parallel postulate from Euclid's other axioms \cite{lobachevskyConciseOutlineFoundations1829, bolyaiAppendixScientiamSpatii1832}. Another early result, more in the style of reverse mathematics, was the demonstration that over $\ZF$, the axiom of choice, Zorn's lemma, and the well-ordering principle are all pairwise equivalent \cite{birkhoffLatticeTheory1940, fraenkelFoundationsSetTheory1958, traylorEquivalenceAxiomChoice1962}.

Traditionally, reverse mathematics is done in second-order arithmetic, in which there are two types of objects: natural numbers $n, m, k, \ldots$, and sets of natural numbers $A, B, C, \ldots$, and quantification is allowed over both types of objects. Restricting oneself to natural numbers may seem unnecessary limiting, but this is not so. In fact, most mathematics deals with countable or ``essentially countable'' objects (such as separable metric spaces), and so can be formalised in second-order arithmetic. This includes virtually all ``classical'' mathematics, or that taught in undergraduate courses \cite[xiv]{simpsonSubsystemsSecondOrder2009}.

In practice, reverse mathematics involves attempting to prove a theorem $\varphi$ of ``ordinary'' mathematics in a weak subsystem $\Ss$ of second-order arithmetic. But, supposing we can do this, how do we know we've found the optimal (weakest) system? The empirical phenomenon is thus:
\begin{center}
	\textit{``When the theorem is proved from the right axioms, the axioms can be proved from the theorem.''}
\end{center}\vspace*{-5.7mm}
\begin{flushright}
	\textit{---Harvey Friedman} \cite{friedmanSystemsSecondOrder1974}
\end{flushright}
This is the ``reverse'' part of reverse mathematics. Having proved $\varphi$ from $\Ss$, to show this is optimal, we want to demonstrate a \textit{reversal} of $\varphi$: a proof of $\Ss$ from $\varphi$. This means that $\varphi$ cannot be proved in a weaker system $\Ss'$, because if it could, then $\Ss'$ would also prove $\Ss$ via $\varphi$, meaning $\Ss'$ is not actually a weaker system after all. Practically speaking, reversals are only possible assuming a weak base system $\B$ (i.e.\ it is really a proof of $\Ss$ from $\B+\varphi$).

The utility of reverse mathematics is abundant. Apart from its obvious use in finding the ``best'' proof of a given statement $\varphi$, it also gives us a way to quantify how nonconstructive or noncomputable $\varphi$ is. The idea is that stronger subsystems correspond to more nonconstructive power, %so the weaker the system $\Ss$ that $\varphi$ can be proved in, the more constructive $\varphi$ is \cite{friedmanCountableAlgebraSet1983}.
so the ``constructiveness'' of $\varphi$ is inversely proportional to the strength of the systems $\Ss$ in which $\varphi$ can be proved \cite{friedmanCountableAlgebraSet1983}.
Similarly, many theorems guarantee a solution to a given problem---reverse mathematics then tells us how complex the solution could be relative to the problem, which can be made precise in terms of computability. For example, in his thesis \cite{miletiPartitionTheoremsComputability2004}, Mileti proved the Erd\H os--Rado canonical Ramsey theorem is equivalent to the system $\ACA$. From the proof, he extracted new bounds on the complexity of the homogeneous set, improving the classical bounds obtained by Erd\H os and Rado.

Here is an example of reverse mathematics in ring theory. The usual way to prove that every commutative ring has a prime ideal is to prove that it has a maximal ideal (Krull's theorem), and then prove every maximal ideal is prime. However, Friedman, Simpson and Smith showed that the existence of maximal ideals is equivalent to the system $\ACA$, whereas the existence of prime ideals is equivalent to the strictly weaker system $\WKL$ \cite{friedmanCountableAlgebraSet1983}. This shows the usual proof strategy is not optimal---there is a ``better'' way to prove the existence of prime ideals, which doesn't require the stronger assumption that maximal ideals exist.

In this report, we examine the reverse-mathematical content of Cousin's lemma, a particular statement in analysis. Cousin's lemma can be viewed as a kind of compactness principle, asserting that every positive valued function $\delta\colon \UI \to \R^+$ has a \textit{partition}---a finite sequence $t_0, \ldots, t_{n-1}$ such that the open balls $B \big( t_i, \delta(t_i) \big)$ cover $\UI$. In particular, we establish the following original results over the weak base theory $\RCA$:
\begin{enumerate}
	\item Cousin's lemma for continuous functions is equivalent to the system $\WKL$;
	
	\item Cousin's lemma for Baire 1 functions is at least as strong as $\ACA$;
	
	\item Cousin's lemma for Baire 2 functions is at least as strong as $\ATR$.
\end{enumerate}

\section*{Notational conventions}

Throughout this report, we abide to the following notational conventions.

\begin{itemize}
	\item We let $\N = \{ 0, 1, 2, \ldots \}$ be the set of \textit{nonnegative} integers, following the usual practice in logic.
	
	\item We will often use $\bar{x}$ to notate a tuple $(x_1,\ldots,x_d)$, where the length should be clear from context.
	
	\item We may use $A^\complement = \{ x: x \notin A \}$ to denote the (absolute) complement of a set $A$, particularly for sets of natural numbers.
	
	\item For mathematical statements $\varphi$ and $\psi$, we use $\varphi \proves \psi$ (``$\varphi$ proves $\psi$'') to mean there is a proof of $\psi$ from $\varphi$. This notation extends to formal systems, e.g.\ $\Ss \proves \varphi$ means there is a proof of $\varphi$ in the formal system $\Ss$.
	
	\item For a statement $\varphi$ and a structure $\M$, we use $\M \models \varphi$ (``$\M$ models $\varphi$'') to mean the statement $\varphi$ is true in $\M$. Similarly, $\M \models \Ss$ means that all axioms of the formal system $\Ss$ are true in $\M$.
\end{itemize}
\chapter{Integration and Cousin's lemma}

The main object of study in this report is \textit{Cousin's lemma}, a compactness principle phrased in terms of positive real-valued functions, rather than open covers. Cousin's lemma arises naturally in the study of the \textit{gauge integral}, a generalisation of the Riemann and Lebesgue integrals due to Kurzweil \cite{kurzweilGeneralizedOrdinaryDifferential1957} and Henstock \cite{henstockTheoryIntegration1963}. In this chapter, we review Riemann integration, before generalising to gauge integration and defining Cousin's lemma.

%Refs:
%
%Kurtz--Swartz, p. 139 onwards
%Gordon, Ch 9

\section{Riemann integration}

\begin{figure}\label{fig:riemann}
\centering

\pgfmathdeclarefunction{f}{1}{%
	\pgfmathparse{
		  (#1<0.3)            * (0.2 + #1^2)
		+ (#1>=0.3)*(#1<0.5)  * (0.29 - (#1-0.3))
		+ (#1>=0.5)*(#1<0.85) * (0.09 -(#1-0.5) +10*(#1-0.5)^2 -20*(#1-0.5)^3))
		+ (#1>=0.85)          * (0.1075 + (#1-0.85))
			 }%
}
\def\Px{{0,     0.14,    0.27,     0.4,     0.55,     0.7,    0.87,     1}}
\def\Pt{{  0.08,     0.2,     0.34,    0.49,     0.63,    0.8,     0.94  }}
\begin{tikzpicture}% coordinates
\begin{axis}[
	width=\linewidth,height=0.2\textheight,
	xmin=-0.05,xmax=1.05,ymin=0,
	xtick=\Px,xtick style={draw=black,thick},
	xticklabels={$0=x_0\phantom{{}=0}$,$\vphantom{0}x_1$,$\vphantom{0}x_2$,$\vphantom{0}x_3$,$\vphantom{0}x_4$,$\vphantom{0}x_5$,$\vphantom{0}x_6$,$\phantom{1={}}x_7=1$},
%	minor xtick=\Pt,minor x tick style={shape=circle,blue},
	axis x line=bottom, x axis line style={thick,latex-latex},
	y axis line style={draw=none},
	ytick style={draw=none},yticklabels=none,
	axis on top, clip=false,
]
	% Rectangles
	\foreach \i in {0,...,6}{
		\pgfmathsetmacro\x{\Px[\i]}
		\pgfmathsetmacro\t{\Pt[\i]}
		\pgfmathsetmacro\y{\Px[\i+1]}
		\edef\temp{\noexpand%
			%\node at (axis cs:0.2*\i,0.1){\x\t\y};%
			\draw[fill=gray,fill opacity=0.2,draw=gray,rounded corners=1] (axis cs:\x,0) rectangle (axis cs:\y,{f(\t)});
		}\temp
	}
	
	% Function
	\addplot[samples=300,domain=0:1,thick]{f(x)};
	
	% Tag points t_i, vertical dashed lines
	\foreach \i in {0,...,6}{
		\pgfmathsetmacro\t{\Pt[\i]}
		\edef\temp{\noexpand
			\draw[dashed,gray] (axis cs:\t,0) -- (axis cs:\t,{f(\t)});
		}\temp
		\edef\temp{\noexpand
			\node[fill=black,circle,scale=0.3,label={[label distance=0.2ex]270:$t_\i$}] at (axis cs:\t,0) {};
		}\temp
		\edef\temp{\noexpand
			\node[fill=black,circle,scale=0.3] at (axis cs:\t,{f(\t)}) {};
		}\temp
	}
\end{axis}
\end{tikzpicture}

\caption{A Riemann sum of a continuous function over a partition of size 7.}
\end{figure}
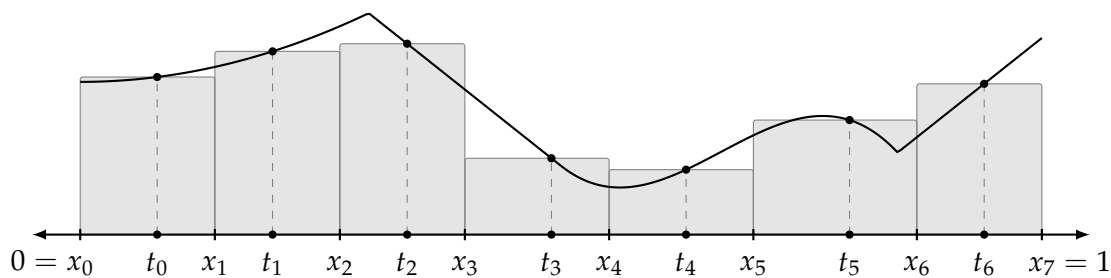

The basic idea of Riemann integration is thus: approximate the area under a curve by a series of rectangles, as in Figure \ref{fig:riemann}. As we increase the number of rectangles, and decrease their width, we hope that this approximation becomes closer and closer to the true area. Here, we will only consider integration over the unit interval $\UI$.

\begin{definition}
	A \textit{tagged partition of $\UI$}\index{tagged partition} is a finite sequence $$P\ =\ \ang{0 = x_0 < t_0 < x_1 < t_1 < \cdots < t_{n-1} < x_n = 1}.$$ We call $n$ the \textit{size} of $P$.
\end{definition}

A tagged partition $P = \ang{x_i,t_i}$ should be interpreted as follows. The $x_i$ are the \textit{partition points} at which the interval $\UI$ is split, and within each subinterval or block $[x_i,x_{i+1}]$, we choose a \textit{tag point} $t_i$. When using $P$ to approximate the area underneath a function $f$, each subinterval $[x_i,x_{i+1}]$ will serve as the base of a rectangle of height $f(t_i)$. This is illustrated in Figure \ref{fig:riemann}.

\begin{definition}
	Let $f\colon \UI \to \R$ be a function, and $P = \ang{x_i,t_i}$ a partition of size $n$. The \textit{Riemann sum of $f$ over $P$} is $$\RS(f,P) = \sum_{i=0}^{n-1} f(t_i)[x_{i+1}-x_i]$$
\end{definition}

\begin{definition}[\cite{riemannUberDarstellbarkeitFunction1854}]\label{defn:riem-int}
	A function $f\colon \UI \to \R$ is \textit{Riemann integrable} if there exists $K \in \R$ such that, for every $\varepsilon > 0$, there exists $\delta>0$ with $\abs{\RS(f,P) - K} < \varepsilon$ whenever each block of $P$ has size $<\delta$. In this case, we say that $K$ is the \textit{Riemann integral of $f$}.
\end{definition}

The Riemann integral can deal with virtually all functions which one may want to integrate in practice. Riemann (and later Lebesgue) gave a characterisation of exactly which bounded functions are Riemann integrable:

\begin{proposition}[\cite{brownProofLebesgueCondition1936, birkhoffSourceBookClassical1973}]
	A bounded function $f\colon \UI \to \R$ is Riemann integrable if and only if its set of discontinuities has (Lebesgue) measure zero.
\end{proposition}

In particular, every continuous function is Riemann integrable. That said, it is not difficult to construct functions which are not Riemann integrable.

\begin{proposition}[\cite{dirichletConvergenceSeriesTrigonometriques1829}]\label{prop:not-riem-intg}
	There are functions which are not Riemann integrable.
\end{proposition}

\begin{proof}
	The characteristic function $\chr_\Q$ of $\Q$, also known as Dirichlet's function, provides an example. Concretely, $\chr_\Q\colon \UI \to \R$ is defined
	\begin{equation*}
		\chr_\Q(x) = \begin{cases}
			0 & \text{ if } x \text{ is irrational} \\
			1 & \text{ if } x \text{ is rational} \\
		\end{cases}
	\end{equation*}
	We will show $\chr_\Q$ is not Riemann integrable. Pick any $K \in \R$, and let $\varepsilon = \nicefrac{1}{3}$. Then, for any $\delta>0$, pick $n > \nicefrac{1}{\delta}$, and consider the partitions $P = \ang{x_i,t_i}$, $P' = \ang{x'_i,t'_i}$ of size $n$ defined by:\vspace*{-1ex}
	\begin{align*}
		x_i &= x'_i = \nicefrac{i}{n}; \\
		t_i &= \text{ some irrational point in } [\nicefrac{i}{n},\nicefrac{(i+1)}{n}]; \\
		t'_i &= \text{ some rational point in } [\nicefrac{i}{n},\nicefrac{(i+1)}{n}].
	\end{align*}
	Then, each block of $P$, $P'$ has size $<\delta$, and $\RS(\chr_\Q,P)=0$ while $\RS(\chr_\Q,P')=1$. Thus, $K$ cannot be within $\varepsilon = \nicefrac{1}{3}$ of both.
\end{proof}

\section{Gauge integration}

So, what failed when trying to integrate $\chr_\Q$? Morally, since almost all real numbers in $\UI$ are irrational (in the measure-theoretic sense), the integral of $\chr_\Q$ ought to be equal to zero. There were many attempts to solve this, the most famous being Lebesgue's measure theory \cite{lebesgueIntegraleLongueurAire1902, lebesgueLeconsIntegrationRecherche1904}. However, Lebesgue integration is not without its issues---in particular, there are derivatives which are not Lebesgue integrable \cite{gordonNonabsoluteIntegrationWorth1996}.

Attempting to remedy this, Denjoy defined an integral which could handle all derivatives \cite{denjoyExtensionIntegraleLebesgue1912}. Shortly after, Luzin \cite{luzinProprietesIntegraleDenjoy1912} and Perron \cite{perronUberIntegralbegriff1914} gave equivalent characterisations of Denjoy's integral. However, all these definitions were complex and highly nonconstructive, making Denjoy's integral impractical for applications \cite{gordonNonabsoluteIntegrationWorth1996}.

In 1957, Kurzweil defined the \textit{gauge integral}, a generalisation of Denjoy's integral. He formulated it in elementary terms similar to the Riemann integral \cite{kurzweilGeneralizedOrdinaryDifferential1957}, thus avoiding the complications of measure theory. Later, Henstock systematically developed the theory of the gauge integral \cite{henstockTheoryIntegration1963}--- as a result, it is sometimes known as the Henstock--Kurzweil integral.
Kurzweil's ingenious solution was to allow the parameter $\delta$ in Definition \ref{defn:riem-int} to be a \textit{variable}, rather than a constant. In effect, this ensures that some partitions are not allowed, such as $P'$ in the proof of Proposition \ref{prop:not-riem-intg}. %We do this using the concept of a gauge:

\begin{definition}[\cite{gordonIntegralsLebesgueDenjoy1994}]
	A \textit{gauge} is a strictly positive-valued function $\delta\colon \UI \to \R^+$.
\end{definition}

The idea is that $\delta$ tells us how \textit{fine} our partition needs to be at any point. At points $x$ where the function is highly discontinuous, or varies greatly, we could make sure $\delta(x)$ is small, so that we only consider partitions which are divided finely enough around $x$.

\begin{definition}[\cite{gordonIntegralsLebesgueDenjoy1994}]
	Given a gauge $\delta$, a partition $P = \ang{x_i,t_i}$ is \textit{$\delta$-fine} if, for any $i<n$, the open ball $B(t_i,\delta(t_i))$ contains $(x_i,x_{i+1})$.
\end{definition}

\begin{definition}[\cite{kurzweilGeneralizedOrdinaryDifferential1957}]\label{defn:gauge-intg}
	A function $f\colon \UI \to \R$ is \textit{gauge integrable} if there exists $K \in \R$ such that, for every $\varepsilon > 0$, there exists a gauge $\delta\colon \UI \to \R^+$ with $\abs{\RS(f,P) - K} < \varepsilon$ whenever $P$ is $\delta$-fine. In this case, we say that $K$ is the \textit{gauge integral of $f$}.
\end{definition}

Note that if $\delta(x) = k$ is constant, a partition $P$ is $\delta$-fine if and only if the blocks of $P$ have size $< 2k$. Thus, Riemann integration is a special case of gauge integration, where we only allow constant gauges. It follows that every Riemann integrable function is gauge integrable.\footnote{Furthemore, the Riemann integral and gauge integral of $f$ will have the same value.} However, the converse does not hold, as we now see:

\begin{proposition}[\cite{kurtzTheoriesIntegrationIntegrals2004}]
	There are functions which are gauge integrable, but not Riemann integrable.
\end{proposition}

\begin{proof}
	Dirichlet's function $\chr_\Q$ is again an example. We saw in Proposition \ref{prop:not-riem-intg} that $\chr_\Q$ is not Riemann integrable---we now show that it is gauge integrable, with integral $K=0$. Pick any $\varepsilon>0$, and let $\Q = \{ q_0, q_1, q_2, \ldots \}$ enumerate the rationals. Define $\delta\colon \UI \to \R^+$ by $\delta(q_m)=2^{-m-2}\,\varepsilon$, and $\delta(x)=1$ for all irrational $x$. Now suppose $P = \ang{x_i,t_i}$ is $\delta$-fine. Then,
	\begin{align*}
		\RS(\chr_\Q,P)\quad &=\quad \sum_{i=0}^{n-1}\ \chr_\Q(t_i)[x_{i+1}-x_i] \\
		&=\quad \sum_{\mathclap{\substack{i<n\\t_i\text{ rational}}}}^{\hphantom{n-1}}\ \chr_\Q(t_i)[x_{i+1}-x_i]\quad +\quad \sum_{\mathclap{\substack{i<n\\t_i\text{ irrational}}}}^{\hphantom{n-1}}\ \chr_\Q(t_i)[x_{i+1}-x_i] \\
		&=\quad \sum_{\mathclap{\substack{i<n\\t_i\text{ rational}}}}^{\hphantom{n-1}}\ (x_{i+1}-x_i)\quad +\quad 0 \\
		&\leq\quad \sum_{\mathclap{\substack{i<n\\t_i=q_{m_i}}}}^{\hphantom{n-1}}\ 2^{-m_i-1}\,\varepsilon\quad <\quad \varepsilon \\[-7ex]
	\end{align*}
\end{proof}

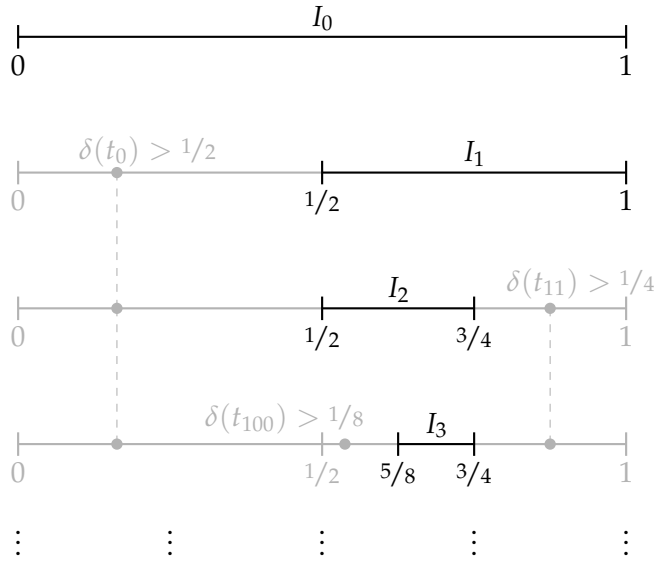
\begin{figure}[t]
\centering
	\begin{tikzpicture}[x=20mm,y=12mm]%[scale=0.8]%[y=20mm]

\definecolor{gg}{gray}{0.7}
\def\he{0.13}
\def\ya{-1.5}
\def\yb{-3}
\def\yc{-4.5}
\def\yd{-5.5}

\draw[thick] (0,0) -- (4,0);
\draw[thick] (0,\he) -- (0,-\he);
\draw[thick] (4,\he) -- (4,-\he);
\node[below,yshift=-1mm] at (0,0) {0};
\node[below,yshift=-1mm] at (4,0) {1};

\draw[thick,gg] (0,\ya) -- (2,\ya);
\draw[thick] (2,\ya) -- (4,\ya);
\draw[thick,gg] (0,\ya+\he) -- (0,\ya-\he);
\draw[thick] (2,\ya+\he) -- (2,\ya-\he);
\draw[thick] (4,\ya+\he) -- (4,\ya-\he);
\node[below,yshift=-1mm,gg] at (0,\ya) {0};
\node[below,yshift=-1mm] at (2,\ya) {$\nicefrac{1}{2}$};
\node[below,yshift=-1mm] at (4,\ya) {1};

\node[circle,gg,fill=gg,inner sep=1.5pt] at (0.65,\ya) {};
%\draw[fill=black] (0.65,\ya) circle (0.05);
\node[above,yshift=-0.5mm,gg] at (0.85,\ya) {$\delta(t_0)>\nicefrac{1}{2}$};

\draw[thick,gg] (0,\yb) -- (2,\yb);
\draw[thick] (2,\yb) -- (3,\yb);
\draw[thick,gg] (3,\yb) -- (4,\yb);
\draw[thick,gg] (0,\yb+\he) -- (0,\yb-\he);
\draw[thick] (2,\yb+\he) -- (2,\yb-\he);
\draw[thick] (3,\yb+\he) -- (3,\yb-\he);
\draw[thick,gg] (4,\yb+\he) -- (4,\yb-\he);
\node[below,yshift=-1mm,gg] at (0,\yb) {0};
%		\node[below,yshift=-1mm] at (1,\yb) {$\nicefrac{1}{4}$};
\node[below,yshift=-1mm] at (2,\yb) {$\nicefrac{1}{2}$};
\node[below,yshift=-1mm] at (3,\yb) {$\nicefrac{3}{4}$};
\node[below,yshift=-1mm,gg] at (4,\yb) {1};

\node[circle,gg,fill=gg,inner sep=1.5pt] at (0.65,\yb) {};
%\draw[fill=black] (0.65,\yb) circle (0.05);
\draw[dashed,thin,gg] (0.65,\ya) -- (0.65,\yb);
\node[circle,gg,fill=gg,inner sep=1.5pt] at (3.5,\yb) {};
\node[above,yshift=0mm,gg] at (3.7,\yb) {$\delta(t_{11})>\nicefrac{1}{4}$};

\draw[thick,gg] (0,\yc) -- (2.5,\yc);
\draw[thick] (2.5,\yc) -- (3,\yc);
\draw[thick,gg] (3,\yc) -- (4,\yc);
\draw[thick,gg] (0,\yc+\he) -- (0,\yc-\he);
\draw[thick,gg] (2,\yc+\he) -- (2,\yc-\he);
\draw[thick] (2.5,\yc+\he) -- (2.5,\yc-\he);
\draw[thick] (3,\yc+\he) -- (3,\yc-\he);
\draw[thick,gg] (4,\yc+\he) -- (4,\yc-\he);
\node[below,yshift=-1mm,gg] at (0,\yc) {0};
%		\node[below,yshift=-1mm] at (1,\yc) {$\nicefrac{1}{4}$};
\node[below,yshift=-1mm,gg] at (2,\yc) {$\nicefrac{1}{2}$};
\node[below,yshift=-1mm] at (2.5,\yc) {$\nicefrac{5}{8}$};
\node[below,yshift=-1mm] at (3,\yc) {$\nicefrac{3}{4}$};
\node[below,yshift=-1mm,gg] at (4,\yc) {1};

\node[circle,gg,fill=gg,inner sep=1.5pt] at (0.65,\yc) {};
%\draw[fill=black] (0.65,\yc) circle (0.05);
\draw[dashed,thin,gg] (0.65,\yb) -- (0.65,\yc);
\node[circle,gg,fill=gg,inner sep=1.5pt] at (3.5,\yc) {};
%\draw[fill=black] (0.65,\yc) circle (0.05);
\draw[dashed,thin,gg] (3.5,\yb) -- (3.5,\yc);
\node[circle,gg,fill=gg,inner sep=1.5pt] at (2.15,\yc) {};
\node[above,yshift=0mm,gg] at (1.75,\yc) {$\delta(t_{100})>\nicefrac{1}{8}$};

\foreach \i in {0,...,4}{
	\node at (\i,\yd) {$\vdots$};}

%\draw[thick,red] (0,0) -- (4,0);
%\draw[thick,red] (0,\he) -- (0,-\he);
%\draw[thick,red] (4,\he) -- (4,-\he);
\node[above,yshift=-0.5mm] at (2,0) {$I_0$};
%
%\draw[thick,red] (2,\ya) -- (4,\ya);
%\draw[thick,red] (2,\ya+\he) -- (2,\ya-\he);
%\draw[thick,red] (4,\ya+\he) -- (4,\ya-\he);
\node[above,yshift=-0.5mm] at (3,\ya) {$I_1$};
%
%\draw[thick,red] (2,\yb) -- (3,\yb);
%\draw[thick,red] (2,\yb+\he) -- (2,\yb-\he);
%\draw[thick,red] (3,\yb+\he) -- (3,\yb-\he);
\node[above,yshift=-0.5mm] at (2.5,\yb) {$I_2$};
\node[above,yshift=-0.5mm] at (2.75,\yc) {$I_3$};
\end{tikzpicture}
\caption{The recursive interval-splitting process used in the proof of Cousin's lemma.}
\end{figure}

\section{Cousin's lemma}

If there were a gauge $\delta$ with no $\delta$-fine partition, then Definition \ref{defn:gauge-intg} could be vacuously satisfied by choosing this gauge. This would present a problem: every $K \in \R$ would then witness that every $f$ is gauge integrable, so we could not uniquely define the value of the gauge integral. Cousin's lemma states that this situation cannot happen. It is originally due to Cousin \cite{cousinFonctionsVariablesComplexes1895}, who proved the statement in a radically different form.

\begin{lemma}[(Cousin's lemma) \cite{cousinFonctionsVariablesComplexes1895, kurtzTheoriesIntegrationIntegrals2004}]\label{lem:cl-og}
	Every gauge $\delta\colon \UI \to \R^+$ has a $\delta$-fine partition.
\end{lemma}

\begin{proof}[Proof\  {\ \normalfont\cite{leeHenstockKurzweilIntegrationEuclidean2011}}]
	If there is $t \in \UI$ with $\delta(t) > 1$, then the partition $\ang{0,t,1}$ is $\delta$-fine, so we are done.
	
	Otherwise, split $\UI$ into halves $[0,\nicefrac{1}{2}]$, $[\nicefrac{1}{2},1]$, and ask if there are $t_0 \in [0,\nicefrac{1}{2}]$, $t_1 \in [\nicefrac{1}{2},1]$ with $\delta(t_0), \delta(t_1) > \nicefrac{1}{2}$. If such a $t_0$ exists, we don't need to split further, since $\delta(t_0)$ covers the subinterval $[0,\nicefrac{1}{2}]$.
	
	If no such $t_0$ exists, we split $[0,\nicefrac{1}{2}]$ into halves $[0,\nicefrac{1}{4}]$, $[\nicefrac{1}{4},\nicefrac{1}{2}]$, and ask if there are $t_{00} \in [0,\nicefrac{1}{4}]$, $t_{01} \in [\nicefrac{1}{4},\nicefrac{1}{2}]$ with $\delta(t_{00}), \delta(t_{01}) > \nicefrac{1}{4}$. Keep repeating this process, and do similar on the side of $t_1$.
	
	We claim this procedure must eventually terminate. Suppose it did not --- then, there is a nested sequence $\UI = I_0 \supsetneq I_1 \supsetneq I_2 \supsetneq \cdots$ of closed intervals, each half the size of the previous. Since $\UI$ is compact, we can pick $r \in \bigcap_{n=0}^\infty I_n$. But $\delta(r)$ is positive, so $\delta(r) > 2^{-n}$ for sufficiently large $n$. Then, the procedure would have terminated at stage $n$, since we would have found $t = r \in I_n$.
\end{proof}

Among other things, Cousin's lemma implies that the value of a gauge integral is unique, if it exists:

\begin{corollary}[\cite{kurtzTheoriesIntegrationIntegrals2004}]
	For gauge integrable $f\colon \UI \to \R$, there is a unique $K$ witnessing the integrability of $f$.
	%If $K_1 = \int_\mathrm{G} f$ and $K_2 = \int_\mathrm{G} f$, then $K_1=K_2$.
\end{corollary}

\begin{proof}
	By contradiction, suppose $K_1 \neq K_2$ both witness that $f$ is gauge integrable. Let $\varepsilon = \nicefrac{\abs{K_1-K_2}}{3}$. By assumption, there are gauges $\delta_1$, $\delta_2$ witnessing $K_1$, $K_2$ respectively for this choice of $\varepsilon$. Let $\delta(x) \defeq \min \{ \delta_1(x), \delta_2(x) \}$ be the pointwise minimum: this is also a gauge. By Cousin's lemma, there is a $\delta$-fine partition $P$, which must also be $\delta_1$-fine and $\delta_2$-fine by definition of $\delta$. So $\abs{\RS(f,P) - K_1} < \varepsilon$ and $\abs{\RS(f,P) - K_2} < \varepsilon$, a contradiction.
\end{proof}

We can view Cousin's lemma as a kind of compactness principle. Effectively, it asserts that the open cover $\big\{ B \big( t,\delta(t) \big)\!: t \in \UI \big\}$ has a finite subcover, corresponding to the tag points of a $\delta$-fine partition.

%\begin{proof}
%	Starting with $I = \UI$ and $P_0 = \ang{x_0=0, x_1=1}$, we recursively build a $\delta$-fine partition by the following procedure:
%	\begin{enumerate}
%		\item Let $I = [a,b]$. If there is no $t \in I$ with $\delta(t) > b-a$, then:
%		\begin{enumerate}
%			\item Change the segment $P = \ang{\, \ldots,\, x_j = a,\, x_{j+1} = b,\, \ldots\, }$ to $$P' = \ang{\, \ldots,\, x_j = a,\, {\color{red} x_{j+1} = (a+b)/2,}\, x_{j+2} = b,\, \ldots\, }$$
%			
%			\item Divide $I$ into intervals $I_1 = [a,(a+b)/2]$, $I_2 = [(a+b)/2,b]$;
%			
%			\item Continue the procedure on $I_1$ and $I_2$.
%		\end{enumerate}
%		
%		\item Else, if there is $t \in I$ with $\delta(t) > b-a$, then the segment $$P = \ang{\, \ldots,\, x_j = a,\, t_j = t,\, x_{j+1} = b,\, \ldots\, }$$ is $\delta$-fine.
%	\end{enumerate}
%	
%	We claim this procedure must eventually terminate. Suppose it did not --- then, there is a nested sequence $\UI = I_0 \supsetneq I_1 \supsetneq I_2 \supsetneq \cdots$ of closed intervals, each half the size of the previous. Since $\UI$ is compact, we can pick $r \in \bigcap_{n=0}^\infty I_n$. But $\delta(r)$ is positive, so $\delta(r) > 2^{-n}$ for sufficiently large $n$, and then the procedure would have terminated at stage $n$.
%\end{proof}				%%  Done
\chapter{Logical prerequisites}

Before we delve into reverse mathematics, it is necessary to have some background in computability and model theory, insofar as they apply to our setting of second-order arithmetic. First, we review the basic concepts of computability, including computable functions, c.e. sets and the universal function. We then develop the model theory of second-order arithmetic, and define the arithmetical and analytical hierarchies.

\section{Computability}
\label{sec:comp}

%%  Good informal account of computability - http://math.andrej.com/wp-content/uploads/2006/05/kleene-tree.pdf

Reverse mathematics is best understood with a background in computability theory. This is because we generally work over the base system $\RCA$, which can be thought of as the ``computable world''. A reversal of a statement $\varphi$ in a formal system $\Ss$ is then a proof of $\Ss$ from $\RCA + \varphi$. In practice, this involves a computable reduction between $\Ss$ and $\varphi$, hence the importance of computability.

The key notion in computability is that of an \textit{algorithm}, an exact method by which something can be computed. Algorithms can be formalised in many ways---Turing machines, the $\lambda$-calculus, $\mu$-recursive functions. All of these formalisations are provably equivalent, and the widely-accepted \textit{Church--Turing thesis} posits that each faithfully captures the idea of something being calculable. Therefore, to avoid unnecessary formality,\footnote{Rigorous definitions of computability are available in any introductory textbook on the subject. See \cite{rogersTheoryRecursiveFunctions1967, soareRecursivelyEnumerableSets1987, soareTuringComputability2016}.} the following intuitive "definition" is sufficient for us:

\begin{qdefinition}\label{defn:alg}
	An \textit{elementary instruction} is one that can be performed mechanically. An \textit{algorithm} is a finite collection of unambiguous elementary instructions, to be carried out in a specified order. Instructions may be repeated.

%	procedure that gives an answer after a finite number of stages. Algorithms may take inputs and return outputs.
	% Need a satisfactory definition of algorithm
\end{qdefinition}

We will allow our algorithms to take natural numbers as input, and act on this input during the computation. Furthermore, we will expect our algorithm to produce a natural number as output if the computation terminates, or \textit{halts}.

``Definition'' \ref{defn:alg} is intentionally very broad---almost all processes arising in mathematics and elsewhere qualify as algorithms. However, there are processes which don't. The archetypal example in computability is the \textit{halting problem}---determining whether a given program halts on a given input. A non-algorithmic process from classical mathematics arises in the proof of the Bolzano--Weierstrass theorem. Given a sequence $(x_n)_{n=1}^\infty$ of real numbers bounded in the interval $I$, we split $I$ into halves, take a half $I'$ which contains infinitely many of the $x_n$, and repeat the splitting on $I'$, ad infinitum.

The problem is that the instruction ``take a half $I'$ containing infinitely many $x_n$'' is \textit{not} elementary, in the sense that it is not possible to carry it out given an arbitrary sequence $(x_n)$. We could start counting along the sequence, noting which half $x_1$ is in, then $x_2$, and so on, but we will never know which one contains \textit{infinitely} many $x_n$. Indeed, a computable version of the Bolzano--Weierstrass theorem fails for this reason \cite{speckerNichtKonstruktivBeweisbare1949}.

\begin{figure}
	\begin{minipage}[b]{0.35\linewidth}
		\begin{algorithm}\label{alg:nohalt}\ % \vspace*{-1.5ex}
			\begin{enumerate}[label={\texttt{\arabic*}},labelsep=2em,labelwidth=5em,itemsep=0ex,parsep=0ex,leftmargin=4em]%,topsep=0ex,partopsep=0ex]
				\item \texttt{let }$x \defeq 0$
				\item \texttt{increment }$x$\texttt{ by }$1$
				\item \texttt{goto line }$2$
				\item \texttt{halt}
			\end{enumerate}
		\end{algorithm}
		\caption{An algorithm which never halts.}\label{fig:nohalt}
	\end{minipage}
	\hfill
	\begin{minipage}[b]{0.6\linewidth}
		\centering
		\begin{tikzpicture}[every path/.style={-stealth,gray},x=15mm,y=10mm]
		\foreach \x in {0,...,4}{
			\pgfmathsetmacro\ymax{5-\x}
			\foreach \y in {0,...,\ymax} 
			{\pgfmathtruncatemacro{\label}{(\x+\y)*(\x+\y+1)/2 + \x}
				\node[black] (\label) at (\y,-\x) {\label};} }
		
		\foreach \x in {2,...,5}{
			\pgfmathsetmacro\ymin{6-\x}
			\foreach \y in {\ymin,...,4}{
				\node at (\x,-\y) {$\cdots$};}}
		
		%%\foreach \y in {1,...,5}{
		%%	\pgfmathsetmacro\xmax{5-\y}
		%%	\foreach \x in {0,...,\xmax} 
		%%	{\draw (\x,-\y) -- ({\x+1},{1-\y});}}
		
		\foreach \n in {0,1,3,4,6,7,8,10,11,12,13,15,16,17,18}{
			\pgfmathtruncatemacro\nt{\n+1}
			\draw%[shorten <=-1mm,shorten >=-1mm]
			(\n) -- (\nt);}
		
		\draw[rounded corners,shorten <=2.5mm,shorten >=2.5mm] (0,-1) -- (0.5,-1) -- (1.5,0) -- (2,0);
		\draw[rounded corners,shorten <=2.5mm,shorten >=2.5mm] (0,-2) -- (0.5,-2) -- (2.5,0) -- (3,0);
		\draw[rounded corners,shorten <=2.5mm,shorten >=3.5mm] (0,-3) -- (0.5,-3) -- (3.5,0) -- (4,0);
		\draw[rounded corners,shorten <=3.5mm,shorten >=3.5mm] (0,-4) -- (0.5,-4) -- (4.5,0) -- (5,0);
		
		%%\pgfmathsetmacro\r{0.5*(1+sqrt(2))}
		%%\draw[->,shorten <=2mm,shorten >=2mm] (1,0) arc (0:-45:\r) -- ({sqrt(2)/4},{-1.5+sqrt(2)/4}) arc (135:180:\r);
		%%\draw[->,shorten <=2mm,shorten >=2mm] (2,0) arc (0:-45:\r) -- ({sqrt(2)/4},{-2.5+sqrt(2)/4}) arc (135:180:\r);
		%%\draw[->,shorten <=2mm,shorten >=2mm] (3,0) arc (0:-45:\r) -- ({sqrt(2)/4},{-3.5+sqrt(2)/4}) arc (135:180:\r);
		%%\draw[->,shorten <=2mm,shorten >=2mm] (4,0) arc (0:-45:\r) -- ({sqrt(2)/4},{-4.5+sqrt(2)/4}) arc (135:180:\r);
		
		\end{tikzpicture}
		\caption{Cantor's pairing function.}
	\end{minipage}
\end{figure}

There is a key difference between ``Definition'' \ref{defn:alg} and our intuitive understanding of an algorithm---for us, algorithms \textit{do not} have to halt. Algorithm \ref{alg:nohalt} in Figure \ref{fig:nohalt} is a simple example. This algorithm indeed satisfies ``Definition'' \ref{defn:alg}---each instruction is elementary and unambiguous, and the order in which they are to be executed (1, 2, 3, 2, 3, 2, 3, ...) is hopefully evident. However, Algorithm \ref{alg:nohalt} never halts, as it will continually increment the variable $x$, never reaching line 4. Notice that while algorithms may not halt, if they do, this must happen in finite time.

To model this idea, we instead consider \textit{partial functions} $f\colon \N \to \N$, i.e.\ functions $f\colon A \to \N$ for some subset $A \subseteq \N$, called the \textit{domain} of $f$.

%\begin{definition}%\label{defn:comp-func}
%	A function $f\colon \N \to \N$ is \textit{computable} if there exists an algorithm which, on input $n$, halts and outputs $f(n)$.
%\end{definition}

\begin{definition}\label{defn:comp-func}
	A partial function $f\colon \N \to \N$ is \textit{computable} if there exists an algorithm which, on input $n \in \N$:
	\begin{enumerate}
		\item If $n \in \dom(f)$: halts and outputs $f(n)$;
		\item If $n \notin \dom(f)$: doesn't halt.
	\end{enumerate}
\end{definition}

\begin{definition}\label{defn:comp-set}
	A set $A \subseteq \N$ is \textit{computable} if its characteristic function $\chr_A\colon \N \to \{ 0, 1 \}$ is computable, in the sense of Definition \ref{defn:comp-func}. Concretely, $A \subseteq \N$ is computable if there exists an algorithm which, given input $n$, always halts, returning 1 if $n \in A$, and 0 if $n \notin A$.
\end{definition}

Informally, a set is computable if there is an algorithm which tells us whether or not any given element is in the set. We can extend Definition \ref{defn:comp-func} to functions $f\colon \N^d \to \N^k$, and Definition \ref{defn:comp-set} to sets $A \subseteq \N^d$, via the \textit{pairing function}:

\begin{definition}[\cite{cantorBeitragZurMannigfaltigkeitslehre1877}]\label{defn:pair-fn}
	The \textit{2-pairing function} $\pi_2\colon \N^2 \to \N$ is the bijection $$\pi_2(m,n) = \dfrac{(m+n)(m+n+1)}{2} + m.$$ Further, we define the \textit{3-pairing function} $\pi_3\colon \N^3 \to \N$, $(m,n,k) \mapsto \pi_2 \big( m,\pi_2(n,k) \big)$, the \textit{4-pairing function} $\pi_4\colon \N^4 \to \N$, $(m,n,k,\ell) \mapsto \pi_2 \big( m,\pi_3(n,k,\ell) \big)$, etc. These are all bijections.
\end{definition}

The pairing functions allow us to treat tuples of natural numbers as single natural numbers, and therefore define computability for tuples of natural numbers. Often, we will implicitly use the pairing functions to think of a $d$-tuple $\bar{x}$ as just a natural number.

\begin{definition}%\label{defn:comp-set}
	A set $B \subseteq \N^d$ is \textit{computable} if $A = \pi_d(B)$ is computable, in the sense of Definition \ref{defn:comp-set}.
\end{definition}

\begin{definition}%\label{defn:comp-func}
	A partial function $f\colon \N^d \to \N^k$ is \textit{computable} if $\tilde{f} = \pi_k^{-1} \circ f \circ \pi_d$ is computable, in the sense of Definition \ref{defn:comp-func}.
\end{definition}

%{\color{red} c.e. = generalisation of computable}

%An algorithm which can always decide if $n \in A$

In computability, it is often useful to consider sets which are almost computable, but not quite. The computably enumerable sets provide examples of such things.

\begin{definition}
	Given a set $B \subseteq \N^2$ of pairs, the \textit{projection of $B$} is the set $$\proj(B) = \{ m \in \N: \exists n\ (m,n) \in B \}.$$ A set $A \subseteq \N$ is \textit{computably enumerable} (c.e.) if there is a computable set $B \subseteq \N^2$ such that $A = \proj(B)$.
\end{definition}

%Informally, $A$ is c.e. if there is an algorithm which, on input $n$, halts when $n \in A$, and doesn't halt if $n \notin A$. {\color{red} Can decide if it's in, but not if it's out - prove equivalence}

%%Given a semi-deciding algorithm

Every computable set $A$ is c.e., since it is the projection of $A \times \{ 0 \} = \{ (n,0): n \in A \}$. The converse does not hold, and there are many noncomputable c.e.\ sets; the archetypal example is the \textit{halting problem}.

\begin{proposition}\label{prop:comp-ce}
	$A \subseteq \N$ is computable if and only if both $A$ and $A^\complement$ are c.e..
\end{proposition}

\begin{proof}
	In the forward direction, $A^\complement$ is also computable, and we have noted that every computable set is c.e.. Conversely, suppose $A = \proj(B)$, $A^\complement = \proj(C)$ for computable $B, C \subseteq \N^2$. Given $n \in \N$, we decide if $n \in A$ as follows: first check if $(n,0) \in B$, then check if $(n,0) \in C$, then if $(n,1) \in B$, then if $(n,1) \in C$, and so on. Since it is true that either $n \in A$ or $n \in A^\complement$, eventually this algorithm will halt.
\end{proof}

We close this section with a fundamental result of computability, originally due to Turing \cite{turingComputableNumbersApplication1937}. By definition, every algorithm admits a finite description. Therefore, we can code algorithms by natural numbers, using a suitable coding scheme. For example, we could code each algorithm in a fixed programming language, and interpret its ASCII code as a natural number written in binary. For a nice coding scheme such as this, we can computably decode these numbers back into functions, and thus compute a function from its code. More formally:

\begin{theorem}[{\cite{turingComputableNumbersApplication1937}}]\label{prop:utm}
	There is a partial computable function $U\colon \N^2 \to \N$ with the following property: for any partial computable function $f\colon \N \to \N$, there is $e \in \N$ such that $U(e,n) = f(n)$ for all $n$.
\end{theorem}

\begin{proof}
	Compute $U$ as follows: given input $(e,n)$, interpret $e$ as a code for a computable function $f$, decode it, and compute $f(n)$.
\end{proof}

We call $U$ a \textit{universal computable function}. Essentially, $U$ can be interpreted as a compiler: it takes in the code $e$ of a function and returns the function itself. A corollary of Proposition \ref{prop:utm} is that the sequence $\varphi_0, \varphi_1, \varphi_2, \ldots$, where $\varphi_e(n) = U(e,n)$, lists all the partial computable functions. Furthermore, this is a \textit{uniformly computable listing}, meaning there is an algorithm taking $(e,n)$ to $\varphi_e(n)$ (namely, the algorithm for $U$). The existence of such a sequence will be useful later.

\section{Second-order arithmetic}
\label{sec:soa}

Now, we develop the necessary model-theoretic tools within the setting of arithmetic. The reader may have heard of first-order arithmetic, more commonly known as Peano arithmetic ($\PA$). The reason $\PA$ is \textit{first-order} is that quantification is only allowed over natural numbers. For example, a number $p$ being prime is expressible in $\PA$ (for all \textit{natural numbers} $m < p$, $m$ divides $p$ iff $m=1$ or $m=p$), but not the well-foundedness of $\N$ (for every \textit{subset} $A \subseteq \N$, $A$ has a least element).

Here, we work in the stronger setting of \textit{second-order arithmetic}, where quantification over subsets is allowed. We review basic model theory in this setting \cite{changModelTheory1990, markerModelTheoryIntroduction2002, simpsonSubsystemsSecondOrder2009}; in short, we consider structures in the \textit{language of second-order arithmetic} $\Ls = \{ 0, 1, +, \cdot, <, \in \}$. This is a two-sorted language, meaning we have two kinds of objects: numbers (denoted by lowercase letters $n,m,k,\ldots$), and sets (denoted in uppercase $A,B,C,\ldots$). The symbols in $\Ls$ are \textit{typed}, e.g. $0$ is a constant symbol of number type, $+$ is a binary operation between two object of number type, $\in$ is a binary relation between an object of number type and one of set type, etc.

We can build terms from symbols in $\Ls$, and we have two kinds of terms: numerical terms and set terms. As is usual, terms may include variable symbols, of number type $x$ or set type $X$. In fact, the only terms of set type are the set variable symbols $X,Y,Z,\ldots$, but there are a wealth of numerical terms:

\begin{definition}
	The collection of \textit{numerical $\Ls$-terms} is defined as follows:
	\begin{enumerate}
		\item 0, 1, and any numerical variable symbol $x$ are numerical terms.
		\item If $s$, $t$ are numerical terms, then $(s+t)$ and $(s \cdot t)$ are numerical terms.
	\end{enumerate}
\end{definition}

\noindent Intuitively, a numerical term represents a natural number. For example, $(1+1)$, $(1+(0+1))$ and $((1 \cdot (1 + 1)) + 0)$ are all numerical terms, all representing the number 2. However, these are all different \textit{terms}, since they do not contain the same arrangement of symbols. Frequently, we will omit brackets where there is no ambiguity---the above terms might be written more concisely as $1+1$, $1+0+1$ and $1 \cdot (1 + 1) + 0$. We will use $k$ to abbreviate the numerical term $\underbrace{1+1+\cdots+1}_{k\text{ times}}$.

\begin{definition}\label{defn:l2-form}
	The collection of \textit{$\Ls$-formulae} is defined as follows:
	\begin{enumerate}
		\item If $s$, $t$ are numerical terms, and $X$ is a set variable symbol, then $(s=t)$, $(s < t)$ and $(s \in X)$ are formulae.
		\item If $\varphi$, $\psi$ are formulae, then $(\lnot \varphi)$, $(\varphi \land \psi)$, $(\varphi \lor \psi)$, $(\varphi \to \psi)$ and $(\varphi \leftrightarrow \psi)$ are formulae.
		\item If $\varphi$ is a formula, then $(\forall x\ \varphi)$ and $(\exists x\ \varphi)$ are formulae.
		\item If $\varphi$ is a formula, then $(\forall X\ \varphi)$ and $(\exists X\ \varphi)$ are formulae.
	\end{enumerate}
\end{definition}

\noindent Intuitively, formulae are statements that may be true or false in a particular situation. Again, unnecessary brackets will often be omitted. We distinguish two types of variables in formulae: \textit{bound variables}, which are preceded by a quantifier over that variable, and \textit{free variables}, which are not. For example, in the formula $\forall x\ (x+y=1)$, the variable $x$ is bound by the quantifier $\forall x$, while $y$ is free.

\begin{definition}
	An \textit{$\Ls$-sentence} is an $\Ls$-formula in which all variables are bound.
\end{definition}

If $\varphi$ contained any free variables, then the truth or falsity of $\varphi$ could conceivably depend on what values were assigned to those free variables. Thus, an $\Ls$-sentence is a formula that can be assigned an unconditional truth value.

\begin{definition}
	Suppose $\varphi(x_1,\ldots,x_n,X_1,\ldots,X_m)$ is an $\Ls$-formula in free number variables $x_1,\ldots,x_n$, and free set variables $X_1,\ldots,X_m$. Then, the \textit{universal closure} of $\varphi$ is the $\Ls$-sentence $$\forall x_1 \cdots \forall x_n\ \forall X_1 \cdots \forall X_m\ \ \varphi(x_1,\ldots,x_n,X_1,\ldots,X_m)$$
\end{definition}

Now, what does it mean for an $\Ls$-sentence $\varphi$ to be true, or false? As is usual in model theory, truth of $\varphi$ is defined relative to a \textit{model}, consisting of a universe $\M$ of elements and interpretations in $\M$ for all symbols in our language. Since we are working with two sorts (numbers $x$ and sets $X$), we need to provide \textit{both} a universe of \textit{numbers} $\A$ and a universe of \textit{sets} $\B$, and interpret the symbols in $\Ll_2$ appropriately. For example, we would interpret $\in$ as a relation between elements of $\A$ and elements of $\B$.

In theory, we could pick any sets $\A$, $\B$ to serve as the universe for our model, and interpret the symbols in $\Ls$ any way we like. However, we will only be interested in the so-called \textit{$\omega$-models}, where $\A = \N$ is the natural numbers, and the symbols $0, 1, +, \cdot, <$ are given their usual interpretations in $\N$.

\begin{definition}\label{defn:om-mod}
	An \textit{$\omega$-model of second-order arithmetic} is a subset $\B \subseteq \Pow(\N)$.
\end{definition}

A priori, there is no reason that our universe of sets $\B$ actually must consist of subsets of $\N$. We could theoretically pick \textit{any} universe $\B$, and interpret the relation $x \in X$ in \textit{any} way we like. However, we can always identify $\B$ with a subset of $\Pow(\N)$ by identifying each $X \in \B$ with the set $\overline{X} = \{ n \in \N:$ the formula ``$n \in X$'' holds$\}$. So, no generality is lost in only considering subsets of $\Pow(\N)$ in Definition \ref{defn:om-mod}.

\begin{definition}\label{defn:truth}
	Given an $\omega$-model $\M$, \textit{truth of a sentence $\varphi$ in $\M$} (notated $\M \vDash \varphi$) is defined in the evident way:
	\begin{itemize}%[itemsep=0pt]
		\item Numerical terms $s$, $t$ are given their standard interpretations $s^\M$ and $t^\M$ in $\N$;
		\item $\M \vDash (s=t)$ if $s^\M$ and $t^\M$ are the same natural number;
		\item $\M \vDash (s<t)$ if $s^\M$ is a smaller natural number than $t^\M$;
		\item The rules for Boolean connectives $\lnot$, $\land$, $\lor$, $\to$, $\leftrightarrow$ are as usual;
		\item $\M \vDash \big( \exists x\ \varphi(x) \big)$ if there is some $n \in \N$ such that $\M \vDash \varphi(n)$;
		\item $\M \vDash \big( \forall x\ \varphi(x) \big)$ if $\M \vDash \varphi(n)$ for any choice of $n \in \N$;
		\item $\M \vDash \big( \exists X\ \varphi(X) \big)$ if there is some set $A \in \M$ such that $\M \vDash \varphi(A)$;
		\item $\M \vDash \big( \forall X\ \varphi(X) \big)$ if $\M \vDash \varphi(A)$ for any set $A \in \M$.
	\end{itemize}
\end{definition}

\noindent The key part of Definition \ref{defn:truth} is that set quantifiers $\forall X$ and $\exists X$ should be interpreted as ranging over exactly the sets \textit{in the model} $\M$. This is the key difference between the different $\omega$-models. As an example, the sentence $\exists X\ (0=0)$ is false in the $\omega$-model $\varnothing$, but is true in any other $\omega$-model.

%Given an $\omega$-model $\M$, we can now talk about what it means for a formula $\varphi$ to be true in $\M$. First, we need to define interpretations of numerical $\Ls$-terms in $\M$.
%
%\begin{definition}
%	Let $t(x_1,\ldots,x_m)$ be a numerical $\Ls$-term, such that every variable symbol in $t$ is one of $x_1,\ldots,x_m$. Let $\M$ be an $\omega$-model, and $n_1, \ldots, n_m \in \N$ be natural numbers. The \textit{interpretation} $t^\M[n_1,\ldots,n_m]$ of $t$ in $\M$ is defined as follows:
%	\begin{enumerate}
%		\item $0^\M$ is the natural number $0 \in \N$, and $1^\M$ is the natural number $1 \in \N$.
%		\item For any variable symbol $x$, $x^\M[n]$ is the natural number $n \in \N$.
%		\item For terms $s$, $t$ with interpretations $s^\M[n_1,\ldots,n_m] = k$, $t^\M[n_1,\ldots,n_m] = \ell$ respectively, we interpret $(s+t)^\M[n_1,\ldots,n_m] = k + \ell$ and $(s \cdot t)^\M[n_1,\ldots,n_m] = k \cdot \ell$.
%	\end{enumerate}
%\end{definition}
%
%\begin{definition}
%	Let $\varphi(x_1,\ldots,x_m,X_1,\ldots,X_d)$ be an $\Ls$-formula, with free number variables ranging in $x_1,\ldots,x_m$, and free set variables ranging in $X_1,\ldots,X_d$. Let $\M$ be an $\omega$-model, $n_1, \ldots, n_m \in \N$ be natural numbers, and $A_1,\ldots,A_d \in \M$ be sets. Inductively, define $\M \vDash \varphi$ (said ``$\varphi$ is true in $\M$'') as follows:
%	\begin{enumerate}
%		\item If $s$, $t$ are 
%	\end{enumerate}
%\end{definition}

\section{The arithmetical and analytical hierarchies}

The collection $\F$ of all $\Ls$-formulae, as in Definition \ref{defn:l2-form}, is an extremely rich and varied class. We wish to stratify $\F$ based on the complexity of formulae it contains. Our chosen measure of complexity will be based on the quantifiers, their type (numerical or set), and the number of alternations between universal ($\forall$) and existential ($\exists$). This way, we classify $\F$ into structures known as the \textit{arithmetical hierarchy} and the \textit{analytical hierarchy}.% Universal formula are given $\Pi$ classifications, and existential formula given $\Sigma$ classifications.

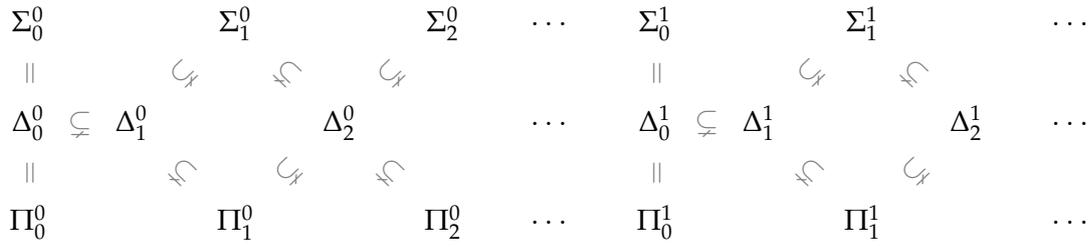
\begin{figure}%[b]%[t]
	\centering
	\begin{tikzpicture}
	\matrix(ah)[matrix of math nodes,row sep=6mm,column sep=6.5mm]{
		\Sigma^0_0 & & \Sigma^0_1 & & \Sigma^0_2 & \cdots & \Sigma^1_0 & & \Sigma^1_1 & & \cdots & \\
		\Delta^0_0 & \Delta^0_1 & & \Delta^0_2 & & \cdots & \Delta^1_0 & \Delta^1_1 & & \Delta^1_2 & \cdots & \\
		\Pi^0_0 & & \Pi^0_1 & & \Pi^0_2 & \cdots & \Pi^1_0 & & \Pi^1_1 & & \cdots \\
	};
	\begin{scope}[every node/.style={sloped,gray}]
	\path (ah-1-1) -- node{$=$} (ah-2-1) -- node{$=$} (ah-3-1);
	\path (ah-2-1) -- node{$\subsetneq$} (ah-2-2) -- node{$\subsetneq$} (ah-1-3) -- node{$\subsetneq$} (ah-2-4) -- node{$\subsetneq$} (ah-1-5);
	\path (ah-2-2) -- node{$\subsetneq$} (ah-3-3) -- node{$\subsetneq$} (ah-2-4) -- node{$\subsetneq$} (ah-3-5);
	\path (ah-1-7) -- node{$=$} (ah-2-7) -- node{$=$} (ah-3-7);
	\path (ah-2-7) -- node{$\subsetneq$} (ah-2-8) -- node{$\subsetneq$} (ah-1-9) -- node{$\subsetneq$} (ah-2-10);% -- node{$\subsetneq$} (ah-1-11);
	\path (ah-2-8) -- node{$\subsetneq$} (ah-3-9) -- node{$\subsetneq$} (ah-2-10);% -- node{$\subsetneq$} (ah-3-11);
	\end{scope}
	\end{tikzpicture}
	\caption[The arithmetical and analytical hierarchies.]{The arithmetical (left) and analytical (right) hierarchies.}
\end{figure}

The lowest level of complexity consists of formulae containing only bounded quantifiers: those of the form $\forall x\ (x < k\, \to\, \psi)$ or $\exists x\ (x < k\, \to\, \psi)$ for some constant $k \in \N$. We will often abbreviate these to $(\forall x < k)\, \psi$ and $(\exists x < k)\, \psi$ respectively. From there, universal formulae are given $\Pi$ classifications, and existential formulae given $\Sigma$ classifications.

\begin{definition}[(arithmetical hierarchy for formulae)]
	Let $\varphi$ be an $\Ls$-formula. We assign classifications to $\varphi$ as follows:
	\begin{enumerate}
		\item $\varphi$ is called $\Sigma^0_0$ and $\Pi^0_0$ if it only contains bounded quantifiers.
		\item $\varphi$ is called $\Sigma^0_{n+1}$ if it is of the form $\varphi = \exists x_1 \cdots \exists x_n\ \psi$, where $\psi$ is $\Pi^0_n$.
		\item $\varphi$ is called $\Pi^0_{n+1}$ if it is of the form $\varphi = \forall x_1 \cdots \forall x_n\ \psi$, where $\psi$ is $\Sigma^0_n$.
	\end{enumerate}
	We say $\varphi$ is \textit{arithmetical} if it receives any of these classifications.
\end{definition}

We also translate the arithmetical hierarchy from formulae to sets defined by those formulae. This gives us a measure of complexity for subsets of $\N$. Here, we obtain additional $\Delta$ classifications.

\begin{definition}[(arithmetical hierarchy for sets)]\
	\begin{enumerate}
		\item A set $A \subseteq \N$ is called $\Sigma^0_n$ if there is a $\Sigma^0_n$ formula $\varphi(x)$ in one free variable such that $A = \{ n \in \N: \varphi(n) \text{ holds} \}$. $\Pi^0_n$ sets are defined analogously.
		\item $A \subseteq \N$ is called $\Delta^0_n$ if it is both $\Sigma^0_n$ and $\Pi^0_n$.
	\end{enumerate}
	We say $A$ is \textit{arithmetical} if it receives any of these classifications, or equivalently, if it is definable by an arithmetical formula.
\end{definition}

For sets, we have that $\Sigma^0_n \subseteq \Pi^0_{n+1}$. If $\varphi(x)$ is a $\Sigma^0_n$ formula defining $A$, and $y$ is a variable not in $\varphi$, then $\psi(x) = \forall y\ \varphi(x)$ is a $\Pi^0_{n+1}$ formula also defining $A$ (since the truth value of $\varphi$ does not depend on $y$). We also get that $\Sigma^0_n \subseteq \Sigma^0_{n+1}$ by placing such ``dummy quantifiers'' after all others, whence $\Sigma^0_n \subseteq \Delta^0_{n+1}$. By taking complements, $\Pi^0_n \subseteq \Delta^0_{n+1}$, thus $\Sigma^0_n \cup \Pi^0_n \subseteq \Delta^0_{n+1}$. In fact, all of these containments are strict, but we will not prove this here.

There is a close relationship between the arithmetical hierarchy and computability:

\begin{proposition}\label{prop:d00-comp}
	Every $\Delta^0_0$ set is computable.
\end{proposition}

\begin{proof}
	To say $A \subseteq \N$ is $\Delta^0_0$ is to say that there is an $\Ls$-formula $\varphi(x)$ such that $n \in A \iff \varphi(n)$ holds, and where all quantifiers in $\varphi(x)$ are bounded. Since $\varphi(x)$ can only contain finitely many quantifiers (say $d$-many), there are only finitely many tuples $(a_1,\ldots,a_d) \in \N^d$ that we need to check to verify whether $\varphi(n)$ holds or not. So, the algorithm to compute $A$ is simply checking all such tuples exhaustively.
\end{proof}

\begin{proposition}[\cite{kleeneRecursivePredicatesQuantifiers1943, postRecursivelyEnumerableSets1944}]\label{prop:s01-ce}
	$A \subseteq \N$ is $\Sigma^0_1$ if and only if it is computably enumerable.
\end{proposition}

\begin{proof}
	We prove the forward direction. If $A$ is $\Sigma^0_1$, then it can be defined by a formula $\varphi(x)$ of the form $\exists y_1 \cdots \exists y_n\ \psi(x,y_1,\ldots,y_n)$, where $\psi$ is $\Delta^0_0$. By Proposition \ref{prop:d00-comp}, the set $B \defeq \{ (x,\bar{y}): \psi(x,\bar{y}) \}$ is computable, and $A = \proj(B)$.
%	{\color{red} c.e. $\implies \Sigma^0_1$ is hard}
\end{proof}

The reverse direction of Proposition \ref{prop:s01-ce} is harder to prove. It requires coding algorithms using $\Ls$-formulae, for which a formal definition of computability is needed. Therefore, we will not complete the proof here, but it can be found in \cite{davisComputabilityUnsolvability1958, rogersTheoryRecursiveFunctions1967}.

\begin{corollary}\label{cor:d01-comp}
	$A \subseteq \N$ is $\Delta^0_1$ if and only if it is computable.
\end{corollary}

\begin{proof}
	$A$ is $\Delta^0_1$ if and only if $A$ is both $\Sigma^0_1$ and $\Pi^0_1$. By the negation rules for $\forall$ and $\exists$ quantifiers, $A$ is $\Pi^0_1$ if and only if $A^\complement$ is $\Sigma^0_1$. By Proposition \ref{prop:s01-ce}, this is if and only if $A$ and $A^\complement$ are c.e.. Hence, the result follows from Proposition \ref{prop:comp-ce}.
\end{proof}

Clearly, there are nonarithmetical formulae---any formula containing a set quantifier is an example. It is less obvious that there are also nonarithmetical sets, but this follows from a simple counting argument. Every arithmetical set is defined by an $\Ls$-formula, of which there are countably many, while there are continuum-many subsets of $\N$.

We can extend the arithmetical hierarchy to the \textit{analytical hierarchy} in much the same manner:

\begin{definition}[(analytical hierarchy for formulae)]\
	\begin{enumerate}
		\item An $\Ls$-formula $\varphi$ is called $\Sigma^1_0$ and $\Pi^1_0$ if it is arithmetical.
		\item $\varphi$ is called $\Sigma^1_{n+1}$ if it is of the form $\varphi = \exists X_1 \cdots \exists X_n\ \psi$, where $\psi$ is $\Pi^1_n$.
		\item $\varphi$ is called $\Pi^1_{n+1}$ if it is of the form $\varphi = \forall X_1 \cdots \forall X_n\ \psi$, where $\psi$ is $\Sigma^1_n$.
	\end{enumerate}
\end{definition}

\noindent So, the analytical formulae are those which allow some level of quantification over sets. $\Sigma^1_n$, $\Pi^1_n$ and $\Delta^1_n$ sets of natural numbers are defined in exactly the same way. As before, we have that $\Sigma^1_n \cup \Pi^1_n \subseteq \Delta^1_{n+1}$ for all $n \in \N$.			%%  Done
\chapter{Subsystems of second-order arithmetic}
\label{chp:soa}

With the tools of computability and model theory in hand, we can now develop the formalism of reverse mathematics. Here, we will define the subsystems $\RCA$, $\WKL$, $\ACA$, $\ATR$ and $\PiCA$ and their $\omega$-models, and see where famous theorems of mathematics show up in this hierarchy. Most of the material of this chapter can be found in \cite{simpsonSubsystemsSecondOrder2009}.

\section{Formal systems}

\begin{definition}
	A \textit{formal system} or \textit{subsystem of second-order arithmetic} is a collection $\Ss$ of $\Ls$-sentences. We refer to the formulae in $\Ss$ as \textit{axioms of $\Ss$.}
\end{definition}

\begin{definition}
	Let $\M$ be an $\omega$-model, and $\Ss$ be a subsystem of second-order arithmetic. We say \textit{$\M$ is a model of $\Ss$} if all the axioms of $\Ss$ are true in $\M$.
\end{definition}

There are infinitely many inequivalent subsystems of second-order arithmetic, but there are five major ones which show up consistently in reverse mathematics. In order of increasing logical strength, these systems are called $\RCA$, $\WKL$, $\ACA$, $\ATR$, and $\PiCA$. These subsystems are affectionately known as the ``Big Five''---their significance comes from the fact that almost all classical theorems turn out to be equivalent to one of the Big Five.

%The system $\ATR$ will not be relevant for us, so we avoid defining it here, but we will give definitions of the other four.
All these subsystems include the following set of basic axioms:

\begin{axioms}
	The \textit{basic axioms} of second-order arithmetic are the following $\Ls$-sentences:
%	\settasks{debug}
	\begin{tasks}[label=(\roman*),label-align=right,label-width=2.2em,column-sep=3em,label-offset=1em,item-indent=4em,after-item-skip=0ex,before-skip=-5pt](3)
		\task $\forall n\ \lnot(n+1 = 0)$
		\task* $\forall n\ \forall m\ \big[ (n+1=m+1) \to (n=m) \big]$
		\task $\forall n\ (n+0 = n)$
		\task* $\forall n\ \forall m\ \big[ n+(m+1) = (n+m)+1 \big]$
		\task $\forall n\ (n \cdot 0 = 0)$
		\task* $\forall n\ \forall m\ \big[ n \cdot (m+1) = (n \cdot m)+n \big]$
		\task $\forall n\ \lnot(n < 0)$
		\task* $\forall n\ \forall m\ \big[ (n < m+1) \leftrightarrow (n < m\ \lor\ n=m) \big]$
	\end{tasks}
%	\begin{multicols}{3}
%	\begin{enumerate}
%		\item $\forall n\ \lnot(n+1 = 0)$
%		\item $\forall n\ \forall m\ \big[ (n+1=m+1) \to n=m \big]$
%		\item $\forall n\ (n+0 = n)$
%%		\item $\forall n\ \forall m\ \big[ n+(m+1) = (n+m)+1 \big]$
%%		\item $\forall n\ (n \cdot 0 = 0)$
%%		\item $\forall n\ \forall m\ \big[ n \cdot (m+1) = (n \cdot m)+n \big]$
%%		\item $\forall n\ \lnot(n < 0)$
%%		\item $\forall n\ \forall m\ \big[ (n < m+1) \leftrightarrow (n < m\ \lor\ n=m) \big]$
%	\end{enumerate}
%	\end{multicols}
\end{axioms}

\noindent Note that the basic axioms are entirely first-order---there is no mention of sets. They are closely related to the Peano axioms $\PA$. The basic axioms formalise the essential properties of $\N$, and are sufficient to prove all basic facts of arithmetic---commutativity, associativity, distributivity, etc. It is easily verified that $\N$ satisfies the basic axioms, whence:

\begin{proposition}\label{prop:om-mod-basic-ax}
	Any $\omega$-model satisfies the basic axioms.
\end{proposition}

The key feature distinguishing the different subsystems is the second-order axioms they contain. Most of the additional axioms we consider will have one of two forms. The first type are \textit{induction axioms}, allowing us to induct over certain statements:

\begin{definition}
	Let $\varphi(n)$ be an $\Ls$-formula in which $n$ appears freely. The \textit{induction axiom for $\varphi$} is the universal closure of $\big[ \varphi(0)\, \land\, \forall n\ \big( \varphi(n) \to \varphi(n+1) \big) \big] \to \forall n\ \varphi(n)$.
\end{definition}

The induction axiom for $\varphi$ allows us to perform induction on $\varphi$. When defining subsystems $\Ss$ of second-order arithmetic, we will generally limit the inductive strength to some point in the arithmetical/analytical hierarchy. For example, $\Ss$ may include the induction axiom for all $\Pi^0_2$ formulae $\varphi$. We do this because we are trying to find the weakest subsystem in which a theorem $\varphi$ is provable; therefore, we don't allow induction beyond what is truly necessary.

As we know, induction in $\N$ is valid for \textit{any} $\Ls$-formula $\varphi$, whence:

\begin{proposition}\label{prop:om-mod-ind}
	Any $\omega$-model satisfies the induction axiom for any $\Ls$-formula $\varphi$.
\end{proposition}

The second type are \textit{comprehension axioms}, guaranteeing that given sets must exist:

\begin{definition}
	Let $\varphi(n)$ be an $\Ls$-formula in which $n$ appears freely, but $X$ does not appear. The \textit{comprehension axiom for $\varphi$} is the universal closure of $\exists X\ \forall n\ \big[ n \in X\, \leftrightarrow\, \varphi(n) \big]$.
\end{definition}

Essentially, the comprehension axiom for $\varphi$ asserts that the set $A_\varphi = \{ n \in \N: \varphi(n) \}$ exists. Again, our subsystems will generally include all comprehension axioms up to some point in the arithmetical/analytical hierarchy.

\section{\texorpdfstring{$\RCA$}{RCA-0}}

\begin{definition}
	$\RCA$ is the subsystem consisting of:
	\begin{enumerate}
		\item the basic axioms;
		\item the induction axiom for every $\Sigma^0_1$ formula $\varphi$;
		\item the \textit{$\Delta^0_1$ comprehension scheme}: the universal closure of $$\big[ \forall n\ \big( \varphi(n) \leftrightarrow \psi(n) \big) \big]\, \to\ \exists X\ \forall n\ \big( n \in X \leftrightarrow \varphi(n) \big)$$ for every $\Sigma^0_1$ formula $\varphi(x)$ and $\Pi^0_1$ formula $\varphi(x)$ containing $x$ as a free variable, but not containing $n$ or $X$.
	\end{enumerate}
\end{definition}

\noindent $\RCA$ stands for ``recursive comprehension axiom'', as it allows comprehension over $\Delta^0_1$ sets (which in $\N$ are the computable sets, as we saw in Proposition \ref{cor:d01-comp}). This is the weakest subsystem we will consider, and intuitively, it should be thought of as corresponding to computable mathematics. Generally, a statement holds in $\RCA$ if and only if a ``computable version'' of the statement is true. Some results from ordinary mathematics do hold in $\RCA$:

\begin{proposition}\label{prop:rca-prov}
	The following theorems are provable in $\RCA$:
	\begin{enumerate}
		\item The Baire category theorem \cite{simpsonSubsystemsSecondOrder2009};
		\item The intermediate value theorem \cite{pour-elComputabilityAnalysisPhysics1989};
		\item The soundness theorem for first order logic \cite{simpsonSubsystemsSecondOrder2009};
		\item Every countable, finite-rank matroid has a basis \cite{hirstReverseMathematicsMatroids2017};
		\item The Weierstrass approximation theorem \cite{pour-elSimpleDefinitionComputable1975}.
	\end{enumerate}
\end{proposition}

\noindent Proposition \ref{prop:rca-prov} boils down to the fact that each of these theorems is computably true. For example, the intermediate value theorem holds in $\RCA$ because the effective IVT is true: if $f\colon \UI \to \R$ is computable and $f(0) < 0 < f(1)$, then there is a computable real number $x \in \UI$ such that $f(x)=0$. But there are many more results that are not computably true, and thus don't hold in $\RCA$:

\begin{proposition}\label{prop:rca-nprov}
	The following theorems are \textit{not} provable in $\RCA$:
	\begin{enumerate}
		\item The Bolzano--Weierstrass theorem \cite{speckerNichtKonstruktivBeweisbare1949};
		\item The Heine--Borel theorem for countable covers \cite{friedmanSystemsSecondOrder1976};
		\item The extreme value theorem \cite{simpsonSubsystemsReverseMathematics1987, simpsonSubsystemsSecondOrder2009};
		\item G\"odel's completeness theorem for first order logic \cite{simpsonSubsystemsSecondOrder2009};
		\item Every continuous functions is Riemann integrable \cite{simpsonSubsystemsSecondOrder2009};
		\item Every countable vector space has a basis \cite{friedmanCountableAlgebraSet1983}.
	\end{enumerate}
\end{proposition}

\noindent Again, interpret Proposition \ref{prop:rca-nprov} as saying those theorems are not computably true. For example, the Bolzano--Weierstrass theorem does not hold in $\RCA$, because the effective Bolzano--Weierstrass theorem fails---there is a computable Cauchy sequence whose limit is not computable \cite{speckerNichtKonstruktivBeweisbare1949}.

The standard $\omega$-model of $\RCA$ is $\REC$, consisting of all recursive, or computable, subsets of $\N$. Given our intuition about $\RCA$, it should be no surprise that $\REC$ actually \textit{is} a model of $\RCA$:

\begin{proposition}\label{prop:rec-rca}
	$\REC$ is a model of $\RCA$.
\end{proposition}

\begin{proof}
	Since $\REC$ is an $\omega$-model, it satisfies the basic axioms by Proposition \ref{prop:om-mod-basic-ax}, and all induction axioms by Proposition \ref{prop:om-mod-ind}. Now, let $\varphi$ be a $\Sigma^0_1$ formula, and $\psi$ be $\Pi^0_1$. If $\varphi$ and $\psi$ are equivalent, then they both define the same set $A_\varphi = A_\psi \subseteq \N$. By Corollary \ref{cor:d01-comp}, $A_\varphi$ is computable, so $A_\varphi \in \REC$. Thus, $A_\varphi$ witnesses $\Delta^0_1$ comprehension for $\varphi$ and $\psi$.
\end{proof}

Recall our discussion of reverse mathematics from the introduction---a key idea was the \textit{reversal}, where, to show that $\Ss$ is the weakest system in which $\varphi$ can be proved, we demonstrate a proof of $\Ss$ from $\varphi$. No single theorem is strong enough to axiomatise mathematics; hence, when doing a reversal in practice, we need to supplement $\varphi$ with a weak base theory $\B$. It is customary to take $\B = \RCA$ (though weaker/stronger systems have been used at times).

%In our study of second-order arithmetic, we will take $\RCA$ as a \textit{base system}. This means that we will always assume $\RCA$ holds in the background. This is essential when proving reversals of a statement $\varphi$ {\color{red} where?}
%
%\begin{definition}
%	A reversal of a theorem $\varphi$ in a subsystem $\Ss$ is a proof of $\Ss$ from $\RCA + \varphi$.
%\end{definition}
%
%We will see two examples of reversals, {\color{red} where?}

\section{\texorpdfstring{$\WKL$}{WKL-0}}
\label{sec:wkl}

\textit{K\H onig's lemma} is a statement about infinite well-founded trees, and \textit{weak K\H onig's lemma} is the restriction of this to binary trees, i.e.\ those where each node has at most two children. It is convenient to define a tree as a certain subset of the following:

\begin{definition}
	\textit{(Finitary) Cantor space} $\CSf$ is the set of all finite binary sequences $\sigma = \sigma_0 \sigma_1 \cdots \sigma_{n-1}$, where each $\sigma_i \in \{ 0, 1 \}$. \textit{(Infinitary) Cantor space} $\CS$ consists of all infinite binary sequences $X = X_0 X_1 X_2 \cdots$, where each $X_i \in \{ 0, 1 \}$.
\end{definition}

$\CSf$ is countable, so we can represent its elements in $\N$ as follows: for $n \in \N$, write $n+1$ in binary and remove the leading 1. For example, 22 in binary is 10110, hence 21 represents the string $0110 \in \CSf$. This is a bijection between $\N$ and $\CSf$ (0 represents the empty string). Henceforth, when we talk about elements $\sigma \in \CSf$ in second-order arithmetic, they will be understood as natural numbers via this coding.

However, $\CS$ is uncountable, so it cannot be represented in $\N$. Instead, we need to represent $\CS$ using subsets of $\N$. The obvious way is to represent $X = X_0 X_1 X_2 \cdots$ by the set $A_X = \{ n \in \N: X_n = 1 \}$, so that $X$ is essentially the characteristic function of $A_X$. Then, given an $\omega$-model $\M$, an element $X \in \CS$ exists in $\M$ if and only if $A_X \in \M$.

\begin{definition}
	For an element $\sigma = \sigma_0 \sigma_1 \cdots \sigma_{n-1} \in \CSf$, its \textit{length} is $\abs{\sigma} = n$. Given $\sigma \in \CSf$ and $\tau \in (\CSf \cup \CS)$, we say \textit{$\tau$ extends $\sigma$} (notated $\sigma \preccurlyeq \tau$) if $\abs{\sigma} \leq \abs{\tau}$ and for all $k < \abs{\sigma}$, $\sigma_k = \tau_k$.
\end{definition}

\begin{definition}
	A tree is a subset $T \subseteq \CSf$ which is closed under initial segments. That is, if $\sigma = \sigma_0 \sigma_1 \cdots \sigma_{n-1} \in T$, then $\substr{\sigma}{k} \defeq \sigma_0 \sigma_1 \cdots \sigma_{k-1} \in T$ for any $k \leq n$.
\end{definition}

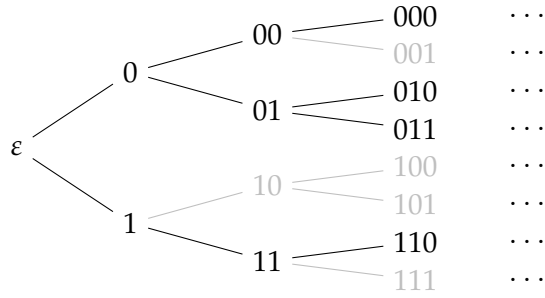
\begin{figure}
	\centering
	\begin{tikzpicture}[x=15mm]
		\def\xa{1}
		\def\xb{2.2}
		\def\xc{3.5}
		\def\xd{4.5}
		
		\node (rt) at (0,0) {$\varepsilon$};
		\node (0) at (\xa,1) {0};
		\node (1) at (\xa,-1) {1};
		\node (00) at (\xb,1.5) {00};
		\node (01) at (\xb,0.5) {01};
		\node[lightgray] (10) at (\xb,-0.5) {10};
		\node (11) at (\xb,-1.5) {11};
		
		\node (000) at (\xc,1.75) {000};
		\node[lightgray] (001) at (\xc,1.25) {001};
		\node (010) at (\xc,0.75) {010};
		\node (011) at (\xc,0.25) {011};
		\node[lightgray] (100) at (\xc,-0.25) {100};
		\node[lightgray] (101) at (\xc,-0.75) {101};
		\node (110) at (\xc,-1.25) {110};
		\node[lightgray] (111) at (\xc,-1.75) {111};
		
		\node at (\xd,1.75) {$\cdots$};
		\node at (\xd,1.25) {$\cdots$};
		\node at (\xd,0.75) {$\cdots$};
		\node at (\xd,0.25) {$\cdots$};
		\node at (\xd,-0.25) {$\cdots$};
		\node at (\xd,-0.75) {$\cdots$};
		\node at (\xd,-1.25) {$\cdots$};
		\node at (\xd,-1.75) {$\cdots$};
		
		\draw            (rt) -- (0) -- (00) -- (000);
		\draw[lightgray]                (00) -- (001);
		\draw                    (0) -- (01) -- (010);
		\draw                           (01) -- (011);
		\draw[lightgray]         (1) -- (10) -- (100);
		\draw[lightgray]                (10) -- (101);
		\draw            (rt) -- (1) -- (11) -- (110);
		\draw[lightgray]                (11) -- (111);
	\end{tikzpicture}
	\caption{A tree $T \subseteq \CSf$ as a subgraph of the full binary tree $\CSf$.}
	\label{fig:tree}
\end{figure}

As in Figure \ref{fig:tree}, we intuitively think of a tree $T \subseteq 2^\omega$ as a certain graph, where the vertices are the elements, and $\sigma \in T$ is connected to its direct extensions $\sigma \ct 0$ and $\sigma \ct 1$, if they are in $T$. Here, $\ct$ denotes concatenation of strings. Every nonempty tree contains the empty sequence $\varepsilon$; this is the ``root'' of the tree.

\begin{definition}
	Given a tree $T \subseteq \CSf$, a \textit{path through $T$} is an element $X \in \CS$ such that $\substr{X}{n} \defeq X_0 X_1 \cdots X_{n-1} \in T$ for any $n \in \N$.
\end{definition}

\begin{definition}
	\textit{Weak K\H onig's lemma} is the statement that every infinite tree $T \subseteq \CSf$ contains a path. $\WKL$ is the subsystem consisting of $\RCA$ plus weak K\H onig's lemma.
\end{definition}

By definition, $\WKL$ is at least as strong as $\RCA$ in terms of logical strength. We can show $\WKL$ is strictly stronger. The idea is to show the standard model $\REC$ of $\RCA$ is \textit{not} a model of $\WKL$, by constructing a computable tree $T \subseteq \CSf$ in $\REC$ with no computable path $X \in \REC$. This result is originally due to Jockusch and Soare \cite{jockuschDegreesMembersPi1972}, though our construction is different to theirs. We use a typical \textit{diagonalisation argument} ubiquitous in computability; we ensure at the $e$th step that $\varphi_e$ is not a branch through $T$.

\begin{proposition}[\cite{jockuschDegreesMembersPi1972}]
	There is a computable tree $T \subseteq \CSf$ with no computable path.
\end{proposition}
%The proof of Proposition \ref{prop:rec-wkl} essentially boils down to the fact that the computable version of weak K\H onig's lemma does not hold---there is an infinite computable tree $T$ with no computable path. This is a typical example of a \textit{diagonalisation argument} ubiquitous in computability: we ensure at the $e$th step that $\varphi_e$ is not a branch through $T$.

\begin{proof}
	Recall $\varphi_0, \varphi_1, \varphi_2, \ldots$ is a uniformly computable listing of all partial computable functions (Proposition \ref{prop:utm}). Construct a tree $T$ as follows. To test if $\sigma \in T$, for every $e < \abs{\sigma}$, run the computation of $\varphi_e(e)$ for $\abs{\sigma}$ steps. If any of these computations halt with $\varphi_e(e) = \sigma_e$, then $\sigma \notin T$; otherwise, $\sigma \in T$.
	
	$T$ is a tree: instead of showing $T$ is downwards closed, we (equivalently) show the complement is upwards closed. Suppose that $\sigma \notin T$---then, there is $e < \abs{\sigma}$ such that $\varphi_e(e)$ halts within $\abs{\sigma}$ steps, and $\varphi_e(e) = \sigma_e$. Then, for any extension $\tau \succcurlyeq \sigma$, $\varphi_e(e)$ also halts within $\abs{\tau} \geq \abs{\sigma}$ steps, and $\varphi_e(e) = \sigma_e = \tau_e$. Hence, $\tau \notin T$.
	
%	 {\color{red} show contrapositive}%since if $\sigma \in T$ with $\abs{\sigma}=n$, this means for every $e < n$, either $\varphi_e(e)$ doesn't halt after $n$ steps, or $\varphi_e(e)$ halts but $\neq\sigma_e$. Then for $k<n$,  $\substr{\sigma}{k} \in T$ for
	
	$T$ is infinite, since every level $n$ is nonempty. To see this, for each $e < n$ such that $\varphi_e(e)$ is defined, we can pick $\sigma_e \neq \varphi_e(e)$. If $\varphi_e(e)$ is undefined, just pick $\sigma_e$ arbitrarily. Then $\sigma \in T$ is on level $n$.
	
	$T$ is computable, since we gave an algorithm to compute it. Hence, $T$ exists in $\REC$. But we claim $T$ has no path in $\REC$. Suppose $T$ did have a path $X \in \REC$: then, $X = \varphi_e$ for some $e$, since the sequence $(\varphi_n)$ lists \textit{all} partial computable functions. Since $X$ is total, the computation of $\varphi_e(e) = X_e$ halts after, say, $s$ steps. But then $\substr{X}{s} \notin T$ by definition, so $X$ is not a path through $T$.
\end{proof}

Since we have constructed an infinite tree $T \in \REC$ with no path $X \in \REC$, it follows that weak K\H onig's lemma does not hold in $\REC$, whence:% $\REC$ is not a model of $\WKL$.

\begin{corollary}\label{prop:rec-wkl}
	$\REC$ is not a model of $\WKL$.
\end{corollary}

Weak K\H onig's lemma is closely related to the finite intersection characterisation of compactness; indeed, it can be viewed as asserting that Cantor space $\CS$, the infinite product of the discrete space $2$, is compact. Therefore, $\WKL$ is generally strong enough to perform compactness arguments. $\WKL$ can prove all the results of Proposition \ref{prop:rca-prov}, and it is equivalent to the following results (in other words, $\WKL$ is the weakest system in which they can be proved):

\begin{proposition}\label{prop:wkl-prov}
	Over $\RCA$, $\WKL$ is equivalent to:
	\begin{enumerate}
		\item The Heine--Borel theorem for countable covers \cite{friedmanSystemsSecondOrder1976};
		\item Continuous functions on $\UI$ are bounded \cite{simpsonSubsystemsReverseMathematics1987, simpsonSubsystemsSecondOrder2009};
		\item Continuous functions on $\UI$ are uniformly continuous \cite{simpsonSubsystemsReverseMathematics1987, simpsonSubsystemsSecondOrder2009};
		\item Continuous functions on $\UI$ are Riemann integrable \cite{simpsonSubsystemsSecondOrder2009};
		\item The extreme value theorem \cite{simpsonSubsystemsReverseMathematics1987, simpsonSubsystemsSecondOrder2009};
		\item G\"odel's completeness theorem for first order logic \cite{simpsonSubsystemsSecondOrder2009};
		\item Every countable commutative ring has a prime ideal \cite{friedmanCountableAlgebraSet1983};
		\item Brouwer's fixed point theorem \cite{shiojiFixedPointTheory1990};
		\item The Hahn--Banach theorem for separable Banach spaces \cite{brownWhichSetExistence1986}.
	\end{enumerate}
\end{proposition}%{\color{red} Give example of a reversal? Heine--Borel?}

\noindent That said, $\WKL$ is still insufficient to prove many important mathematical results, such as the completeness of $\R$, the Bolzano--Weierstrass theorem, and the existence of bases for vector spaces. We will see below that each of these statements is equivalent to the stronger system $\ACA$.

\section{\texorpdfstring{$\ACA$}{ACA-0}}

Now, we move on to stronger subsystems of second-order arithmetic.
% define two significantly stronger subsystems of second-order arithmetic.
The first is $\ACA$ (arithmetical comprehension), which guarantees the existence of any arithmetical set:

\begin{definition}
	$\ACA$ is the subsystem consisting of:
	\begin{enumerate}
		\item the basic axioms;
		\item the induction axiom for every arithmetical formula $\varphi$;
		\item the comprehension axiom for every arithmetical formula $\varphi$.
	\end{enumerate}
\end{definition}

\noindent $\ACA$ has a standard $\omega$-model $\ARITH$, consisting of all arithmetical subsets of $\N$.

\begin{proposition}
	$\ARITH$ is a model of $\ACA$.
\end{proposition}

\begin{proof}
	Very similar to Proposition \ref{prop:rec-rca}, so omitted.
\end{proof}

Clearly $\ACA$ implies $\RCA$, and therefore all of the results of Proposition \ref{prop:rca-prov}. It is less obvious that $\ACA$ implies $\WKL$, and thereby the results of Proposition \ref{prop:wkl-prov}. In fact, $\ACA$ is strong enough to prove almost all the results of classical mathematics (algebra, analysis, etc.), and virtually all the theorems taught in undergraduate mathematics.

\begin{proposition}\label{prop:aca-prov}
	Over $\RCA$, $\ACA$ is equivalent to:
	\begin{enumerate}
		\item The sequential completeness of the reals \cite{friedmanSystemsSecondOrder1976};
		\item The Bolzano--Weierstrass theorem \cite{friedmanSystemsSecondOrder1976};
		\item Every countable commutative ring has a maximal ideal \cite{friedmanCountableAlgebraSet1983};
		\item Every countable vector space has a basis \cite{friedmanCountableAlgebraSet1983};
		\item K\H onig's lemma: every infinite, finitely branching tree has an infinite path \cite{friedmanSystemsSecondOrder1974, friedmanSystemsSecondOrder1976};
		\item Ramsey's theorem for $k$-tuples, for fixed $k \geq 3$ \cite{jockuschRamseyTheoremRecursion1972, simpsonSubsystemsSecondOrder2009}.
	\end{enumerate}
\end{proposition}

\noindent There are a few mathematical theorems still out of reach for $\ACA$; for instance, in areas such as set theory, Ramsey theory and descriptive set theory, where strong set existence axioms are required. We will see some examples in the next two sections.

\section{\texorpdfstring{$\ATR$}{ATR-0}}

$\ATR$ (arithmetical transfinite recursion) comprises $\ACA$ plus the assertion that any ``arithmetical operator'' can be iterated along any countable ordinal, starting with any set. Let's try to understand what all of this means.

\begin{definition}%[($\ACA$)]
	A (countable) ordinal is a set $\alpha \subseteq \N$ with a linear order $<_\alpha$ that is \textit{well-founded}: there is no infinite descending sequence $a_0 >_\alpha a_1 >_\alpha a_2 >_\alpha \cdots$ in $\alpha$.
\end{definition}

Ordinals are important in mathematics because we can induct on them: if $\varphi(\alpha)$ is a statement about ordinals, such that $\varphi(0)$ holds, and $\varphi(\beta)$ for all $\beta < \gamma$ implies $\varphi(\gamma)$, then $\varphi(\alpha)$ holds for \textit{all} ordinals $\alpha$. One classic and important result of set theory is that any two ordinals are \textit{comparable}: either they are isomorphic, or one is isomorphic to a strict initial segment of the other. $\ACA$ is not strong enough to prove this result, which is part of the motivation for introducing $\ATR$.

Now, let $\theta(n,X)$ be an arithmetical formula, with one free number variable $n$ and one free set variable $X$. $\theta$ defines an ``arithmetical operator'' $\Theta\colon \Pow(\N) \to \Pow(\N)$ by $$\Theta(X) = \{ n \in \N: \theta(n,X) \text{ holds} \}$$
For a set $Y \subseteq \N\!\times\! \alpha$, for each $\beta \in \alpha$, we let $Y^{[\beta]} = \{ n \in \N: (n,\beta) \in Y \}$ be the $\beta$th column of $Y$, and $Y^{[<\beta]} = \{ (n,\gamma) \in Y: \gamma <_\alpha \beta \}$ be all the columns up to $\beta$.

\begin{definition}\label{defn:atr}
	For a countable ordinal $\alpha \subseteq \N$ and set $X \subseteq \N$, let $\Theta^\alpha(X)$ be the subset\footnote{$\ACA$ is required to prove the set $\Theta^\alpha(X)$ is uniquely defined.} $Y \subseteq \N\!\times\! \alpha$ such that $Y^{[0]} = X$ and $Y^{[\beta]} = \Theta \big( Y^{[<\beta]} \big)$. Then, $\ATR$ is the subsystem consisting of $\ACA$, plus the assertion that $\Theta^\alpha(X)$ exists, for every arithmetical operator $\Theta$, countable ordinal $\alpha$, and set $X \subseteq \N$.
\end{definition}

%{\color{red} Leave this until we start talking about Baire 2 fns}

Of course, $\ATR$ implies $\ACA$, and thereby all the results of the previous section. Generally, $\ATR$ is needed to prove theorems which use ordinals in an essential way. An example is \textit{Ulm's theorem} on countable abelian $p$-groups: we assign each such group a countable ordinal $\alpha$, and a sequence in $\N \cup \{ \infty \}$ of length $\alpha$ called its \textit{Ulm invariant}. Ulm's theorem states that two countable abelian $p$-groups are isomorphic if and only if they have the same Ulm invariant; this theorem is equivalent to $\ATR$ \cite{friedmanCountableAlgebraSet1983}.

\begin{proposition}\label{prop:atr-prov}
	Over $\RCA$, $\ATR$ is equivalent to:
	\begin{enumerate}
		\item Any two countable ordinals are comparable \cite{steelDeterminatenessSubsystemsAnalysis1977};
		\item Every uncountable closed set in $\R$ has a perfect subset \cite{friedmanSystemsSecondOrder1974, friedmanSystemsSecondOrder1976};
		\item Determinacy for open or clopen sets in $\N^\N$ \cite{steelDeterminatenessSubsystemsAnalysis1977};
		\item Ramsey's theorem for open or clopen sets in $\N^\N$ \cite{friedmanFiniteCombinatorialPrinciple1982}.
	\end{enumerate}
\end{proposition}

\section{\texorpdfstring{$\PiCA$}{Pi-1-1-CA-0}}

The strongest system we will discuss is $\PiCA$, guaranteeing the existence of any $\Pi^1_1$ set:

\begin{definition}
	$\PiCA$ is the subsystem consisting of:
	\begin{enumerate}
		\item the basic axioms;
		\item the induction axiom for every $\Pi^1_1$ formula $\varphi$;
		\item the comprehension axiom for every $\Pi^1_1$ formula $\varphi$.
	\end{enumerate}
\end{definition}

\noindent In terms of subsystems of second-order arithmetic, $\PiCA$ is ``way up in the stratosphere''; it can prove almost any mathematical theorem the reader can imagine. $\PiCA$ implies $\ATR$, and hence all the results of the previous sections. Here are some further results equivalent to $\PiCA$; thus, they require some level of quantification over sets:

\begin{proposition}\label{prop:pica-prov}
	Over $\RCA$, $\PiCA$ is equivalent to:
	\begin{enumerate}
		\item The Cantor--Bendixson theorem \cite{friedmanSystemsSecondOrder1976};
		\item Every countable abelian group is the direct sum of a divisible group and a reduced group \cite{friedmanCountableAlgebraSet1983};
		\item Determinacy for sets of the form $U \setminus U'$, where $U$, $U'$ open in $\N^\N$ \cite{tanakaWeakAxiomsDeterminacy1991};
		\item Ramsey's theorem for $\Delta^0_2$ sets in $\N^\N$ \cite{solovayHyperarithmeticallyEncodableSets1978, simpsonSubsystemsSecondOrder2009};
		\item The minimal bad sequence lemma \cite{marconeLogicalStrengthNashWilliams1996};
		\item Maltsev's theorem: every countable ordered group has order type $\Z^\alpha$ or $\Z^\alpha \Q$ \cite{solomonPiCAOrderTypes2001};
		\item Every countable ring has a prime radical \cite{conidisComplexityRadicalsNoncommutative2009}.
	\end{enumerate}
\end{proposition}

\noindent As strong as $\PiCA$ is, there are still a few results which manage to escape it. These are generally restricted to select theorems in infinitary Ramsey theory, WQO theory, and set theory, where $\Pi^1_2$ or $\Pi^1_3$ comprehension might be required.
\chapter{Analysis in second-order arithmetic}
\label{chp:analysis-soa}

As described in Section \ref{sec:soa}, second-order arithmetic only includes two types of objects: natural numbers $n,m,k,\ldots$ and sets $A,B,C,\ldots$ thereof. Therefore, any other objects which we want to discuss must be \textit{coded} using natural numbers or subsets of $\N$. We've already seen an example in Section \ref{sec:wkl}---coding finite binary strings $\sigma \in \CSf$ by natural numbers, and infinite binary strings $X \in \CS$ using sets. In this section, we code the basic number systems $\Z$, $\Q$, $\R$ in second-order arithmetic, which then allows us to formalise basic concepts of analysis.

\section{Number systems}

The smallest number system is $\N$, and this is already given in second-order arithmetic, as the collection of all objects of number type. To code the integers $\Z$, we imitate the usual construction of $\Z$ from $\N$, where we use a pair $(a,b) \in \N \times \N$ to represent $a-b \in \Z$, and then quotient $\N \times \N$ by a suitable equivalence relation. To perform this construction, we first need to code pairs of natural numbers.

\begin{definition}
	For $m,n \in \N$, define the \textit{pair} $(m,n)$ as the natural number $(m+n)^2+m$.
\end{definition}

The reason we use this pairing function, rather than Cantor's pairing function of Definition \ref{defn:pair-fn}, is that the definition is more elementary, not requiring division, and thus easier to reason about. It has the disadvantage of not being a bijection.

\begin{definition}[{\cite[\S II.4]{simpsonSubsystemsSecondOrder2009}}]
	For pairs $(m,n)$, $(p,q) \in \N^2$, say $(m,n) =_\Z (p,q)$ if $m+q = n+p$. Then, an \textit{integer} is a pair $(m,n)$ which is minimal in its $=_\Z$-equivalence class.
\end{definition}

Again, the pair $(m,n)$ should be interpreted as the integer $m-n$. Instead of taking the equivalence classes as objects, we instead take minimal elements, as this way, integers can be represented by single natural numbers, rather than sets thereof. We can also define the standard arithmetic operations $+_\Z$, $-_\Z$, $\cdot_\Z$ and ordering $<_\Z$ on integers in the evident way. For example, $(m,n) -_\Z (p,q)$ is the pair $(m+q,n+p)$.

Having defined the integers $\Z$, we can now define the rationals $\Q$ from $\Z$, via the usual field of fractions construction:

\begin{definition}
	Let $\Z^+ = \{ x \in \Z: x >_\Z 0_\Z \}$. For pairs $(a,b)$, $(c,d) \in \Z \times \Z^+$, say $(a,b) =_\Q (c,d)$ if $a \cdot_\Z d =_\Z b \cdot_\Z c$. Then, a \textit{rational number} is a pair $(a,b)$ which is minimal in its $=_\Q$-equivalence class.
\end{definition}

Here, we interpret the pair $(a,b)$ as the rational number $a/b$. Again, the standard operations and relations $+_\Q$, $-_\Q$, $\cdot_\Q$, $<_\Q$, $\abs{\cdot}_\Q$ are defined as expected.

Moving to the real numbers $\R$, we have to change strategy, as we are now moving from countable to uncountable. It will not be possible to define real numbers as pairs of rationals, or even finite sequences of rationals, as there are too many reals. Instead, we define real numbers as certain \textit{infinite} sequences of rationals, mirroring the familiar Cauchy construction of $\R$ from $\Q$, with a small twist.

\begin{definition}
	Given sets $X, Y \subseteq \N$, a \textit{function} $f\colon X \to Y$ is a set of pairs $(x,y) \in X \times Y$, such that for all $x \in X$, there is a unique $y \in Y$ with $(x,y) \in f$.
\end{definition}

\begin{definition}\label{defn:real}
	A \textit{sequence of rationals} is a function $f\colon \N \to \Q$. By convention, we will instead denote sequences by $(q_i)$, where $q_i = f(i)$. A \textit{real number} is a sequence of rationals $(q_i)$ such that for all $m \leq n \in \N$, $\abs{q_m - q_n} \leq 2^{-m}$. We say $(q_i) =_\R (q'_i)$ if for all $k$, $\abs{q_k - q'_k} \leq 2^{-k+1}$.
\end{definition}

One may wonder why we require $\abs{q_m - q_n} \leq 2^{-m}$, and not just the usual Cauchy condition: for all $\varepsilon \in \Q^+$, there is $N$ such that for all $m, n \geq N$, $\abs{q_m - q_n} \leq \varepsilon$. The reason is that Definition \ref{defn:real} is modelled on the definition of a \textit{computable} real number---hence, it is the more suitable definition in weak systems such as $\RCA$ and $\WKL$. $\ACA$ is needed to prove the equivalence between Definition \ref{defn:real} and the usual Cauchy definition.

Also note that we are not picking a representative from each $=_\R$-class---this would require strong comprehension/choice axioms which we may not have access to. The standard arithmetic operations $+_\R$, $-_\R$, $\cdot_\R$, $\abs{\cdot}_\R$ can be defined ``pointwise'', and we say $(q_i) \leq_\R (q'_i)$ if for all $k$, $q_k \leq q'_k + 2^{-k+1}$. Furthermore, any rational $q$ can be identified with the real number $r_q = (q,q,q,\ldots)$.

%We now move to the setting of reverse mathematics. We assume $\RCA$ as a base system. We formalise analogues in $\RCA$ of the concepts developed thus far. Refer to [Simpson] for an overview of the basic concepts, such as the definitions of Cartesian products, functions, $\N^k$, $\Z$, $\Q$, $\R$, $<_\R$ etc. within $\RCA$.

%\begin{definition}
%	The \textit{rational unit interval} is the set $\UIQ \defeq \{ q \in \Q: 0 \leq_\Q q \leq_\Q 1 \}$, which exists by $\Delta^0_1$-comprehension.
%\end{definition}

\begin{definition}
	Given a real number $r$, we say $r \in [0,1]$ if $0 \leq_\R r \leq_\R 1$.
\end{definition}

\section{Open sets}

A key topological property of $\R$ is that it is \textit{second-countable}, i.e.\ its topology has a countable basis, consisting of open intervals $(p,q)$ with rational endpoints. This property is essential in allowing us to code open sets of $\R$ in second-order arithmetic. We first use natural numbers to code a basis of rational intervals for $\R$ and $\UI$:

\begin{definition}
	\phantom{hi}
	\begin{enumerate}
		\item A pair $(p,q) \in \Q \times \Q$, where $p < q$, codes the open interval $V_{p,q} \defeq (p,q) \subseteq \R$.
		
		\item For a real number $r$, we say $r \in V_{p,q}$ if $p < r < q$.
		
		\item We say $V_{p,q} \cap V_{p',q'} \neq \varnothing$ if $(p < q') \land (p' < q)$.
		
		\item We say $V_{p,q} \subseteq V_{p',q'}$ if $(p \geq p') \land (q \leq q')$.
		
		\item The \textit{length} of $V_{p,q}$ is $\ell \big( V_{p,q} \big) \defeq q-p$.
		
		\item $\B_\R \subseteq \Q \times \Q$ denotes the set of all such intervals.
	\end{enumerate}
\end{definition}

\begin{definition}
	\phantom{hi}
	\begin{enumerate}
		\item Given $p, q \in \Q$, we define $\overline{p} \defeq \max \{ p,0 \}$ and $\overline{q} \defeq \min \{ q,1 \}$.
		
		\item We also use $(p,q)$, where $\overline{p} < \overline{q}$, to code the open\footnote{In $[0,1]$ with the standard subspace topology.} interval $U_{p,q} \defeq (p,q) \cap [0,1]$.
		
		\item For a real number $r$, we say $r \in U_{p,q}$ if $r \in \UI$ and $p < r < q$.
		
		\item We say $U_{p,q} \cap U_{p',q'} \neq \varnothing$ if $(\overline{p} < \overline{q}') \land (\overline{p}' < \overline{q})$.
		
		\item We say $U_{p',q'} \subseteq U_{p,q}$ if $(p<0\, \lor\, p \leq p') \land (q>1\, \lor\, q \geq q')$.
		
		\item The \textit{length} of $U_{p,q}$ is $\ell \big( U_{p,q} \big) \defeq \max \{ \overline{q} - \overline{p}, 0 \}$.
		
		\item $\B_\UI \defeq \Q \times \Q$ denotes the set of all such intervals.
	\end{enumerate}
\end{definition}

\noindent Having coded the basis elements into the model, we can now define arbitrary open sets:

\begin{definition}
	An \textit{open set} $O \subseteq \R$ is a sequence $(U_i)$ of open intervals $U_i \in \B_\R$, i.e.\ a function $f\colon \N \to \B_\R$.
\end{definition}

The sequence $(U_i)$ should be interpreted as the open set $O = \bigcup_{i \in \N} U_i$. We will use the same sequence to code the closed set $C = \R\setminus O$. Relatively open and closed sets in $\UI$ are defined the same way, starting from $\B_\UI$.

\section{The proof of Cousin's lemma, revisited}
\label{sec:cl-revis}

Let's again look at the proof of Cousin's lemma (Lemma \ref{lem:cl-og}), and attempt to formalise it in second-order arithmetic. We have not yet given a formal definition of a \textit{function} $f\colon \UI \to \R$ in second-order arithmetic; for now, let us take it to be a primitive, undefined notion.

\begin{definition}
	A function $f\colon \UI \to \R$ is a \textit{gauge} if for all $x \in \UI$, we have $f(x) > 0$. This property will be denoted $f\colon \UI \to \R^+$.
\end{definition}

\begin{definition}
	A \textit{tagged partition of $\UI$} is a finite, odd-length sequence of reals $$P = \ang{x_0, t_0, x_1, t_1, \ldots, x_{\ell-1}, t_{\ell-1}, x_\ell}\ \subseteq\ \R^{2\ell+1}$$ such that $x_0 = 0$, $x_\ell = 1$, and for all $j < \ell$, we have $x_j < t_j < x_{j+1}$. The number $\ell$ is called the \textit{size of $P$}. 
\end{definition}

\begin{definition}
	Let $\delta\colon \UI \to \R^+$ be a gauge, and $P$ be a tagged partition of size $\ell$. Then, we say $P$ is \textit{$\delta$-fine} if for all $j < \ell$, $t_j - \delta(t_j) \leq x_j$ and $t_j + \delta(t_j) \geq x_{j+1}$.
\end{definition}

%\begin{definition}
%	Let $\CLc$ be the following statement in $\RCA$: every total continuous gauge $\delta\colon \UI \to \R^+$ has a $\delta$-fine partition.
%\end{definition}

We can now conduct the proof of Lemma \ref{lem:cl-og} in second-order arithmetic.

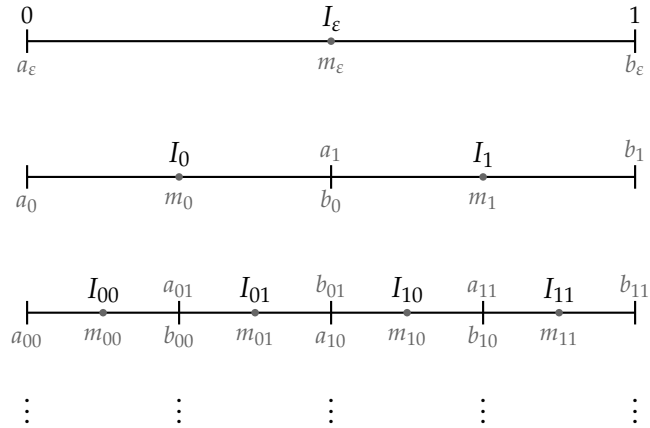
\begin{figure}[b]
	\centering
	\begin{tikzpicture}[x=80mm,y=12mm]%[scale=0.8]%[y=20mm]
	
	\definecolor{gg}{gray}{0.4}
	\def\he{0.13}
	\def\ya{-1.5}
	\def\yd{-4}
	
	\foreach \s/\a/\l/\p in {
		\varepsilon/0/0/0,
		0/0/1/0,
		1/0.5/1/1,
		00/0/2/0,
		01/0.25/2/1,
		10/0.5/2/0,
		11/0.75/2/1
	}{
		\pgfmathsetmacro\y{\l*\ya}
		\pgfmathsetmacro\b{\a+2^-\l}
		\pgfmathsetmacro\m{(\a+\b)/2}
		\pgfmathsetmacro\ys{-9+19*\p}
		
		\draw[thick] (\a,\y) -- (\b,\y);
		\draw[thick] (\a,\y+\he) -- (\a,\y-\he);
		\node[yshift=\ys,gg] at (\a,\y) {\footnotesize$a_{\s}\vphantom{b}$};
		\node[yshift=\ys,gg] at (\b,\y) {\footnotesize$b_{\s}$};
		\node[circle,gg,fill=gg,inner sep=1pt] at (\m,\y) {};
		\node[above,yshift=0mm] at (\m,\y) {$I_{\s}$};
		\node[yshift=-8,gg] at (\m,\y) {\footnotesize$m_{\s}\vphantom{b}$};
%		\draw[thick] (\a,{\n*\ya}) -- (\a+0.1,{\n*\ya});
%		\node at (\a,\n*\ya) {$k_\s$};
	}

	\node[above,yshift=1mm] at (0,0) {\footnotesize$0$};
	\node[above,yshift=1mm] at (1,0) {\footnotesize$1$};

	\foreach \l in {0,1,2}{
		\pgfmathsetmacro\y{\l*\ya}
		\draw[thick] (1,\y+\he) -- (1,\y-\he);
	}

	\foreach \i in {0,...,4}{
		\node at (\i/4,\yd) {$\vdots$};}
%	
%	
%	%\draw[thick,red] (0,0) -- (4,0);
%	%\draw[thick,red] (0,\he) -- (0,-\he);
%	%\draw[thick,red] (4,\he) -- (4,-\he);
%	\node[above,yshift=-0.5mm] at (2,0) {$I_0$};
%	%
%	%\draw[thick,red] (2,\ya) -- (4,\ya);
%	%\draw[thick,red] (2,\ya+\he) -- (2,\ya-\he);
%	%\draw[thick,red] (4,\ya+\he) -- (4,\ya-\he);
%	\node[above,yshift=-0.5mm] at (3,\ya) {$I_1$};
%	%
%	%\draw[thick,red] (2,\yb) -- (3,\yb);
%	%\draw[thick,red] (2,\yb+\he) -- (2,\yb-\he);
%	%\draw[thick,red] (3,\yb+\he) -- (3,\yb-\he);
%	\node[above,yshift=-0.5mm] at (2.5,\yb) {$I_2$};
%	\node[above,yshift=-0.5mm] at (2.75,\yc) {$I_3$};
	\end{tikzpicture}
	\caption[The definition of $a_\sigma$, $b_\sigma$, $I_\sigma$, $m_\sigma$ in the proof of Theorem \ref{thm:cl-pica}.]{The definition of $a_\sigma$, $b_\sigma$, $I_\sigma$, $m_\sigma$ in the proof of Theorem \ref{thm:cl-pica}; $\varepsilon$ is the empty string.}
	\label{fig:i-sigma}
\end{figure}

\begin{theorem}[($\PiCA$)]\label{thm:cl-pica}
	Any gauge $\delta\colon \UI \to \R^+$ has a $\delta$-fine partition.
\end{theorem}

\begin{proof}[``Proof'']
	For each $\sigma \in 2^\omega$, we define rationals $a_\sigma \defeq \sum_{i<n} \sigma_i \cdot 2^{-i-1}$ and $b_\sigma \defeq a_\sigma + 2^\abs{\sigma}$, and the interval $I_\sigma = (a_\sigma,b_\sigma)$. For convenience, we also let $m_\sigma = ( a_\sigma + b_\sigma )/2$, the midpoint of $I_\sigma$. So the strings $\sigma$ of length $n$ partition $\UI$ into $2^n$ subintervals $I_\sigma$ of equal length $2^{-n}$, as shown in Figure \ref{fig:i-sigma}.
	
	We define a tree $T$ in levels as follows. For each $n \in \N$, inductively define
	\begin{equation*}
	T_n\ \defeq\ \Big\{ \sigma \in \CSf :\ \abs{\sigma}=n,\ \forall k < n\ (\substr{\sigma}{k} \in T_k),\ \forall r \in I_\sigma\ \big( \delta(r)\leq 2^{-n} \big) \Big\}
	\end{equation*}
	Each $T_n$ exists by $\Pi^1_1$ comprehension. Take $T = \bigcup_{n \in \N} T_n$. Then, $T \subseteq \CSf$ is downward closed by construction, hence a tree.
	
	If $T$ is finite, then $T$ defines a $\delta$-fine partition $P_T$ of $\UI$ as follows: let $(\sigma^{(0)},\ldots,\sigma^{(n-1)})$ be a lexicographically sorted list\footnote{There are finitely many such $\sigma$ if $T$ is finite.} of all the $\sigma \notin T$ such that $\substr{\sigma}{k} \in T$ for all $k < \abs{\sigma}$. $P_T = \ang{x_0<t_0<\ldots<x_n}$ is defined by letting $x_i = a_{\sigma^{(i)}} = b_{\sigma^{(i-1)}}$, and $t_i = m_{\sigma^{(i)}}$.
	
	Now, we claim that $T$ must be finite. If not, then $\WKL$ proves there is an infinite path $X = X_0 X_1 \cdots$ through $T$. Define the real $r_X = (q_n)_{n \in \N}$, where each $q_n = m_{\substr{X}{n}}$. Note that $r_X \in I_{\substr{X}{n}}$ for every $n$. Hence, by the definition of $T$, $\delta(r_X) \leq 2^{-n}$ for every $n$, whence $\delta(r_X)=0$, contradicting the fact that $\delta$ is a gauge.
\end{proof}

$\PiCA$ was required when defining the $T_n$; we used universal quantification over real numbers. The reason this is a ``proof'', and not a proof, is that second-order arithmetic is unable to talk about \textit{arbitrary} functions $f\colon \UI \to \R$; these require uncountably much information to specify (i.e.\ where $f$ sends each point in $\UI$). Furthermore, there are $2^\cc$-many functions $\UI \to \R$; too many to code even using subsets of $\N$!

However, we will see in Sections \ref{chp:cl-cont} and \ref{chp:cl-baire} that second-order arithmetic can describe certain types of functions $f\colon \UI \to \R$. Essentially, we can formalise any class of functions that can be specified by countable information---examples include continuous functions, functions of a given Baire class, and Borel functions. The above ``proof'' shows that $\PiCA$ is an upper bound on the axiomatic strength of Cousin's lemma, for \textit{any} class of functions which can be defined in second-order arithmetic. We will see that this upper bound is often far from tight; in particular, Cousin's lemma for continuous functions can be proven in $\WKL$, a much weaker system than $\PiCA$.
\chapter{Cousin's lemma for continuous functions}
\label{chp:cl-cont}

\newcommand{\spk}{\mathrm{sp}}
\newcommand{\HB}{\mathsf{HB}}%_\Q}

Having outlined the main concepts of reverse mathematics in Chapter \ref{chp:soa}, and begun to formalise analysis in Chapter \ref{chp:analysis-soa}, we are now ready for a formal discussion of Cousin's lemma in second-order arithmetic. First, we will define continuous functions in second-order arithmetic, then determine the axiomatic strength of $\CLc$, Cousin's lemma for continuous functions. We will construct an explicit example showing that $\CLc$ fails in $\RCA$, and then prove the equivalence between $\CLc$ and $\WKL$ over $\RCA$.

%Many definitions given here are taken from \cite{simpsonSubsystemsSecondOrder2009}. As is customary in reverse mathematics, we will take $\RCA$ as a base system, hence we will make all definitions in $\RCA$ where possible.

\section{Continuous functions}

Since we wish to do real analysis in second-order arithmetic, we would now hope to be able to code \textit{functions} $\UI \to \R$. But, as discussed in the previous section, second-order arithmetic cannot describe \textit{arbitrary} functions $f\colon \UI \to \R$; for this, we would require \textit{third-order} arithmetic. However, certain types of functions $f\colon \R \to \R$ can be. It is known that any continuous function $f\colon A \to B$, with $A$ separable and $B$ Hausdorff, is uniquely determined by its values on a countable dense subset of $A$. As $\UI$ is separable and $\R$ is Hausdorff, this means continuous functions $f\colon \UI \to \R$ can be specified by countable information, and thus coded by subsets of $\N$.

We will only be concerned with continuous functions $f\colon \UI \to \R$, but exactly the same idea can be used to code continuous functions $f\colon \R \to \R$. Our method will be to code $f\colon \UI \to \R$ by the collection of \textit{pairs} of rational open intervals $(U,V)$ such that $f(U) \subseteq V$.

%In the following definitions, we abuse notation by using $f$ to denote both the code for a continuous function (i.e. a set of paired intervals), as well as the actual function itself.

\begin{definition}[{\cite[Defn II.6.1]{simpsonSubsystemsSecondOrder2009}}]\label{defn:cont-func}
	A \textit{(partial) continuous function} $f\colon \UI \to \R$ is a subset $f \subseteq \B_\UI \times \B_\R$ satisfying the following:
	
	\begin{enumerate}
		\item $(U,V) \in f$ and $(U,V') \in f \implies V \cap V' \neq \varnothing$;
		
		\item $(U,V) \in f$ and $U' \subseteq U \implies (U',V) \in f$;
		
		\item $(U,V) \in f$ and $V \subseteq V' \implies (U,V') \in f$;
		%		
		%		\item $\forall x \in \UI\ \forall \varepsilon \in \Q^+\ \exists U \in \B_\UI\ \exists V \in \B_\R\ \big( x \in U\, \land\, \ell(V) < \varepsilon\, \land\, f(U) \subseteq V \big)$ \label{defn:func-rca-tot}
	\end{enumerate}
\end{definition}

\noindent Again, we should interpret $(U,V) \in f$ (in the formal sense) to mean ``$f(U) \subseteq V$'' (in the colloquial sense). Such functions may be partial because they may not give us enough information to define $f(x)$ at a point $x \in \UI$. For example, the collection $g = \{ (U,V): U \in \B_\UI,\ V \supseteq (0,1) \}$ meets the conditions of Definition \ref{defn:cont-func}, but for any real $x \in \UI$, we only know that $g(x) \in (0,1)$; $g$ does not give us enough information to localise $g(x)$ more than this.

\begin{definition}
	Let $f\colon \UI \to \R$ be a partial continuous function. A real $x \in \UI$ is \textit{in the domain of $f$} if, for all $\varepsilon \in \Q^+$, there is a pair $(U,V) \in f$ such that $x \in U$ and $\ell(V) \leq \varepsilon$. If all $x \in \UI$ are in the domain of $f$, we say $f$ is \textit{total}.
\end{definition}

If $x$ is in the domain of $f$, we define $f(x)$ as the real $(q_n)_{n \in \N}$ obtained as follows. For each $n$, let $\varepsilon = 2^{-n}$, and for the least $V = V_{p,q}$ witnessing the above, let $q_n = (p+q)/2$. Using the assumptions in Definition \ref{defn:cont-func}, we can verify that $(q_n)$ satisfies Definition \ref{defn:real}.

Let's formalise some basic examples of continuous functions in second-order arithmetic, and check that they satisfy Definition \ref{defn:cont-func}.

\begin{proposition}[($\RCA$)]\label{prop:lin-cont}
	For any $m,c \in \Q$, the linear function $f\colon x \mapsto mx + c$ is total continuous.
\end{proposition}

\begin{proof}
	If $m=0$, this is simply the constant function $x \mapsto c$. Then $f = \{ (U,V): U \in \B_\UI,\ c \in V \}$ is continuous, total, and has $f(x) =_\R c$ for all $x \in \UI$.
	
	Now, suppose $m > 0$. Define $f\colon \UI \to \R$ by letting $\big( U_{p,q},V_{r,s} \big) \in f$ if and only if $m\overline{p} + c > r$ and $m\overline{q} + c < s$. We verify Definition \ref{defn:cont-func}:
	\begin{enumerate}
		\item If $\big( U_{p,q},V_{r,s} \big), \big( U_{p,q},V_{r',s'} \big) \in f$, then $r < m\overline{p} + c < m\overline{q} + c < s'$ since $m>0$ and $\overline{p} < \overline{q}$. We have $r'<s$ similarly, hence $V_{r,s} \cap V_{r',s'} \neq \varnothing$.
		
		\item Suppose $\big( U_{p,q},V_{r,s} \big) \in f$ and $U_{p',q'} \subseteq U_{p,q}$. Then, either $p<0$, in which case $\overline{p} = 0 \leq \overline{p}'$ or $p \leq p'$, in which case $\overline{p} \leq \overline{p}'$ also. Either way, we get $m\overline{p}' + c \geq m\overline{p} + c > r$. It follows similarly that $m\overline{q}' + c < s$, hence $\big( U_{p',q'},V_{r,s} \big) \in f$ as required.
		
		\item If $\big( U_{p,q},V_{r,s} \big) \in f$ and $V_{r,s} \subseteq V_{r',s'}$, then $m\overline{p} + c > r \geq r'$ and $m\overline{q} + c < s \leq s'$. It follows that $\big( U_{p,q},V_{r',s'} \big) \in f$.
		
		\item $f$ total: pick a real $x = (x_n)_{n \in \N} \in \UI$ and $\varepsilon \in \Q^+$. Let $n$ be least such that $2^{-n} < \varepsilon/2m$. Let $p \defeq x_n - 2^{-n}$, $q \defeq x_n + 2^{-n}$, $r \defeq mx_n + c - \varepsilon/2$ and $s \defeq mx_n + c + \varepsilon/2$. Then, the pair $\big( U_{p,q}, V_{r,s} \big) \in f$ is as required.
	\end{enumerate}
	
	If $m<0$, then we let $\big( U_{p,q},V_{r,s} \big) \in f$ if and only if $m\overline{q} + c > r$ and $m\overline{p} + c < s$, and the proof follows similarly.
\end{proof}

Within $\RCA$, we also have some methods of constructing new continuous functions from old:

\begin{lemma}[($\RCA$)]\label{lem:sum-cont}
	If $\sum_{n=0}^\infty \alpha_n$ is a convergent series of nonnegative real numbers, and $(f_n)_{n \in \N}$ is a sequence of continuous functions $\UI \to \R$ such that $\abs{f_n(x)} \leq \alpha_n$ for all $x \in \UI$, $n \in \N$, then $f = \sum_{n=0}^\infty f_n$ is continuous. Furthermore, $f$ is total if all the $f_n$ are total.
\end{lemma}

\begin{proof}
	\cite[Lemma II.6.5]{simpsonSubsystemsSecondOrder2009}.
\end{proof}

\begin{lemma}[($\RCA$)]\label{lem:piecewise-cont}
	Let $0 = d_0 < d_1 < d_2 < \cdots < d_k = 1$ be a finite, increasing sequence of rationals, and let $f_1, \ldots, f_k\colon \UI \to \R$ be continuous such that $f_i(d_i) = f_{i+1}(d_i)$ whenever $0 < i < k$. Then, the piecewise function $f$ defined $f(x) = f_i(x)$ for $d_{i-1} \leq x \leq d_i$ is also continuous. Furthermore, $f$ is total if all of the $f_i$ are total.
\end{lemma}

\begin{proof}
	Formally, we construct the code for $f$ in $\RCA$, as per Definition \ref{defn:cont-func}. Let $(U,V) \in f$ if and only if $(U,V) \in f_i$ for all $i$ such that $[d_{i-1},d_i]$ intersects $U$. One can easily prove that $f$ is partial continuous.
	
	Now, suppose each $f_i$ is total. To prove $f$ is total, pick some $x \in \UI$ and $\varepsilon \in \Q^+$. There are three cases:
	\begin{description}
		\item[Case 1:] $x \neq d_i$ for any $i \leq k$. Then, there is a unique $i$ such that $d_{i-1} < x < d_i$. Pick $(U,V) \in f_i$ witnessing that $f_i$ is total for this $x$ and $\varepsilon$, i.e. $\ell(V) \leq \varepsilon$. If $U = U_{p,q}$, then let $p' \defeq \max\{ p, d_{i-1} \}$, $q' \defeq \max\{ q, d_i \}$, and $U' \defeq U_{p',q'}$. Then $U' \subseteq U$, so by Definition \ref{defn:cont-func}.(ii), $(U',V) \in f_i$ also, and $x \in U'$. Since $U'$ only intersects $[d_{i-1}, d_i]$, it follows that $(U',V) \in f$ is as required.

		\item[Case 2:] $x = d_i$ for $0 < i < k$. We pick $(U_0,V_0) \in f_i$ and $(U_1,V_1) \in f_{i+1}$ witnessing that $f_i$ (resp. $f_{i+1}$) is total for this $x$ and $\varepsilon/2$. $U \defeq U_0 \cap U_1$ is a nonempty interval since $x \in U$, and $V \defeq V_0 \cup V_1$ is an interval since $f(x) \in V_0 \cap V_1$. Now let $U' = U \cap (d_{i-1},d_{i+1})$. We then have $\ell(V) \leq \varepsilon$ and $(U',V) \in f_i \cap f_{i+1} \implies (U',V) \in f$.

		\item[Case 3:] $x = d_0$ or $x = d_k$. A similar argument to \textbf{{Case 1}} works.\qedhere
	\end{description}
	%	\begin{itemize}
	%		\item $(U,V) \in f_i$, if $U = U_{p,q}$ and $0 < i \leq k$ are such that $d_{i-1} \leq \overline{p}$ and $\overline{q} \leq d_i$; or
	%		
	%		\item $(U,V) \in f_i$, if $U = U_{p,q}$ and $0 < i \leq k$ are such that $d_{i-1} \leq \overline{p}$ and $\overline{q} \leq d_i$;
	%	\end{itemize}
\end{proof}

%The proof of Theorem \ref{thm:cl-markov} can be carried out in $\RCA$ {\color{red} (???)}, which shows:
%
%\begin{theorem}
%	$\CLc$ is not provable in $\RCA$.
%\end{theorem}
%

\section{\texorpdfstring{$\CLc$}{CLc} fails in \texorpdfstring{$\RCA$}{RCA-0}}

We are now ready to discuss Cousin's lemma for continuous functions. The usual definitions of gauge, partition, and $\delta$-fine (as in Section \ref{sec:cl-revis}) are still valid here.

\begin{definition}
	Let $\CLc$ be the following statement in $\RCA$: every total continuous gauge $\delta\colon \UI \to \R^+$ has a $\delta$-fine partition.
\end{definition}

\begin{remark}
	When we say that a theorem $\varphi$ \textit{holds} in a model $\M$ of second-order arithmetic, all the quantification and interpretations should be made relative to \textit{objects in $\M$}. So, to say that $\CLc$ holds in $\M$ is to say that for every object $\delta$ \textit{in $\M$}, which $\M$ \textit{believes} to be total, continuous and a gauge, there is a object $P$ \textit{in $\M$}, which $\M$ \textit{believes} to be a finite, $\delta$-fine sequence of reals.
	
	If we wanted to show $\CLc$ fails in $\M$, then we would need to construct an object $\delta$ \textit{in $\M$}, which $\M$ \textit{believes} to be total, continuous and a gauge, and so that there is no object $P$ \textit{in $\M$} that $\M$ \textit{thinks} is a $\delta$-fine partition. To show that $\CLc$ fails in a subsystem $\Ss$ of second-order arithmetic, we need to demonstrate a model $\M$ of $\Ss$ where $\CLc$ fails, in the sense just described.
\end{remark}

We want to show that $\RCA$ does not prove $\CLc$; to do this, we exhibit a model of $\RCA$ where $\CLc$ doesn't hold. In fact, this is true in the standard model $\REC$ of recursive sets (Proposition \ref{prop:rec-rca}). Recall $\varphi_0, \varphi_1, \varphi_2, \ldots$ is a standard enumeration of the partial computable functions (Proposition \ref{prop:utm}). Via the coding of $\Q$ into $\N$, we can assume WLOG that the $\varphi_e$ take values in $\Q$.% We also fix computable enumerations $q_0, q_1, q_2, \ldots$ of $\Q$, and $U_0, U_1, U_2, \ldots$ of $\B_\UI$, the rational open intervals in $\UI$.

To construct our counterexample, we will use the idea of a \textit{$\Pi^0_1$ class} from classical computability.

\begin{definition}[\cite{jockuschDegreesMembersPi1972}]%[Cenzer--Remmel, \S 4]
	A \textit{$\Pi^0_1$ class in $\UI$} is a set of the form
	\begin{equation*}
		K_\psi\ \defeq\ \UI \setminus \bigcup_{n=0}^\infty U_{\psi(n)}
	\end{equation*}
	for some computable function $\psi\colon \N \to \Q^2$.
\end{definition}

Recall that in second-order arithmetic, real numbers are defined as fast-converging sequences $(q_i)_{i \in \N}$ of rational numbers. We say a real number is \textit{computable} if this sequence is computable, considered as a function $f\colon \N \to \Q$. These are exactly the real numbers that exist in $\REC$. The following result is closely related to Proposition \ref{prop:rec-wkl}.

\begin{lemma}[\cite{jockuschDegreesMembersPi1972}]\label{lem:pi-0-1-disjoint}
	There exists a nonempty $\Pi^0_1$ class which contains no computable reals.
\end{lemma}

\begin{proof}
	We define $\psi$ as follows: search over all pairs $(e,s) \in \N^2$ until we find the next one such that $\varphi_e(e+3)$ halts after $s$ steps. When we find such a pair, let $\psi(n)$ be the code for the rational open ball $B_e \defeq B(\varphi_e(e+3),2^{-e-3})$. We can always find another such pair, so in particular, $\psi$ is total computable.
	
%	Let $W \subseteq \N$ be the set of all codes for balls $B_e \defeq B(\varphi_e(e+3),2^{-e-3})$ where $\varphi_e(e+3)$ is defined. Then, $W$ is a c.e. set of natural numbers, so let $\psi$ be a total recursive function enumerating $W$.
	%	Construct a recursive function $\psi$ as follows: for each $e \in \N$, let $y_e \defeq \varphi_e(e+3)$ if this is defined, else $y_e = \nicefrac{1}{2}$. Then, let $\psi(e)$ be the code for the open ball $B(y_e,2^{-e-3})$.
	We claim that $K_\psi$ has the required properties.	
	First, note that each $B_e$ has Lebesgue measure $\lambda(B_e) = 2^{-e-2}$, so their union has measure at most $\frac{1}{2}$. In particular, the complement $K_\psi$ must be nonempty.
	
	%	Construct a recursive function $\psi$ as follows: for each $e \in \N$, compute $y \defeq \varphi_e(e+3)$. If $y$ is defined, let $\psi(e)$ be the code for the open ball $B(y,2^{-e-3})$. Otherwise, let $\psi(e)$ be the code for some open ball of small radius $\varepsilon$, e.g. $(1.234,1.235)$. We then claim that $K_\psi$ has the required properties.
	
	%	First, note that each $U_{\psi(e)}$ has Lebesgue measure $\lambda(U_{\psi(e)}) = \max \{ 2^{-e-2}, 2\varepsilon \}$, so their union has measure at most $\nicefrac{1}{2} + 2\varepsilon$. In particular, the complement $K_\psi$ must be nonempty.
	
	Now, suppose $r = (q_n)_{n \in \N}$ is a \textit{computable} real number; then $r$ is computed by some $\varphi_e$, i.e. $\varphi_e(n) = q_n$ for all $n$. By the definition of real number, $q_{e+3}$ is an approximation of $r$ to within $2^{-e-3}$. Thus, $r \in B(q_{e+3},2^{-e-3}) = B_e$, so $r \notin K_\psi$.
	%We also see that $K_\psi$ is disjoint from $\UIc$. Fix $r \in \UIc$, and suppose $\varphi_e$ computes a Cauchy name of $r$. Then $q \defeq \varphi_e(e+3)$ is a rational approximation of $r$ to within $2^{-e-3}$. Thus, $r \in B(q,2^{-e-3}) = B_e$, implying $r \notin K_\psi$.
\end{proof}

\begin{figure}[b]
	\centering
	\begin{tikzpicture}
	\begin{axis}[
	xmin=-0.1,xmax=1.1,
	ymin=-0.2,ymax=0.27,
	width=8cm,height=4cm,
	axis x line=center,axis y line=none,
	axis line style={black!30!white,{Latex[width=1.5mm]}-{Latex[width=1.5mm]}},
	%		every x tick/.style={lightgray},
	xtick={0,1},
	axis equal,
	]
		\definecolor{gg}{gray}{0.4}
	\addplot[samples at={0,0.3,0.5,0.7,1},black] {max(0.2-abs(x-0.5),0)};
	%		\addplot[samples at={0,0.22,0.45},red,dashed] {abs(x-0.22)};
	\draw[gg,dashed] (axis cs:0.22,0) -- (axis cs:-0.04,0.26) {};
	\draw[gg,dashed] (axis cs:0.22,0) -- (axis cs:0.48,0.26) {};
	\node[label={[text=gg]270:$b$},circle,fill,gg,inner sep=0pt,minimum size=1mm] at (axis cs:0.22,0) {};
	\draw[{Parenthesis[length=1.5mm,width=4mm,color=black]}-{Parenthesis[length=1.5mm,width=4mm,color=black]},black!30!white] (axis cs:0.3,-0.08) -- (axis cs:0.7,-0.08) node[pos=0.5,below,black] {$U_{\psi(n)}$};
	\end{axis}
	\end{tikzpicture}
	\caption{The $n$th spike $\spk_n$, compared to $z \mapsto \abs{z-b}$ for a point $b \notin U_{\psi(n)}$.}
	\label{fig:spk}
\end{figure}
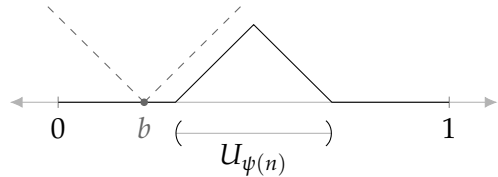

Now, we are ready to construct our counterexample, to show $\CLc$ fails in $\RCA$. The construction was inspired by \cite[Thm 3.1]{koComplexityTheoryReal1991}, and the idea is as follows. Given a $\Pi^0_1$ class $K_\psi$ as in Lemma \ref{lem:pi-0-1-disjoint}, we construct a continuous gauge $\delta$ which is positive \textit{exactly} on the complement of $K_\psi$, and furthermore is 1-Lipschitz ($\abs{\delta(x)-\delta(y)} \leq \abs{x-y}$). Then, $\REC$ will think that $\delta$ is a gauge, since it is positive on all computable reals. However, for any point $b \in K_\psi$ and $t_i \neq b$, it is not possible for $\delta(t_i)$ to cover $b$, since $\delta$ is 1-Lipschitz; therefore, there are no $\delta$-fine partitions in $\REC$. The formal proof follows.

\begin{theorem}\label{thm:rca-clc}
	$\RCA$ does not prove $\CLc$.
\end{theorem}

\begin{proof}
	Let $K_\psi$ be a $\Pi^0_1$ class as in Lemma \ref{lem:pi-0-1-disjoint}. For each $n$, define the \textit{$n$th spike} $\spk_n\colon \UI \to \R$ by
	%	\begin{equation*}
	%		\spk_n(z) = \begin{cases}
	%			0 &\quad z \leq \overline{p} \\
	%			\abs{z-\overline{p}} &\quad \overline{p} < z \leq m \\
	%			\abs{z-\overline{q}} &\quad m < z < \overline{q} \\
	%			0 &\quad z \geq \overline{q} \\
	%		\end{cases}
	%	\end{equation*}
	%	
	\setlength\arraycolsep{2pt}
	\begin{equation*}
	\spk_n(z) = \left\{ \begin{array}{lrcccl}
	0						&\qquad 0			&\leq &z &\leq &\overline{p}	\\
	\abs{z-\overline{p}}	&\qquad \overline{p}	&\leq &z &\leq &m			\\
	\abs{z-\overline{q}}	&\qquad m			&\leq &z &\leq &\overline{q}	\\
	0						&\qquad \overline{q}	&\leq &z &\leq &1			\\
	\end{array} \right.
	\end{equation*}
	where $(p,q)$ is the code of $U_{\psi(n)}$, and $m \defeq (\overline{p}+\overline{q})/2$ is the midpoint of $U_{\psi(n)}$.
	
	As in Figure \ref{fig:spk}, $\spk_n$ is graphically a spike whose base is exactly $U_{\psi(n)}$, and whose sides have gradient $\pm 1$. For each fixed $n$, $\spk_n$ is 1-Lipschitz (i.e. $\abs{\spk_n(x)-\spk_n(y)} \leq \abs{x-y}$) and bounded above by $\frac{1}{2}$. By Proposition \ref{prop:lin-cont}, each part is total continuous, so $\spk_n$ is total continuous by Lemma \ref{lem:piecewise-cont}.
	
	Then, we define $\delta\colon \UI \to \R$ by 
	\begin{equation*}
	\delta(x) = \sum_{n=0}^\infty\, 2^{-n-2} \cdot \spk_n(x)
	\end{equation*}
	which is total continuous by Lemma \ref{lem:sum-cont}. In $\REC$, $\delta$ is a gauge; for any real $x \in \REC$, $x$ is a computable real number, so $x$ is in some $U_\psi(n)$ by definition of $K_\psi$. Then, $\delta(x) \geq \spk_n(x) > 0$, since $U_{\psi(n)}$ is open.

We claim there is no $\delta$-fine partition in $\REC$. Suppose, by contradiction, that $P = \ang{x_i,t_i}$ is such a partition, of size $\ell$. By assumption, $K_\psi$ is nonempty, so pick any point $b \in K_\psi$, which is necessarily noncomputable. There is unique $m < \ell$ such that $b \in (x_m, x_{m+1})$; then, we claim $\delta(t_m) < \abs{t_m - b}$. This would imply that $b \notin B \big( t_m, \delta(t_m) 
\big)$ and thus $(x_m, x_{m+1}) \nsubseteq B \big( t_m, \delta(t_m) 
\big)$, so $P$ is not $\delta$-fine after all.

Note that for all $n \in \N$, $b \notin U_{\psi(n)}$, and so $\spk_n(t_m) \leq \abs{t_m - b}$ (see Figure \ref{fig:spk}). We compute:
\begin{equation*}
\delta(t_m)\ =\ \sum_{n=0}^\infty\, 2^{-n-2} \cdot \spk_n(t_m)\ \leq\ \sum_{n=0}^\infty\, 2^{-n-2} \abs{t_m - b}\ \leq\ \tfrac{1}{2} \abs{t_m - b}\ <\ \abs{t_m - b}.\vspace*{-5mm}%\qedhere
\end{equation*}

\end{proof}

\section{\texorpdfstring{$\WKL$}{WKL-0} proves \texorpdfstring{$\CLc$}{CLc}}
We saw in the previous section that $\RCA$ is not strong enough to prove Cousin's lemma for continuous functions, $\CLc$. In this section, we show that $\WKL$ \textit{is} strong enough to prove $\CLc$. In the next section, we demonstrate a reversal of $\CLc$ in $\WKL$, thus showing that $\CLc$ and $\WKL$ are equivalent, and that $\WKL$ is the weakest subsystem of second-order arithmetic in which $\CLc$ can be proved.

The idea is a variation on the ``proof'' of Theorem \ref{thm:cl-pica}. We avoid using $\Pi^1_1$ comprehension by only considering the \textit{midpoint} of each $I_\sigma$, rather than \textit{all} real numbers in $I_\sigma$. The proof then proceeds exactly the same. To deduce a contradiction, we will use the fact (provable in $\WKL$) that every continuous function $f\colon \UI \to \R$ is \textit{uniformly continuous}:

\begin{lemma}[($\WKL$)]\label{lem:wkl-unif-cont}
	Let $f\colon \UI \to \R$ be a total continuous function. Then $f$ has a \textit{modulus of uniform continuity}, i.e.\ a function $h\colon \N \to \N$ such that for all $x,y \in \UI$, $$\abs{x-y} \leq 2^{-h(n)} \implies \abs{f(x)-f(y)} \leq 2^{-n}$$
\end{lemma}

\begin{proof}
	\cite[Thm IV.2.2]{simpsonSubsystemsSecondOrder2009}.
\end{proof}

%, and we can use sequential continuity to deduce a contradiction:

%\begin{proposition}
%	If $f\colon \UI \to \R$ is continuous, then it is 
%\end{proposition}

\begin{theorem}\label{thm:wkl->clc}
	$\WKL$ proves $\CLc$.
\end{theorem}

\begin{proof}
	Let $\delta\colon \UI \to \R^+$ be a total continuous gauge. For each $x \in \UI$, $\delta(x)$ is a real number, hence a sequence of rationals; let $\substr{\delta}{n}(x)$ denote the $n$th term in this sequence. For each $\sigma \in \CSf$, let $$m_\sigma = 2^{\abs{\sigma}-1} + \sum_{i<n} \sigma_i \cdot 2^{-i-1}$$ be as in the ``proof'' of Theorem \ref{thm:cl-pica}.
	
	We define a tree $T \subseteq \CSf$ in levels. For each $n \in \N$, inductively define
	\begin{equation*}
		T_n\ \defeq\ \Big\{ \sigma \in \CSf :\ \abs{\sigma}=n,\ \forall k < n\ (\restr{\sigma}{k} \in T_k),\ \substr{\delta}{n}(m_\sigma) \leq 2^{-n+1} \Big\}
	\end{equation*}
	The $\leq$ relation between rationals is computable, so by $\Delta^0_1$ comprehension, each $T_n$ exists. Then, $T = \bigcup_{n=0}^\infty T_n$ is a tree by construction; if it is finite, we construct a $\delta$-fine partition as in Theorem \ref{thm:cl-pica} (this can all be done in $\RCA$).
	
%	If $T$ is finite, then we construct a tagged partition $P$ by ...
	
	Now, we claim $T$ must be finite, so the above construction always works. Suppose by contradiction that $T$ is infinite. By $\WKL$, there is an infinite path $X$ through $T$. Again, we define the real $r_X = (q_n)_{n \in \N}$, where each $q_n = m_{\substr{X}{n}}$. By definition of $T$, we have $\substr{\delta}{n}(q_n) \leq 2^{-n+1}$, and $\abs{\vphantom{\big(}\delta(q_n) - \substr{\delta}{n}(q_n) } \leq 2^{-n}$ by definition of $\substr{\delta}{n}(x)$, so for each $n$,
	\begin{equation}\label{eq:delta-qn}
		\delta (q_n)\ \leq\ 3 \cdot 2^{-n}
	\end{equation}
	
	By Lemma \ref{lem:wkl-unif-cont}, pick $h\colon \N \to \N$ a modulus of uniform continuity for $\delta$. By definition of $r_X$, for each $n, k \in \N$, we have
	$$\abs{q_{h(n)+k} - r_X}\ \leq\ 2^{-h(n)-k}\ \leq\ 2^{-h(n)}$$
	hence by definition of $h$,
	$$\abs{\delta(q_{h(n)+k}) - \delta(r_X)}\ \leq\ 2^{-n}$$
	Combining this with equation \eqref{eq:delta-qn}, we get that for all $n, k \in \N$, $$\delta(r_X)\ \leq\ 2^{-n}\ +\ 3 \cdot 2^{-h(n)-k}$$
	We can make this arbitrarily small by picking the right $n$ and $k$; hence $\delta(r_X) = 0$, contradicting the fact that $\delta$ is a gauge.% {\color{red} Use uniform continuity? Simpson Thm IV.2.2}
	%	\begin{equation*}
	%		T_0 = \begin{cases}
	%			\{ \varepsilon \} & \text{ if } $\restr{\delta}{n}(\nicefrac{1}{2}) > 1$ \\
	%			\varnothing & \text{ otherwise.}
	%		\end{cases}
	%	\end{equation*}
\end{proof}

\section{\texorpdfstring{$\CLc$}{CLc} is equivalent to \texorpdfstring{$\WKL$}{WKL-0}}
\label{sec:clc}

In the previous section we showed that $\CLc$ can be proved in $\WKL$. Now, we show that $\WKL$ is the weakest system having this property, by demonstrating a reversal of $\CLc$ in $\WKL$. The reversal goes through the Heine--Borel theorem $\HB$, which is known to be equivalent to $\WKL$.

\begin{definition}
	An \textit{open cover of $\UI$} is a (finite or infinite) sequence $(U_i)$ of open intervals in $\B_\UI$, such that every $x \in \UI$ is in some $U_i$.
\end{definition}

\begin{definition}[{\cite[Lem IV.1.1]{simpsonSubsystemsSecondOrder2009}}]
	Let $\HB$ be the following statement in $\RCA$: for every infinite open cover $(U_i)_{i \in \N}$ of $\UI$, there is $n$ such that $(U_i)_{i \leq n}$ is a finite open cover of $\UI$.
\end{definition}

\begin{proposition}[($\RCA$) {\cite[Lem IV.1.1]{simpsonSubsystemsSecondOrder2009}}]\label{prop:wkl-hb}
	$\WKL$ is equivalent to $\HB$.
\end{proposition}

Now, we show that over $\RCA$, $\CLc$ implies $\HB$, and thereby $\WKL$. The idea of the proof is similar to Theorem \ref{thm:rca-clc}; given an open cover $(U_i)_{i \in \N}$, we define $\delta$ on the $U_i$ in the same way. This time, we assume $\CLc$, so there is a $\delta$-fine partition $P = \ang{x_j, t_j}$. Then, for each $t_j$, we can find some $U_i$ such that $B \big( t_j, \delta(t_j) \big) \subseteq U_i$. Since the balls $\big( t_j, \delta(t_j) \big)$ cover $\UI$, it follows that the corresponding $U_i$ also cover $\UI$, so we get a finite subcover.

\begin{theorem}[($\RCA$)]\label{thm:cl->hb}
	$\CLc$ implies $\HB$.
\end{theorem}

% continuous piecewise functions
% weighted sum of cont functions is cont

\begin{proof}\renewcommand{\qedsymbol}{}
	Let $(U_i)_{i \in \N}$ be an open cover of $\UI$. We define $\delta$ as in the proof of Theorem \ref{thm:rca-clc}, where for each $i$, $\psi(i)$ is the code for $U_i$. Formally, define $\spk_i \colon \UI \to \R$ by
	%	\begin{equation*}
	%		\spk_n(z) = \begin{cases}
	%			0 &\quad z \leq \overline{p} \\
	%			\abs{z-\overline{p}} &\quad \overline{p} < z \leq m \\
	%			\abs{z-\overline{q}} &\quad m < z < \overline{q} \\
	%			0 &\quad z \geq \overline{q} \\
	%		\end{cases}
	%	\end{equation*}
	%	
	\setlength\arraycolsep{2pt}
	\begin{equation*}
	\spk_n(z) = \left\{ \begin{array}{lrcccl}
	0						&\qquad 0			&\leq &z &\leq &\overline{p}	\\
	\abs{z-\overline{p}}	&\qquad \overline{p}	&\leq &z &\leq &m			\\
	\abs{z-\overline{q}}	&\qquad m			&\leq &z &\leq &\overline{q}	\\
	0						&\qquad \overline{q}	&\leq &z &\leq &1			\\
	\end{array} \right.
	\end{equation*}
	where $(p,q)$ is the code for $U_i$, and $m \defeq (\overline{p}+\overline{q})/2$ is the midpoint of $U_i$. Again, each $\spk_n$ is total continuous.
	
	Then, we define $\delta\colon \UI \to \R$ by 
	\begin{equation*}
	\delta(x) = \sum_{n=0}^\infty\, 2^{-n-2} \cdot \spk_n(x)
	\end{equation*}
	which is total continuous by Lemma \ref{lem:sum-cont}. Furthermore, $\delta$ is a gauge, since for any $x \in \UI$, $x \in U_n$ for some $n$, so $\delta(x) \geq \spk_n(x) > 0$.
\end{proof}
	
\begin{thmclaim}\label{clm:delta-approx}
	Let $r = (q_n)_{n \in \N}$ be a real number, and for each $e \in \N$,
	\begin{equation*}%\label{eq:part-sum}
		y_e\ =\ \sum_{n=0}^e\, 2^{-n-2} \cdot \spk_n(q_e)
	\end{equation*}
	Then, for each $e \in \N$, $\abs{y_e - \delta(r)} \leq 2^{-e}$.
\end{thmclaim}

\begin{proof}[Proof of Claim \ref{clm:delta-approx}]
	Recall that each $\spk_n$ is 1-Lipschitz, i.e.\ $\abs{\spk_n(x) - \spk_n(y)} \leq \abs{x-y}$ for all $x,y \in \UI$; and $\spk_n(x) \leq \nicefrac{1}{2}$ for all $x \in \UI$. Using these facts, we compute:
	
%		$\delta$ is Markov computable: pick $r \in \Rc$, and a computable Cauchy name $(q_e)$ of $r$. We claim the sequence of partial sums
%		\begin{equation*}%\label{eq:part-sum}
%		y_e = \sum_{n=0}^e\, 2^{-n-2} \cdot \spk_n(q_e)
%		\end{equation*}
%		is a computable Cauchy name for $\delta(r)$. The sequence $(y_e)$ is computable, since $(q_e)$ and $\spk_n$ are. Now we show $\abs{y_e - \delta(r)} \leq 2^{-e}$ for each $e \in \N$:
	\begin{align*}
	\abs{y_e - \delta(r)}\ &=\ \abs{\ \sum_{n=0}^e\, 2^{-n-2} \big( \spk_n(q_e) - \spk_n(r) \big)\ -\ \sum_{n=e+1}^\infty\, 2^{-n-2} \cdot \spk_n(r)\ } \\[2mm]
	&\leq\ \sum_{n=0}^e\, 2^{-n-2} \abs{ \spk_n(q_e) - \spk_n(r) }\ +\ \sum_{n=e+1}^\infty\, 2^{-n-2} \cdot \spk_n(r) \\[2mm]
	&\leq\ \sum_{n=0}^e\, 2^{-n-2} \abs{ q_e - r }\ +\ \sum_{n=e+1}^\infty\, 2^{-n-2} \cdot \tfrac{1}{2} \\[2mm]
	&\leq\ 2^{-e}\, \sum_{n=0}^e\, 2^{-n-2}\ +\ 2^{-e-3} \\[2mm]
	&\leq\ 2^{-e} \cdot \tfrac{1}{2}\ +\ 2^{-e-3}\ =\ 2^{-e} \cdot \tfrac{5}{8}\ \leq\ 2^{-e}.\qedhere
	\end{align*}
\end{proof}

By assumption, there exists a $\delta$-fine partition $P = \ang{x_j,t_j}$. Let $\ell$ be the size of $P$.

\begin{thmclaim}\label{clm:finite-subc}
	For each $j<\ell$, there exists $m = m_j \in \N$ such that $\delta(t_j) < \spk_m(t_j)$.
\end{thmclaim}

\begin{proof}[Proof of Claim \ref{clm:finite-subc}]
	Let $y_e$ be as in the previous claim, for $r = t_j$, and let $e \defeq \min \{ k: y_k \geq (3k+1)2^{-k-1} \}$. This set is nonempty since $\delta(t_j) > 0$, so we can find $e$ in $\RCA$ by minimisation \cite[Thm II.3.5]{simpsonSubsystemsSecondOrder2009}. We must also have $e \geq 1$.
	
	We claim there is $m \leq e$ such that $\spk_m(q_e) > (3e+1)2^{-e}$. If there were not (i.e. $\spk_n(q_e) \leq (3e+1)2^{-e}$ for all $n \leq e$), then
	\begin{equation*}
	y_e\ =\ \sum_{n=0}^e\, 2^{-n-2} \cdot \spk_n(q_e)\ \leq\ \sum_{n=0}^e\, 2^{-n-2} (3e+1)2^{-e}\ <\ (3e+1)2^{-e-1}
	\end{equation*}
	contradicting the definition of $e$.
	
	We take $m_j$ to be the least such $m$, and claim this is as required. Because $\spk_m$ is 1-Lipschitz, we have $\abs{\spk_m(q_e) - \spk_m(t_j)} \leq \abs{q_e - t_j} \leq 2^{-e}$. We compute:
	\begin{align*}
	\delta(t_j)\ &\leq\ y_{e-1} + 2^{-e+1} &&\text{since $\abs{\delta(t_j) - y_{e-1}} \leq 2^{-e+1}$} \\
	&<\ (3e-2)2^{-e} + 2(2^{-e}) &&\text{by definition of $e$} \\
	&=\ 3e \cdot 2^{-e} \\
	&=\ (3e+1)2^{-e} - 2^{-e} \\
	&<\ \spk_m(q_e) - 2^{-e} &&\text{by definition of $m$} \\
	&\leq\ \spk_m(t_j) &&\text{since $\abs{\spk_m(q_e) - \spk_m(t_j)} \leq 2^{-e}$}\qedhere
	\end{align*}
\end{proof}
	
\begin{proof}[Proof of Theorem \ref{thm:cl->hb}, continued] 
	For each $j$, fix $m_j$ as in the claim. Then, $n \defeq \max \{ m_j : j < \ell \}$ gives a finite subcover. Taking $z \in \UI$, there is some $j < \ell$ such that $x_j \leq x \leq x_{j+1}$. Then, $t_j - \delta(t_j) \leq z \leq t_j + \delta(t_j)$ since $P$ is $\delta$-fine, i.e. $\abs{z - t_j} \leq \delta(t_j) < \spk_{m_j}(t_j)$ by the claim. It follows that $x \in U_{m_j}$.
	%for each $j < \ell$, let $t_j = (q_n)_{n \in \N}$ and $\delta(t_j) = (y_n)_{n \in \N}$. Then, there exists a pair $(i,n)$ such that $c_i <_\Q q_n - y_n - 2^{-n+1}$ and $d_i >_\Q q_n + y_n + 2^{-n+1}$. {\color{red} How to see this is true}
	%By minimisation [Simpson, Thm II.3.5], we can find the least pair $(i,n)$ satisfying the above. Then, define $f(j) = i$.
\end{proof}

\begin{theorem}[($\RCA$)]
	$\CLc$ is equivalent to $\WKL$.
\end{theorem}

\begin{proof}
	The forward direction is Theorem \ref{thm:cl->hb} and Proposition \ref{prop:wkl-hb}, while the reverse direction is Theorem \ref{thm:wkl->clc}.
\end{proof}
\chapter{Cousin's lemma for Baire functions}
\label{chp:cl-baire}

\newcommand{\SC}{\mathsf{SC}}

In Chapter \ref{chp:cl-cont}, we completely characterised the axiomatic strength of Cousin's lemma for continuous functions $\CLc$, showing its equivalence to $\WKL$ over $\RCA$. In this chapter, we will define the Baire classes of functions, and study the strength of $\CLB{n}$, Cousin's lemma for functions of a given Baire class $n$. In contrast to $\CLc$, the reverse mathematics of $\CLB{n}$ appears much harder, and so far has resisted complete characterisation for any $n \geq 1$.

\section{Baire classes of functions}

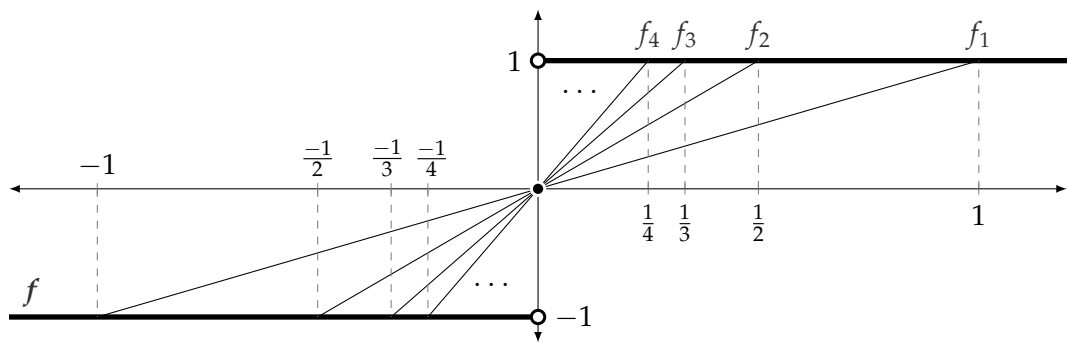
\begin{figure}
	\centering
	\begin{tikzpicture}% coordinates
	\begin{axis}[
		axis lines=center,
		axis line style={latex-latex},
		ymin=-1.2,ymax=1.4,
		xmin=-1.2,xmax=1.2,
		height=6cm,width=\textwidth,
		ytick={1}, extra y ticks = {-1},
		extra y tick style={y tick label style={right, xshift=0.4em}},
		xtick={1,1/2,1/3,1/4}, extra x ticks = {-1,-1/2,-1/3,-1/4},
		xticklabels={$1$,$\frac{1}{2}$,$\frac{1}{3}$,$\frac{1}{4}$},
		extra x tick labels={$-1$,$\frac{-1}{2}$,$\frac{-1}{3}\hspace*{1mm}$,$\hspace*{1mm}\frac{-1}{4}$},
		extra x tick style={x tick label style={above, yshift=0.6ex}},
%		axis on top,
%		xtick={1,-1,1/2,-1/2,1/3,-1/3},
	]
		% Limit - step function
		\addplot[line width=0.7mm,domain=-2:0] {-1};
		\addplot[line width=0.7mm,domain=0:2] {1};
		\node[circle,draw=black,line width=0.4mm, fill=white, inner sep=0pt,minimum size=5pt] at (axis cs:0,1) {};
		\node[circle,draw=black,line width=0.4mm, fill=white, inner sep=0pt,minimum size=5pt] at (axis cs:0,-1) {};
		\node[above,white!20!black] at (axis cs:-1.15,-1) {$\boldsymbol{f}$};
		
		% Convergents
		\addplot[domain=-1:1] {x};
		\addplot[dashed,gray] coordinates {(1,0) (1,1)};
		\addplot[dashed,gray] coordinates {(-1,0) (-1,-1)};
		\node[above,white!20!black] at (axis cs:1,1) {$f_1$};
		
		\addplot[domain=-1/2:1/2] {2*x};
		\addplot[dashed,gray] coordinates {(1/2,0) (1/2,1)};
		\addplot[dashed,gray] coordinates {(-1/2,0) (-1/2,-1)};
		\node[above,white!20!black] at (axis cs:1/2,1) {$f_2$};

		\addplot[domain=-1/3:1/3] {3*x};
		\addplot[dashed,gray] coordinates {(1/3,0) (1/3,1)};
		\addplot[dashed,gray] coordinates {(-1/3,0) (-1/3,-1)};
		\node[above,white!20!black] at (axis cs:1/3,1) {$f_3$};
		
		\addplot[domain=-1/4:1/4] {4*x};
		\addplot[dashed,gray] coordinates {(1/4,0) (1/4,1)};
		\addplot[dashed,gray] coordinates {(-1/4,0) (-1/4,-1)};
		\node[above,white!20!black] at (axis cs:1/4,1) {$f_4$};
		
		\node at (axis cs:0.1,0.75) {$\cdots$};
		\node at (axis cs:-0.1,-0.75) {$\cdots$};
		\node[circle,draw=white,line width=0.4mm, fill=black, inner sep=0pt,minimum size=5pt] at (axis cs:0,0) {};
	\end{axis}
	\end{tikzpicture}
	\caption[A sequence of continuous functions converging to a discontinuous function.]{A sequence of continuous functions $f_n$ converging to a discontinuous function $f$.}
	\label{fig:step}
\end{figure}

It is well-known that, while \textit{uniform} limits of continuous functions remain continuous, general pointwise limits don't have to be. A famous example are the functions
\setlength\arraycolsep{2pt}
\begin{equation*}
	f_n(x) = \left\{ \begin{array}{lrcl}
		-1 &\quad x &\leq& \nicefrac{-1}{n} \\
		nx &\quad \abs{x} &<& \nicefrac{1}{n} \\
		1 &\quad x &\geq& \nicefrac{1}{n} \\
	\end{array} \right.
\end{equation*}
which converge pointwise, non-uniformly, to the Heaviside step function (Figure \ref{fig:step}):
\begin{equation*}
	f(x) = \left\{ \begin{array}{lrcl}
		-1 &\quad x< 0 \\
		0 &\quad x=0 \\
		1 &\quad x>0 \\
	\end{array} \right.
\end{equation*}
%$$f_n(x) = \begin{cases}
%	-1 & x \leq -1/n \\
%	nx & \abs{x} < 1/n \\
%	1 & x \geq 1/n \\
%\end{cases}$$

Taking all pointwise limits of continuous functions gives the \textit{Baire 1 functions}. The Baire 1 functions aren't closed under pointwise limits either, so again taking \textit{their} pointwise limits gives the Baire 2 functions.

We can continue this process transfinitely up to $\omega_1$, at which point the Baire hierarchy collapses. In this way, the Baire classes assign a measure of complexity to the Borel functions. Indeed, an equivalent definition of Baire class $\alpha$ is that the preimage of any open set is $\boldsymbol{\Sigma}^0_{\alpha+1}$ in the Borel hierarchy. We will only be concerned with finite Baire classes here.

The Baire classes were introduced by Baire in his PhD thesis \cite{baireFonctionsVariablesReelles1899}, as a natural generalisation of the continuous functions. One motivation for Baire functions is that many functions arising in analysis are not continuous, such as step functions \cite{heavisideElectromagneticTheory1893}, Walsh functions \cite{walshClosedSetNormal1923}, or Dirichlet's function \cite{dirichletConvergenceSeriesTrigonometriques1829}. However, all such ``natural'' functions generally have low Baire class; for example, the derivative of any differentiable function is Baire 1, as are functions arising from Fourier series \cite{kechrisClassificationBaireClass1990}.%, the Baire 1 functions can be characterised as exactly the derivatives of differentiable functions \cite{arg2}.

The Baire class functions have previously been studied with respect to computability \cite{kuyperEffectiveGenericityDifferentiability2014, porterNotesComputableAnalysis2017}. In particular, Kuyper and Terwijn showed a real number $x$ is 1-generic (``random'') if and only if every \textit{effective} Baire 1 functions is continuous at $x$ \cite{kuyperEffectiveGenericityDifferentiability2014}.

Because continuous functions can be specified by countable information, so can Baire 1 functions (countably many continuous functions), and by induction, so can functions of any Baire class. Therefore, we can define Baire functions in second-order arithmetic, coding them using subsets of $\N$.

%\begin{minipage}{0.5\textwidth}
%	In this section, we continue the discussion of the previous section in the new setting of Baire $n$ functions. We prove that for any $n \geq 1$, Cousin's lemma for Baire $n$ functions is equivalent to $\ACA$ over $\RCA$, and this is independent of whether we require the function values to actually exist ($\CLB{'n}$) or not ($\CLB{n}$).
%	
%	On the right is a diagram of implications for various forms of Cousin's lemma (over base system $\RCA$). The implications in {\color{red} red} are proven in [Simpson]; those in {\color{green!70!black} green} were proven in Section \ref{sec:clc}; those in black are trivial (Proposition \ref{prop:clbn-triv}); those in {\color{blue} blue} form the main results of this section. None of the implications can be reversed above the dashed line.
%\end{minipage}
%\hfill
%\begin{minipage}{0.45\textwidth}
%%	\begin{tikzcd}[row sep=large,% column sep = small, %shorten=1mm,
%%	execute at end picture={
%%		\draw[dashed,lightgray] (-3.3,-2.2) -- (3.2,-2.2);
%%	}
%%	]
%%	& \vdots \arrow[d] & \vdots \arrow[d] \\
%%	& \CLB{3} \arrow[r] \arrow[d] & \CLB{'3} \arrow[d] \\
%%	\ACA \arrow[uur,blue,dashed] \arrow[ur,blue,dashed] \arrow[r,blue,dashed] \arrow[dr,blue,dashed] \arrow[dd,bend left=-35,red] & \CLB{2} \arrow[r] \arrow[d] & \CLB{'2} \arrow[d] \\
%%	\SC \arrow[leftrightarrow,u,red] & \CLB{1} \arrow[r] & \CLB{'1} \arrow[ll,bend left=15,blue] \arrow[d] \\
%%	\WKL \arrow[rr,bend left=-15,green!70!black] & \HB \arrow[leftrightarrow,l,red] & \CLB{0} = \CLc \arrow[l,green!70!black] \\
%%	\end{tikzcd}
%\end{minipage}

\begin{definition}\label{defn:baire-fn}
	The following definitions proceed simultaneously and inductively on $n$.
	
	\begin{enumerate}
		\item The \textit{Baire 0 functions} are exactly the (total) continuous functions of Definition \ref{defn:cont-func}.

		\item For each $n \in \N$, a \textit{Baire $n+1$ function} $f\colon \UI \to \R$ is a countable sequence $(f_n)_{n \in \N}$ of Baire $n$ functions $\UI \to \R$ which is pointwise Cauchy. That is, for each $x \in \UI$ %, there is some $y \in \R$ such that for all 
		and $\varepsilon \in \Q^+$, %there is some $N \in \N$ such that for all $m \geq N$,
		$\abs{f_m(x) - f_n(x)} \leq \varepsilon$ for sufficiently large $m,n$.

		\item Given two Baire $n+1$ functions $f = (f_n)_{n \in \N}$, $g = (g_n)_{n \in \N}\colon \UI \to \R$, a point $x \in \UI$, and $\varepsilon \in \Q^+$, we say $\abs{f(x) - g(x)} \leq \varepsilon$ if for any $\delta \in \Q^+$, we have $\abs{f_n(x) - g_n(x)} \leq \varepsilon+\delta$ for sufficiently large $n$.% there exists $N \in \N$ such that for all $n \geq N$, {\color{red} u w0t m9}% $f_n(x)$, $g_n(x)$ are defined and $\abs{f_n(x) - g_n(x)} \leq \varepsilon$.

%		\item For a given (partial) Baire $n+1$ function $f\colon \UI \to \R$ and a point $x \in \UI$, $x$ is \textit{in the domain of $f$} (symbolically, $x \in \dom f$) if there exists $y \in \R$ such that for all $\varepsilon \in \Q^+$, there is some $N \in \N$ such that for all $m \geq N$, $f_m(x)$ is defined and $\abs{f_m(x) - y} < \varepsilon$.

%		\item If $x \in \dom f$, such a $y$ (which must be unique up to $=_\R$-equivalence) is denoted $f(x)$. We then say $f(x)$ \textit{exists} or \textit{is defined}.

%		\item We say $f$ is \textit{total} if all $x \in \UI$ have $x \in \dom f$.
	\end{enumerate}
\end{definition}

\noindent Here, we have only required our Baire functions to be pointwise \textit{Cauchy}. This means, in weak subsystems such as $\RCA$ and $\WKL$ which can't prove the completeness of the reals, the function value $f(x)$ may not actually \textit{exist}. We could have made the stronger requirement that $(f_n)$ is pointwise \textit{convergent}; however, Definition \ref{defn:baire-fn} proves to be the right one for a reverse-mathematical analysis.

% Equivalent definitions of Baire 1: (need to prove equiv)
%  - Pointwise limits of a sequence of continuous functions
%  - Derivatives of continuous functions
%  - Continuous except on a countable/null/meagre set

% Even if partial, can talk about it being a gauge / delta-fine
% Kuyper - continuous at generics
%   Jump function - for each r, f(r) = r'

%\begin{definition}
%	\begin{enumerate}
%		\item $f_n \big( B(a,r) \big) \subseteq B(b,s)$ and $f_n \big( B(a,r) \big) \subseteq B(b',s') \implies B(b,s) \cap B(b',s') \neq \varnothing$;
%		
%		\item $f_n \big( B(a,r) \big) \subseteq B(b,s)$ and $B(a',r') \subseteq B(a,r) \implies f_n \big( B(a,r) \big) \subseteq B(b,s)$;
%		
%		\item $f_n \big( B(a,r) \big) \subseteq B(b,s)$ and $B(b',s') \subseteq B(b,s) \implies f_n \big( B(a,r) \big) \subseteq B(b',s')$;
%		
%		\item $\forall \varepsilon \in \Q^+\ \exists n \in \N\ \forall m \geq n\ {\color{red}\mu(())}$
%	\end{enumerate}
%	where $f_n \big( B(a,r) \big) \subseteq B(b,s)$ abbreviates $(n,a,b,c,d) \in f$, $B(b,s) \cap B(b',s') \neq \varnothing$ abbreviates $(c < d') \land (c' < d)$, and $B(a',r') \subseteq B(a,r)$ abbreviates $(a \leq a') \land (b' \leq b)$.
%\end{definition}

\begin{example}\label{exm:cont-b1}
	Any Baire $n$ function $f$ can be identified with a Baire $n+1$ function $\tilde{f} = (f)_{n \in \N}$. Therefore, the Baire classes are nested: $\mathsf{B}0 \subseteq \mathsf{B}1 \subseteq \mathsf{B}2 \subseteq \cdots$.
\end{example}

Since the function values $f(x)$ may not actually exist, we must take some care when making definitions concerning Baire functions. We can see this already in Definition \ref{defn:baire-fn}, where we had to define what $\abs{f(x) - g(x)} \leq \varepsilon$ means, despite the fact that both values may not exist. These difficulties can generally be overcome with a little caution.

\begin{definition}\label{defn:baire-geq}\
	\begin{enumerate}
		\item The following definition proceeds inductively on $n$. Given a Baire $n+1$ function $f = (f_n)_{n \in \N}\colon \UI \to \R$, a point $x \in \UI$, and a point $y \in \R$, we say $f(x) \geq y$ if for all rational $q < y$, $f_m(x) > q$ for sufficiently large $m$. $f(x) \leq y$ is defined similarly.
		
		\item We say $f(x) > y$ if it is not true that $f(x) \leq y$. %We say $f(x) > y$ if $f(x) \geq y'$ for some real $y' > y$.
		$f(x) < y$ is defined similarly.
		
		\item A Baire $n$ function $\delta\colon \UI \to \R$ is a \textit{gauge} if $\delta(x) > 0$ for all $x \in \UI$. This property will be denoted $\delta\colon \UI \to \R^+$.
	\end{enumerate}
\end{definition}

\begin{definition}
	Let $\delta\colon \UI \to \R^+$ be a Baire $n$ gauge, and $P$ be a tagged partition of size $\ell$. Then, we say $P$ is \textit{$\delta$-fine} if for all $j < \ell$, we have $\delta(t_j) \geq t_j - x_j$ and $\delta(t_j) \geq x_{j+1} - t_j$, in the sense of Definition \ref{defn:baire-geq}.
\end{definition}

Having defined gauges and $\delta$-fine partitions, we are now ready to define Cousin's lemma for Baire $n$ functions:

\begin{definition}
	For each $n \in \N$, let $\CLB{n}$ be the following statement in $\RCA$: every Baire $n$ gauge $\delta\colon \UI \to \R^+$ has a $\delta$-fine partition.
\end{definition}

Note that $\CLc = \CLB{0}$, and for each $m \geq n$, $\CLB{m} \proves \CLB{n}$, since the Baire classes are nested. Combining these with the results of the previous section, immediately we get:

\begin{theorem}[($\RCA$)]
	For each $n \in \N$, $\CLB{n}$ implies $\WKL$.
\end{theorem}

As we will see, for $n \geq 1$, this is far from optimal; $\CLB{n}$ is much stronger than $\WKL$.

\section{\texorpdfstring{$\CLB{1}$}{CLB1} proves \texorpdfstring{$\ACA$}{ACA-0}}

Having completely classified $\CLc = \CLB{0}$, the natural next step would be to study the reverse-mathematical strength of $\CLB{1}$. Our first result about $\CLB{1}$ is a reversal; we show that over $\RCA$, $\CLB{1}$ proves $\ACA$. In other words, to prove $\CLB{1}$ in second-order arithmetic, we need a system \textit{at least} as strong as $\ACA$. The reversal goes through the sequential completeness of $\R$, which is known to be equivalent to $\ACA$.

\begin{definition}
	Let $\SC$ be the following statement in $\RCA$: every Cauchy sequence of real numbers in $\UI$ has a limit.
\end{definition}

\begin{theorem}[($\RCA$)]\label{thm:clb1->aca}
	$\CLB{1}$ implies $\SC$.
\end{theorem}

The idea of the proof is as follows. Supposing we have a Cauchy sequence $(z_n)$ with no limit, we look at the sequence of functions $\delta_n\colon x \mapsto \tfrac{1}{2}\abs{x-z_n}$. This is pointwise Cauchy, hence Baire 1, and it is a gauge since $(z_n)$ has no limit. But $\delta = (\delta_n)$ can't have a $\delta$-fine partition, since no partition $P$ can cover the gap where $\lim z_n$ should be. Here are the details.

\begin{proof}[Proof of Theorem \ref{thm:clb1->aca}]
	By contradiction: suppose there is a Cauchy sequence $(z_n)_{n \in \N} \subseteq [0,1]$ that has no limit in $[0,1]$. For each $n \in \N$, let $\delta_n\colon x \mapsto \tfrac{1}{2}\abs{x-z_n}$, which defines a continuous function by earlier lemmas.
	
	The sequence $(\delta_n)_{n \in \N}$ is pointwise Cauchy; fixing $\varepsilon \in \Q^+$, we have $\abs{z_m - z_n} \leq \varepsilon$ for sufficiently large $m,n$. But $\abs{z_m - z_n} \leq \varepsilon$ implies $\abs{\delta_m(x) - \delta_n(x)} \leq \varepsilon$ by the reverse triangle inequality:
	\begin{equation*}
		\abs{\delta_m(x) - \delta_n(x)}\ =\ \abs{\tfrac{1}{2}\abs{x-z_m} - \tfrac{1}{2}\abs{x-z_n}}\ \leq\ \tfrac{1}{2}\abs{z_n - z_m}\ \leq\ \varepsilon.
	\end{equation*}
	It follows that $\abs{\delta_m(x) - \delta_n(x)} \leq \varepsilon$ for sufficiently large $m,n$, whence $\delta \defeq (\delta_n)_{n \in \N}$ is Baire 1.
	
	Now, because $(z_n)$ does not have a limit, we claim $\delta$ is a gauge. For any $x \in \UI$, since $x$ is not a limit for $(z_n)$, there is some $\varepsilon \in \Q^+$ such that $\delta_n(x) = \abs{x - z_n} \geq \varepsilon$ eventually. Thus $\delta(x) > 0$. % But since $(z_n)$ is Cauchy, pick $N \in \N$ where beyond $N$, all $z_n$ are within $\varepsilon/2$. The triangle inequality now implies that all $n \geq N$ have $\delta_n(x) = \abs{x - z_n} \geq \varepsilon/2$ as required.
	Let $P = \ang{x_j,t_j}$ be a partition of $\UI$; we will show $P$ is not $\delta$-fine.
	
	Fix $j < \ell$. Since $x_j$ is not the limit of $(z_n)$, we must eventually have $\abs{x_j - z_n} \geq \varepsilon$ for some $\varepsilon \in \Q^+$. But since $(z_n)$ is Cauchy, all the terms are eventually within $2\varepsilon$ of each other; from this point on, we must have $z_n < x_j$ for all $n$, or $z_n > x_j$ for all $n$. As $(z_n) \subseteq \UI$, we can't have $z_n < x_0 = 0$ or $z_n > x_\ell = 1$. It follows that there is $j < \ell$ such that eventually $x_j < z_n < x_{j+1}$.
	
	For this $j$, we claim that $\delta(t_j)$ cannot cover $(x_j,x_{j+1})$. By the same argument to the previous paragraph, either $z_n < t_j$ eventually, or $z_n > t_j$ eventually; let us suppose WLOG that $z_n < t_j$. Then, $$\delta_n(t_j)\ =\ \frac{1}{2} \abs{t_j - z_n}\ <\ \abs{t_j - z_n}\ <\ \abs{t_j - x_j}$$
	for sufficiently large $n$, so $\delta(t_j) < t_j - x_j$. Thus, $P$ is not $\delta$-fine.%
%	Take $\varepsilon \defeq \min \left\{ \abs{t_j - x_j}, \abs{x_{j+1} - t_j}: j < \ell \right\}$ to be the minimum distance between points in $P$. Then, $\varepsilon>0$, and since $
\end{proof}

Having demonstrated the reversal $\CLB{1} \proves \ACA$, it seems natural to see if we can get a proof of $\CLB{1}$ in $\ACA$. %Here, we don't rule out that possibility, but give an argument as to why one obvious way to do it will not work.
The most natural way to do this is as follows. A Baire 1 function $\delta$ is a pointwise limit of continuous gauges $\delta_n$, and we have already seen that $\ACA$ (in fact, $\WKL$) can construct a $\delta_n$-fine partition for each $n$. Therefore, one might expect that there would be some inductive way to combine the partitions for $\delta_n$ into a partition for $\delta$.

The following proposition suggests that this is not possible. Looking back to the proof of $\CLc$, for each continuous gauge $\delta_n$, we in fact constructed a \textit{dyadic} $\delta_n$-fine partition---one whose partition points and tag points were all \textit{dyadic rationals}, i.e.\ of the form $j/2^n$. If there were a way to combine these partitions into one for $\delta$, it would follow that every Baire 1 gauge has a dyadic partition. Now, we present a Baire 1 gauge with no dyadic partition.

%Since $\delta$ above is Baire 1, it is a limit of continuous gauges $\delta_n$, each of which we saw has a dyadic  However, $\delta$ itself does not have a dyadic $\delta$-fine partition. This suggests that in general, we cannot hope to combine the partitions for $\delta_n$ into a partition for $\delta$.

\begin{proposition}
	There is a Baire 1 gauge $\delta\colon \UI \to \R$ with no dyadic $\delta$-fine partition.
\end{proposition}

\begin{proof}
	The desired gauge is
	\begin{equation*}
		\delta(x) = \left\{ \begin{array}{lrcl}
			\frac{1}{2} \abs{x - \frac{1}{3}} &\qquad  x \neq \frac{1}{3} \\[2mm]
			1 &\qquad  x = \frac{1}{3}\\
		\end{array} \right.
	\end{equation*}
	It is the pointwise limit of the following sequence $\delta_n\colon \UI \to \R$, where $n \geq 3$:
	\begin{equation*}
		\delta_n(x) = \left\{ \begin{array}{lrcl}
			\frac{1}{2} \abs{x - \frac{1}{3}} &\qquad \abs{x-\frac{1}{3}} &\geq& \frac{1}{n} \\[2mm]
			1 - \left( n - \frac{1}{2} \right) \abs{x - \frac{1}{3}} &\qquad \abs{x-\frac{1}{3}} &<& \frac{1}{n} \\
		\end{array} \right.
	\end{equation*}
	By Lemma \ref{lem:piecewise-cont}, each $\delta_n$ is continuous, so $\delta$ is Baire 1. Now, any $\delta$-fine partition $P$ must have some $t_i = \tfrac{1}{3}$, because otherwise $\delta(t_i) = \frac{1}{2} \abs{t_i - \frac{1}{3}} < \abs{t_i - \frac{1}{3}}$. Then $\tfrac{1}{3} \notin B \big( t_i,\delta(t_i) \big)$, so $(x_j,x_{j+1}) \nsubseteq B \big( t_j,\delta(t_j) \big)$ for the $j < \ell$ such that $\frac{1}{3} \in (x_j,x_{j+1})$. It follows that a $\delta$-fine partition $P$ cannot be dyadic.
\end{proof}

%Show $\CLB{n}$ fails in $\ARITH$ - same proof that $\CLB{n} \proves \ACA^+$

\section{\texorpdfstring{$\CLB{2}$}{CLB2} fails in \texorpdfstring{$\ACA$}{ACA-0}}

In attempting to show the implication in Theorem \ref{thm:clb1->aca} is strict, we tried to construct an arithmetical Baire 1 gauge with no arithmetical partition. This would show that $\CLB{1}$ fails in the standard model $\ARITH$ of $\ACA$, and thus that $\ACA \nproves \CLB{1}$.

We have not yet been able to construct such a Baire 1 gauge, but we were able to construct an arithmetical \textit{Baire 2 gauge} with no arithmetical partition. We present the construction in this section, showing that $\ACA \nproves \CLB{2}$. In the next section, we will slightly generalise the ideas of the proof to show $\CLB{2} \proves \ATR$.

Instead of working in $\UI$ as before, we will actually construct our gauges in Cantor space $\CS$. The definitions are generally analogous to the $\UI$ case:
\begin{itemize}
	\item A \textit{gauge} on $\CS$ is a positive real-valued function $\delta\colon \CS \to \R^+$.
	
	\item $\CS$ is a metric space under the distance function\vspace*{-2mm}
	\begin{equation*}
		d(X,Y) = \begin{cases}
			0 & X = Y \\
			2^{-n} & n \text{ least such that } X_n \neq Y_n
		\end{cases}
	\end{equation*}\vspace*{-3mm}
	
	\item A \textit{$\delta$-fine partition} is a finite set $P \subseteq \CS$ such that $\big\{ B \big( X,\delta(X) \big)\!: X \in P \big\}$ is an open cover of $\CS$.
	
	\item \textit{Cousin's lemma} says that every gauge $\delta\colon \CS \to \R^+$ has a $\delta$-fine partition $P \subseteq \CS$.
	
	\item $\CS$ has a countable basis of basic open sets $[\sigma] = \{ X \in \CS: \sigma \preceq X \}$ for each $\sigma \in \CSf$. We can code basic open sets $[\sigma]$ by natural numbers, then give the same definitions of open sets, continuous/Baire functions $f\colon \CS \to \R$, etc... in second-order arithmetic.
\end{itemize}

\noindent There is a well-known embedding $g\colon \CS \to \UI$ defined by $$g(X)\ =\ \sum_{n=0}^\infty \dfrac{2X_n}{3^{n+1}}$$
The range of $g$ is the Cantor middle-thirds set $C \subseteq \UI$, and topologically, $g$ is a homeomorphism $\CS \to C$. Via the embedding $g$, we can map any Baire $n$ gauge $\delta\colon \UI \to \R^+$ to a Baire $n$ gauge $\bar{\delta}\colon \CS \to \R^+$, and vice versa. Furthermore, we can do this in such a way to preserve covering, i.e. $d(X,Y) < \bar{\delta}(Z)$ if and only if $\abs{g(X) - g(Y)} < \delta(g(Z))$.

There is one difficulty to contend with: when going from $\CS$ to $\UI$, this correspondence only defines a gauge $\delta$ on $C \subseteq \UI$. However, since $C \subseteq \UI$ is closed, any point $x \notin C$ has positive distance $r$ to $C$, so we can just choose $\delta(x)<r$. This ensures the aforementioned covering property is preserved, and we can make this choice in a Baire 1 way. It follows that, for $n \geq 1$, $\CLB{n}$ for gauges on $\UI$ and $\CLB{n}$ for gauges on $\CS$ are equivalent.

Before we see the proofs, we need to introduce a bit more computability. In Section \ref{sec:comp}, we focused on \textit{absolute} computability---the existence of an algorithm to solve some problem (e.g.\ membership in a set $A \subseteq \N$). However, many natural problems are not computable in this sense. This leads us to a more general notion of \textit{relative} computability.

The idea is we allow our computations access to an \textit{oracle}---a (noncomputable) set $A$. While performing our algorithm, we are allowed to query $A$ at any point, and ask if it contains some element or not. We say $B$ is \textit{Turing reducible} to $A$ ($B \leq_\mathrm{T} A$) if there is an algorithm which can compute $B$, with $A$ as an oracle. $\leq_\mathrm{T}$ is a preorder on $\Pow(\N)$, forming a hierarchy known as the \textit{Turing degrees}. Intuitively, one should think of $B \leq_\mathrm{T} A$ as meaning that $A$ has more computational power than $B$, or that $A$ is \textit{less computable} than $B$.

The \textit{Turing jump} is an operation assigning to every set $A \subseteq \N$ a set $A' >_\mathrm{T} A$ which is strictly higher in the Turing degrees, i.e.\ $A'$ (read ``$A$-jump'') is less computable than $A$. We can iterate this operation, getting a sequence $A <_\mathrm{T} A' <_\mathrm{T} A'' <_\mathrm{T} A^{(3)} <_\mathrm{T} A^{(4)} <_\mathrm{T} \cdots$. This is a countable sequence of countable sets, so we can combine them all into a single countable set $A^{(\omega)} = \big\{ (e,n): e \in A^{(n)} \big\}$, called the \textit{arithmetic jump} or \textit{$\omega$-jump} of $A$. The Turing jumps of the empty set $\varnothing$ provide some useful milestones in the arithmetical hierarchy:

\begin{proposition}\label{prop:zero-jump}\
	\begin{itemize}
		\item $\varnothing^{(n)}$ is strictly $\Sigma^0_n$ in the arithmetical hierarchy.
		
		\item $\varnothing^{(\omega)}$ is nonarithmetical.
	\end{itemize}
\end{proposition}

\noindent We are now ready to prove the first result: that $\CLB{2}$ fails in $\ACA$.

\begin{theorem}\label{thm:aca-clb2}
	$\ACA$ does not imply $\CLB{2}$.
\end{theorem}

The idea is similar to Theorem \ref{thm:rca-clc}, and proceeds as follows. We work in the standard model $\ARITH$ of $\ACA$. In Proposition \ref{prop:zero-jump}, we saw that $X = \varnothing^{(\omega)}$ is nonarithmetical; hence, it does not exist in $\ARITH$. However, the singleton set $\{ X \} \subseteq \CS$ is \textit{effectively $G_\delta$}; i.e.\ it can be written $\{ X \} = \bigcap_{n \in \N} O_n$ for a computable sequence $O_n \subseteq \CS$ of open sets \cite[Prop XII.2.19]{odifreddiClassicalRecursionTheory1999}. Implicitly using this result, we can construct (arithmetically) a Baire 2 function $\delta\colon \CS \to \R^{\geq 0}$ such that:
\begin{enumerate}
	\item $\delta(Y) = 0 \iff Y = X$;
	\item For all $Y \neq X \in \CS$, $\delta(Y) < d(Y,X)$.
\end{enumerate}
By property (i), $\ARITH$ believes that $\delta$ is a gauge, since it is positive everywhere except $X$; in particular, at every arithmetical point. However, by property (ii), no point $Y \neq X$ can $\delta$-cover $X$, and hence there is no arithmetical $\delta$-fine partition.

In what follows, we freely identify an element $Y \in \CS$ of Cantor space with the set $\{ n \in \N: Y_n = 1 \}$.

\begin{proof}[Proof of Theorem \ref{thm:aca-clb2}]
	Let $X = \varnothing^{(\omega)}$. We define a function $\delta\colon \CS \to \R^{\geq 0}$ as follows. For every $Y \in \CS$, we want to find a position $k$ where $Y_k \neq X_k$, and define $\delta(Y)$ accordingly.
	
	First, we consider the columns $Y^{[n]} = \{ e: (e,n) \in Y \}$, and try to find the least column where $Y^{[n]} \neq X^{[n]}$. Ask if $Y^{[0]} = \varnothing$; then if $Y^{[1]} = \big( Y^{[0]} \big)'$, then if $Y^{[2]} = \big( Y^{[1]} \big)'$, etc. The desired column is the first one where the answer is ``no''. Having found this column $Y^{[n]}$, we simply search along it to find $k$, the first point of difference from $\big( Y^{[n-1]} \big)'$ or $\varnothing$. Then, let $\delta(Y) = 2^{-k-1}$.
	
	Since $Y^{[n]}$ is the \textit{least} column of difference, by induction we have $Y^{[m]} = \varnothing^{[m]}$ for all $m < n$. Thus, we indeed have $Y_k \neq X_k$. So for any $Y \neq X$, we \textit{will} find a point of difference, whence property (i) above holds. Also, $d(Y,X) \geq 2^{-k} > 2^{-k-1} = \delta(Y)$, giving property (ii).
	
	Now, the question $A = B'$ can be answered by the double-jumps $A''$ and $B''$ \cite[Prop XII.2.19]{odifreddiClassicalRecursionTheory1999}, so it follows that $\delta$ is computable from the double-jump function $X \mapsto X''$. Since the double-jump is Baire 2 \cite{porterNotesComputableAnalysis2017}, and computable reductions are always continuous \cite{pour-elComputabilityAnalysisPhysics1989}, this implies that $\delta$ is Baire 2.
	
	So, we have constructed a Baire 2 gauge $\delta$ in $\ARITH$; now we claim that it has no $\delta$-fine partition. The argument is as sketched---suppose $P \subseteq \CS$ is a finite subset. Property (ii) implies that for all $Y \in P$, $\delta(Y) < d(Y,X)$, and so $X \notin B \big( Y,\delta(Y) \big)$. Thus, $P$ is not a $\delta$-fine partition.
\end{proof}

Viewing this proof in a different light, it can be construed as a proof from $\CLB{2}$ that $\varnothing^{(\omega)}$ exists. The argument is by contradiction: if $\varnothing^{(\omega)}$ doesn't exist, then the function $\delta$ constructed in the proof of Theorem \ref{thm:aca-clb2} \textit{is} a gauge. As we essentially argued, any $\delta$-fine partition then must include $X = \varnothing^{(\omega)}$, hence this set exists.

The same argument works replacing $\varnothing$ by \textit{any} arithmetical set $A$. As a corollary, then, we see that $\CLB{2}$ implies the stronger system $\ACA^+$, consisting of $\ACA$ plus the assertion that the $\omega$-jump of any set exists. The system $\ACA^+$ has arisen previously in reverse mathematics, first with the work of Blass, Hirst and Simpson in combinatorics and topological dynamics \cite{blassLogicalAnalysisTheorems1987}. Later, it surfaced in Shore's work on Boolean algebras \cite{shoreInvariantsBooleanAlgebras2005}, and Downey and Kach's work on Euclidean domains \cite{downeyEuclideanFunctionsComputable2011}.

\section{\texorpdfstring{$\CLB{2}$}{CLB2} proves \texorpdfstring{$\ATR$}{ATR-0}}

% is zero \textit{only} at the point $X$.% Then, $\ARITH$ will believe 

Using a similar idea to the proof of Theorem \ref{thm:aca-clb2}, we can show that $\CLB{2}$ implies $\ATR$. We don't prove the existence of $\Theta^\alpha(X)$ for \textit{every} arithmetical operator $\Theta$, countable ordinal $\alpha$, and set $X \subseteq N$ (Definition \ref{defn:atr}). Instead, it is enough to show this for $\Theta = \mathrm{TJ}\colon A \mapsto A'$, the Turing jump operator \cite[Thm VIII.3.15]{simpsonSubsystemsSecondOrder2009}. This is because the Turing jump $A \mapsto A'$ is a \textit{universal $\Sigma^0_1$ operator} \cite[Defn VIII.1.9]{simpsonSubsystemsSecondOrder2009}, and so any arithmetical operator can be expressed using a finite number of Turing jumps. Otherwise, the proof proceeds along similar lines; the details follow.

%We use the fact that  and so $\ATR$ is equivalent to the existence of the $\alpha$th Turing jump $A^{(\alpha)}$ for every set $A \subseteq \N$ and countable ordinal $\alpha$  The details follow.

\begin{theorem}
	$\CLB{2}$ implies $\ATR$.
\end{theorem}

\begin{proof}
	Let $\M$ be a model of $\CLB{2}$, so in particular, $\M \models \CLB{1}$ and hence $\ACA$ by Theorem \ref{thm:clb1->aca}. By contradiction, suppose that $\ATR$ fails in $\M$. Then, there is a set $A \subseteq \N$ and countable ordinal $\alpha$ such that the $\alpha$th Turing jump $X = A^{(\alpha)}$ doesn't exist in $\M$.
	
	Using the same construction as the proof of Theorem \ref{thm:aca-clb2}, for each $Y \in \CS$, we can find the least column $\beta = \beta_Y < \alpha$ such that $Y^{[\beta]} \neq \big( Y^{[<\beta]} \big)'$, then find the first point of difference $k = k_Y$. Furthermore, the function $\delta(Y) = 2^{-k-1}$ is Baire 2, as before. %This $\delta$ has the same properties as earlier, so
	
	$\CLB{2}$ gives a $\delta$-fine partition $P$; let $\beta^* = \max\{ \beta_Y: Y \in P \}$. Then, $Z = X^{[\beta^*]}$ exists in $\M$ by arithmetical comprehension; we claim $Z$ is not covered by $P$. For any $Y \in P$, $\beta_Y$ is defined so that $Y$ disagrees with $Z$ on column $\beta_Y \leq \beta^*$. Hence, $k_Y$ is a point of disagreement between $Y$ and $Z$, so the \textit{first} point of disagreement is at most $k_Y$. It follows that $\delta(Y) < d(Y,Z)$, thus $P$ is not $\delta$-fine; contradiction.
%
%Let $Y_* \in P$ be such that $\beta_{Y_*} = \beta_*$. By induction, we have $Y_*^{[<\beta_*]} = X^{[<\beta_*]} \eqdef Z$, but $Y_*^{[\beta_*]} \neq Z'$. But by $\ACA$, $Z'$ exists in $\M$
%	
%	
	 %there is a Baire 2 function $f\colon \CS \to \alpha$ such that for all $Y \neq X \in 2^\omega$, $f(Y)$ picks out a point of difference from $\Theta^\beta(X)$, i.e. $$Y_{f(Y)}\ \neq\ \big( \Theta^\beta(X) \big)_{f(Y)}$$ By definition of $\beta$, $\Theta^\gamma(X)$ exists in $\M$ for all $\gamma < \beta$; hence, the function $f$ exists in $\M$.
\end{proof}

\chapter{Conclusion}\label{C:con}

\newcommand{\CLK}{\mathsf{CL}_\K}
\renewcommand{\thetheorem}{\arabic{chapter}.\arabic{theorem}}

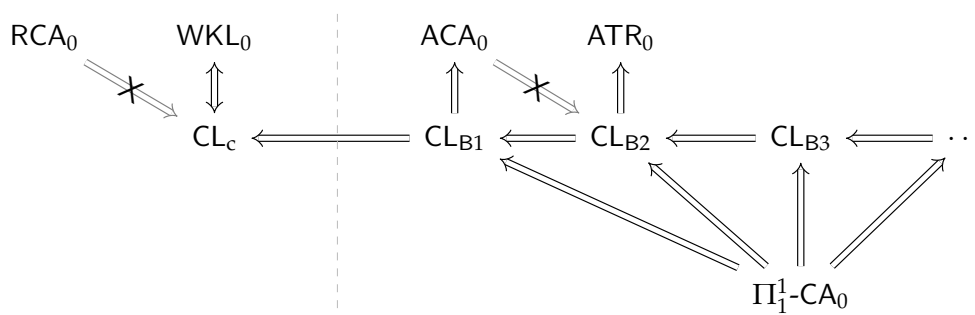
\begin{figure}
	\centering
	\begin{tikzcd}[execute at end picture={\draw[dashed,lightgray] (-2.1,-1.5) -- (-2.1,2.5);},arrows=Rightarrow]
		\RCA \arrow[dr,"\text{\large\sffamily X}"{anchor=center,sloped,black},gray] & \WKL \arrow[d,Leftrightarrow] && \ACA \arrow[dr,"\text{\large\sffamily X}"{anchor=center,sloped,black},gray] & \ATR \\
		& \CLc && \CLB{1} \arrow[ll] \arrow[u] & \CLB{2} \arrow[l] \arrow[u] & \CLB{3} \arrow[l] & \cdots \arrow[l] \\
		\\
		& & & & & \PiCA \arrow[lluu] \arrow[luu] \arrow[uu] \arrow[ruu] \\
	\end{tikzcd}
	\caption{A graphical summary of our contributions to the reverse mathematics zoo.}
\end{figure}

In this report, we introduced reverse mathematics and Cousin's lemma, and then began a reverse-mathematical analysis of Cousin's lemma for various classes of functions. We have established many original results in this direction: here is the summary of our knowledge so far.

\begin{theorem}[(summary of results)]
	All implications are over $\RCA$.
	\begin{enumerate}
		\item Cousin's lemma for continuous functions is equivalent to $\WKL$.
		
		\item Cousin's lemma for Baire 1 functions is provable in $\PiCA$, and it implies $\ACA$.
		
		\item For $n \geq 2$, Cousin's lemma for Baire $n$ functions is provable in $\PiCA$, and it implies $\ATR$.
	\end{enumerate}
\end{theorem}

\noindent Cousin's lemma for continuous functions, $\CLc$, is the only theorem for which we have been able to \textit{completely} determine the axiomatic strength. Naturally, there is further work to be done on classifying Cousin's lemma for Baire $n$ functions. We are still most interested in the case $n=1$; therefore, our main open question is:

\begin{question}
	Where does Cousin's lemma for Baire 1 functions, $\CLB{1}$, fall in the reverse-mathematical hierarchy?
\end{question}

Here is a heuristic reason to believe $\CLB{1}$ implies $\PiCA$, and is thus equivalent to it. Recall an alternative characterisation of Baire 1 functions is that the preimage of any open set is $\boldsymbol{\Sigma}^0_2$ in the Borel hierarchy. Similarly, \textit{effectively} Baire 1 functions can be characterised as those where the preimage of any lightface $\Sigma^0_1$ class is lightface $\Sigma^0_2$.

Now, the only proof of Cousin's lemma we currently know is the ``proof'' of Theorem \ref{thm:cl-pica}, and variations thereof (Theorem \ref{thm:wkl->clc}). To make this proof work for a function $f$, we need to decide if $f^{-1}\big( (2^{-n},\infty) \big)$ is empty or not. For continuous functions, this set is $\Sigma^0_1$, i.e. open, so it is enough to check if it contains any dyadic rational points (which is what we did in the proof of Theorem \ref{thm:wkl->clc}).

For Baire 1 functions, this set is $\Sigma^0_2$, as mentioned. Unfortunately, there is no easy way to determine whether an arbitrary $\Sigma^0_2$ set is empty or not; this problem is $\Pi^1_1$-hard in general. This means that $\Pi^1_1$ comprehension is \textit{required} to make the proof work for Baire 1 functions. So, this suggests that $\CLB{1}$ is equivalent to $\PiCA$, unless there is a smarter way to prove Cousin's lemma (and we don't believe there is).

If it turns out that $\CLB{1}$ and $\PiCA$ are equivalent, this would also imply the equivalence of $\PiCA$ to Cousin's lemma for any class of functions containing the Baire 1 functions. Otherwise, there are many more classes of functions $\K$ for which one could explore the strength of Cousin's lemma. The general question is thus:

\begin{question}
	For a specified class of functions $\K$ definable in second-order arithmetic, where does Cousin's lemma for functions in $\K$, $\CLK$, fall in the reverse-mathematical hierarchy?
\end{question}

Theorem \ref{thm:cl-pica} shows that $\PiCA$ proves $\CLK$ for any class of functions $\K$ definable in second-order arithmetic. Presumably, for large enough $\K$, $\CLK$ becomes equivalent to $\PiCA$; it would be interesting to know where exactly this threshold is. Here are some other classes $\K$ for which one could study the strength of $\CLK$:

\begin{itemize}
	\item Of course, the Baire $n$ functions for $n \geq 2$.
	
	\item On that note, we mentioned that the Baire hierarchy can be iterated transfinitely, so one could equally study the Baire class $\alpha$ functions, for $\omega \leq \alpha < \omega_1$. Defining these in second-order arithmetic can be quite messy, but it is possible.
	
	\item The Borel functions, as the limit of all the Baire classes. There has been some study into Borel sets and functions in reverse mathematics \cite{simpsonSubsystemsSecondOrder2009}.
	
	\item One could look at the \textit{strong} Baire classes $\mathsf{B}'\alpha$, where the sequence is required to be pointwise \textit{convergent}, rather than just pointwise Cauchy. We would hope for \textit{stability} here, i.e.\ $\CLB{'\alpha} \equiv \CLB{\alpha}$, but this is not immediately clear.
	
	\item The Fine continuous functions, which are those continuous with respect to the metric introduced by Fine \cite{fineWalshFunctions1949}. These fall strictly between continuous and Baire 1 functions, and have been studied with respect to computability \cite{moriComputabilitiesFinecontinuousFunctions2001, moriComputabilityWalshFunctions2002, brattkaNotesFineComputability2002}.
\end{itemize}

\noindent Finally, we have only studied a single theorem about gauge integration in this paper. There is a whole theory of gauge integration, with many results waiting to be analysed reverse-mathematically. For those so inclined, the following results would be interesting to study:
\begin{itemize}
	\item The equivalence between different characterisations of the gauge integral, as given by Denjoy, Perron, Luzin and others;
	
	\item Hake's theorem: $\displaystyle \int_a^bf(x)\dif x=\lim_{c\to b^-}\int_a^cf(x)\dif x$
	
	\item If $f\colon \UI \to \R$ is bounded, $f$ is gauge integrable if and only if it is Lebesgue integrable;
	
	\item Basic facts about the gauge integral, such as convergence properties.
\end{itemize}

%%%%%%%%%%%%%%%%%%%%%%%%%%%%%%%%%%%%%%%%%%%%%%%%%%%%%%%

\backmatter

%%%%%%%%%%%%%%%%%%%%%%%%%%%%%%%%%%%%%%%%%%%%%%%%%%%%%%%

%\printindex
\phantomsection
\addcontentsline{toc}{chapter}{Bibliography}
\printbibliography

@thesis{baireFonctionsVariablesReelles1899,
  title = {Sur Les Fonctions de Variables Réelles},
  author = {Baire, René-Louis},
  date = {1899},
  institution = {{École Normale Supérieure}},
  location = {{Paris}},
  type = {PhD thesis}
}

@book{birkhoffLatticeTheory1940,
  title = {Lattice {{Theory}}},
  author = {Birkhoff, Garrett},
  date = {1940},
  publisher = {{American Mathematical Society}},
  location = {{Providence, RI}},
  isbn = {978-0-8218-1025-5 978-1-4704-3173-0},
  langid = {english},
  number = {25},
  pagetotal = {431},
  series = {Colloquium {{Publications}}}
}

@book{birkhoffSourceBookClassical1973,
  title = {A {{Source Book}} in {{Classical Analysis}}},
  author = {Birkhoff, Garrett},
  date = {1973},
  publisher = {{Harvard University Press}},
  eprint = {djA8yQEACAAJ},
  eprinttype = {googlebooks},
  isbn = {978-0-674-82245-0},
  langid = {english},
  pagetotal = {470}
}

@article{blassLogicalAnalysisTheorems1987,
  title = {Logical Analysis of Some Theorems of Combinatorics and Topological Dynamics},
  author = {Blass, Andreas R. and Hirst, Jeffry L. and Simpson, Stephen G.},
  date = {1987},
  journaltitle = {Contemporary Mathematics},
  volume = {65},
  pages = {125--156}
}

@inbook{bolyaiAppendixScientiamSpatii1832,
  title = {Appendix Scientiam Spatii Absolute Veram Exhibens},
  booktitle = {Tentamen Juventutem Studiosam in Elementa Matheseos Purae},
  author = {Bolyai, János},
  date = {1832},
  publisher = {{Kali}},
  bookauthor = {Bolyai, Farkas}
}

@book{bolzanoReinAnalytischerBeweis1817,
  title = {Rein Analytischer {{Beweis}} Des {{Lehrsatzes}}, Dass Zwischen Je Zwey {{Wert}}-Hen, Die Ein Entgegengesetztes {{Resultat}} Gewähren, Wenigstens Eine Reelle {{Wurzel}} Der {{Gleichung}} Liege},
  author = {Bolzano, Bernard},
  date = {1817},
  publisher = {{Gottlieb Haase}},
  location = {{Prague}}
}

@book{booleInvestigationLawsThought1854,
  title = {An {{Investigation}} of the {{Laws}} of {{Thought}}},
  shorttitle = {An {{Investigation}} of the {{Laws}} of {{Thought}}},
  author = {Boole, George},
  date = {1854},
  publisher = {{Walton and Maberly}},
  location = {{London}},
  eprint = {SWgLVT0otY8C},
  eprinttype = {googlebooks},
  langid = {english},
  pagetotal = {688}
}

@article{brattkaNotesFineComputability2002,
  title = {Some Notes on {{Fine}} Computability},
  author = {Brattka, Vasco},
  date = {2002},
  journaltitle = {Journal of Universal Computer Science},
  volume = {8},
  pages = {382--395},
  langid = {english},
  number = {3}
}

@article{brownProofLebesgueCondition1936,
  title = {A Proof of the {{Lebesgue}} Condition for {{Riemann}} Integrability},
  author = {Brown, A. B.},
  date = {1936},
  journaltitle = {The American Mathematical Monthly},
  volume = {43},
  pages = {396--398},
  publisher = {{Mathematical Association of America}},
  issn = {0002-9890},
  doi = {10.2307/2301737},
  eprint = {2301737},
  eprinttype = {jstor},
  number = {7}
}

@article{brownWhichSetExistence1986,
  title = {Which Set Existence Axioms Are Needed to Prove the Separable {{Hahn}}-{{Banach}} Theorem?},
  author = {Brown, Douglas K. and Simpson, Stephen G.},
  date = {1986-01-01},
  journaltitle = {Annals of Pure and Applied Logic},
  shortjournal = {Annals of Pure and Applied Logic},
  volume = {31},
  pages = {123--144},
  issn = {0168-0072},
  doi = {10.1016/0168-0072(86)90066-7},
  url = {http://www.sciencedirect.com/science/article/pii/0168007286900667},
  urldate = {2020-10-27},
  langid = {english}
}

@article{cantorBeitragZurMannigfaltigkeitslehre1877,
  title = {Ein Beitrag zur Mannigfaltigkeitslehre},
  author = {Cantor, Georg},
  date = {1877},
  journaltitle = {Journal für die reine und angewandte Mathematik},
  volume = {84},
  pages = {242--258},
  issn = {0075-4102; 1435-5345/e},
  url = {https://eudml.org/doc/148353},
  urldate = {2020-10-27},
  langid = {und}
}

@article{cantorUeberEigenschaftInbegriffs1874,
  title = {Ueber Eine {{Eigenschaft}} Des {{Inbegriffs}} Aller Reellen Algebraischen {{Zahlen}}},
  author = {Cantor, Georg},
  date = {1874-01-01},
  journaltitle = {Journal für die reine und angewandte Mathematik},
  volume = {77},
  pages = {258--262},
  publisher = {{De Gruyter}},
  issn = {0075-4102, 1435-5345},
  doi = {10.1515/crll.1874.77.258},
  url = {https://www.degruyter.com/view/journals/crll/1874/77/article-p258.xml},
  urldate = {2020-10-25}
}

@article{cantorUeberUnendlicheLineare1883,
  title = {Ueber unendliche, lineare Punktmannichfaltigkeiten},
  author = {Cantor, Georg},
  date = {1883-12-01},
  journaltitle = {Mathematische Annalen},
  shortjournal = {Math. Ann.},
  volume = {21},
  pages = {545--591},
  issn = {1432-1807},
  doi = {10.1007/BF01446819},
  url = {https://doi.org/10.1007/BF01446819},
  urldate = {2020-10-25},
  langid = {german},
  number = {4}
}

@book{cauchyCoursAnalyseEcole1821,
  title = {Cours d'Analyse de l'Ecole Royale Polytechnique},
  author = {Cauchy, Augustin-Louis},
  date = {1821},
  publisher = {{L'Imprim-erie Royale}},
  eprint = {_mYVAAAAQAAJ},
  eprinttype = {googlebooks},
  langid = {french},
  pagetotal = {604}
}

@book{changModelTheory1990,
  title = {Model {{Theory}}},
  author = {Chang, C. C. and Keisler, H. J.},
  date = {1990-06-12},
  edition = {3},
  publisher = {{North Holland}},
  location = {{Amsterdam}},
  eprint = {uiHq0EmaFp0C},
  eprinttype = {googlebooks},
  isbn = {978-0-08-088007-5},
  langid = {english},
  number = {73},
  pagetotal = {667},
  series = {Studies in {{Logic}} and the {{Foundations}} of {{Mathematics}}}
}

@article{conidisComplexityRadicalsNoncommutative2009,
  title = {On the Complexity of Radicals in Noncommutative Rings},
  author = {Conidis, Chris J.},
  date = {2009-11},
  journaltitle = {Journal of Algebra},
  shortjournal = {Journal of Algebra},
  volume = {322},
  pages = {3670--3680},
  issn = {00218693},
  doi = {10.1016/j.jalgebra.2009.07.039},
  url = {https://linkinghub.elsevier.com/retrieve/pii/S0021869309004712},
  urldate = {2020-09-25},
  langid = {english},
  number = {10}
}

@article{cousinFonctionsVariablesComplexes1895,
  title = {Sur les fonctions de $n$ variables complexes},
  author = {Cousin, Pierre},
  date = {1895},
  journaltitle = {Acta Mathematica},
  volume = {19},
  pages = {1--61},
  number = {1}
}

@book{davisComputabilityUnsolvability1958,
  title = {Computability \& {{Unsolvability}}},
  author = {Davis, Martin},
  date = {1958},
  publisher = {{McGraw-Hill}},
  eprint = {85cEAQAAIAAJ},
  eprinttype = {googlebooks},
  langid = {english},
  pagetotal = {248}
}

@book{demorganFormalLogic1847,
  title = {Formal {{Logic}}},
  shorttitle = {Formal {{Logic}}},
  author = {De Morgan, Augustus},
  date = {1847},
  publisher = {{Taylor and Walton}},
  location = {{London}},
  eprint = {HscAAAAAMAAJ},
  eprinttype = {googlebooks},
  langid = {english},
  pagetotal = {376}
}

@article{denjoyExtensionIntegraleLebesgue1912,
  title = {Une Extension de l'intégrale de {{Lebesgue}}},
  author = {Denjoy, Arnaud},
  date = {1912},
  journaltitle = {Comptes rendus de l'Aca-démie des Sciences},
  volume = {134},
  pages = {859--862}
}

@article{dirichletConvergenceSeriesTrigonometriques1829,
  title = {Sur La Convergence Des Séries Trigonométriques Qui Servent à Représenter Une Fonction Arbitraire Entre Des Limites Données},
  author = {Dirichlet, P. G. L.},
  date = {1829},
  journaltitle = {Journal für die reine und angewandte Mathematik},
  volume = {4},
  pages = {157--169}
}

@article{downeyEuclideanFunctionsComputable2011,
  title = {Euclidean Functions of Computable {{Euclidean}} Domains},
  author = {Downey, Rodney G. and Kach, Asher M.},
  date = {2011-04},
  journaltitle = {Notre Dame Journal of Formal Logic},
  shortjournal = {Notre Dame J. Formal Logic},
  volume = {52},
  pages = {163--172},
  issn = {0029-4527},
  doi = {10.1215/00294527-1306172},
  url = {http://projecteuclid.org/euclid.ndjfl/1303995712},
  urldate = {2020-09-25},
  langid = {english},
  number = {2}
}

@book{evesIntroductionHistoryMathematics1969,
  title = {An {{Introduction}} to the {{History}} of {{Mathematics}}},
  author = {Eves, Howard},
  date = {1969},
  edition = {3},
  publisher = {{Holt, Rinehart and Winston}},
  eprint = {LIsuAAAAIAAJ},
  eprinttype = {googlebooks},
  langid = {english},
  pagetotal = {488}
}

@article{fineWalshFunctions1949,
  title = {On the {{Walsh}} Functions},
  author = {Fine, N. J.},
  date = {1949-03-01},
  journaltitle = {Transactions of the American Mathematical Society},
  shortjournal = {Trans. Amer. Math. Soc.},
  volume = {65},
  pages = {372--372},
  issn = {0002-9947},
  doi = {10.1090/S0002-9947-1949-0032833-2},
  url = {http://www.ams.org/jourcgi/jour-getitem?pii=S0002-9947-1949-0032833-2},
  urldate = {2020-10-20},
  langid = {english},
  number = {3}
}

@book{fraenkelFoundationsSetTheory1958,
  title = {Foundations of {{Set Theory}}},
  author = {Fraenkel, Abraham A. and Bar-Hillel, Yehoshua},
  date = {1958},
  publisher = {{North Holland}},
  location = {{Amsterdam}}
}

@article{friedmanBarInductionPi1969,
  title = {Bar induction and {$\Pi^1_1$-$\mathsf{CA}$}},
  author = {Friedman, Harvey},
  date = {1969-09},
  journaltitle = {The Journal of Symbolic Logic},
  shortjournal = {J. Symbolic Logic},
  volume = {34},
  pages = {353--362},
  publisher = {{Association for Symbolic Logic}},
  issn = {0022-4812, 1943-5886},
  url = {https://projecteuclid.org/euclid.jsl/1183736849},
  urldate = {2020-10-26},
  langid = {english},
  mrnumber = {MR250877},
  number = {3},
  zmnumber = {0182.00904}
}

@article{friedmanCountableAlgebraSet1983,
  title = {Countable Algebra and Set Existence Axioms},
  author = {Friedman, Harvey and Simpson, Stephen G. and Smith, Rick L.},
  date = {1983-11-01},
  journaltitle = {Annals of Pure and Applied Logic},
  shortjournal = {Annals of Pure and Applied Logic},
  volume = {25},
  pages = {141--181},
  issn = {0168-0072},
  doi = {10.1016/0168-0072(83)90012-X},
  url = {http://www.sciencedirect.com/science/article/pii/016800728390012X},
  urldate = {2020-09-23},
  langid = {english},
  number = {2}
}

@inproceedings{friedmanFiniteCombinatorialPrinciple1982,
  title = {A Finite Combinatorial Principle Which Is Equivalent to the 1-Consistency of Predicative Analysis},
  booktitle = {Patras {{Logic Symposion}}},
  author = {Friedman, Harvey and McAloon, Kenneth and Simpson, Stephen G.},
  editor = {Metakides, George},
  date = {1982-01-01},
  volume = {109},
  pages = {197--230},
  publisher = {{North Holland}},
  doi = {10.1016/S0049-237X(08)71365-X},
  url = {http://www.sciencedirect.com/science/article/pii/S0049237X0871365X},
  urldate = {2020-10-27},
  langid = {english},
  series = {Studies in {{Logic}} and the {{Foundations}} of {{Mathematics}}}
}

@thesis{friedmanSubsystemsSetTheory1967,
  title = {Subsystems of Set Theory and Analysis},
  author = {Friedman, Harvey},
  date = {1967},
  institution = {{MIT}},
  url = {https://dspace.mit.edu/handle/1721.1/33486},
  urldate = {2020-10-26},
  langid = {english},
  type = {PhD thesis}
}

@inproceedings{friedmanSystemsSecondOrder1974,
  title = {Some Systems of Second Order Arithmetic and Their Use},
  booktitle = {Proceedings of the {{International Congress}} of {{Mathematicians}}},
  author = {Friedman, Harvey},
  date = {1974},
  pages = {235--242},
  location = {{Vancouver}}
}

@article{friedmanSystemsSecondOrder1976,
  title = {Systems of Second Order Arithmetic with Restricted Induction, {{I}}, {{II}} (Abstracts)},
  author = {Friedman, Harvey},
  date = {1976},
  journaltitle = {The Journal of Symbolic Logic},
  volume = {41},
  pages = {557--559}
}

@thesis{godelUberVollstandigkeitLogikkalkuls1929,
  title = {Über die Vollständigkeit des Logikkalküls},
  author = {Gödel, Kurt},
  date = {1929},
  institution = {{University of Vienna}},
  langid = {german},
  pagetotal = {33},
  type = {PhD thesis}
}

@article{goedelUeberFormalUnentscheidbare1931,
  title = {Über formal unentscheidbare Sätze der Principia Mathematica und verwandter Systeme I},
  author = {Gödel, Kurt},
  date = {1931-12-01},
  journaltitle = {Monatshefte für Mathematik und Physik},
  shortjournal = {Monatsh. f. Mathematik und Physik},
  volume = {38},
  pages = {173--198},
  issn = {1436-5081},
  doi = {10.1007/BF01700692},
  url = {https://doi.org/10.1007/BF01700692},
  urldate = {2020-10-25},
  langid = {german},
  number = {1}
}

@book{gordonIntegralsLebesgueDenjoy1994,
  title = {The {{Integrals}} of {{Lebesgue}}, {{Denjoy}}, {{Perron}}, and {{Henstock}}},
  author = {Gordon, Russell A.},
  date = {1994},
  publisher = {{American Mathematical Society}},
  eprint = {VN8RCgAAQBAJ},
  eprinttype = {googlebooks},
  isbn = {978-0-8218-3805-1},
  number = {4},
  pagetotal = {409},
  series = {Graduate {{Studies}} in {{Mathematics}}}
}

@article{gordonNonabsoluteIntegrationWorth1996,
  title = {Is Nonabsolute Integration Worth Doing?},
  author = {Gordon, Russell A.},
  date = {1996},
  journaltitle = {Real Analysis Exchange},
  volume = {22},
  pages = {23--33},
  publisher = {{Michigan State University Press}},
  issn = {0147-1937},
  doi = {10.2307/44152707},
  eprint = {44152707},
  eprinttype = {jstor},
  number = {1}
}

@book{heavisideElectromagneticTheory1893,
  title = {Electromagnetic {{Theory}}},
  author = {Heaviside, Oliver},
  date = {1893},
  volume = {1},
  publisher = {{The Electrician Printing and Publishing Company}},
  location = {{London}}
}

@book{henstockTheoryIntegration1963,
  title = {Theory of {{Integration}}},
  author = {Henstock, R.},
  date = {1963},
  publisher = {{Butterworths}},
  location = {{London}}
}

@book{hilbertGrundlagenGeometrie1899,
  title = {Grundlagen der Geometrie},
  author = {Hilbert, David},
  date = {1899},
  publisher = {{B.G. Teubner}},
  location = {{Leipzig}},
  eprint = {d4lKAAAAYAAJ},
  eprinttype = {googlebooks},
  langid = {german},
  pagetotal = {102}
}

@book{hilbertGrundlagenMathematik1934,
  title = {Grundlagen der Mathematik},
  author = {Hilbert, David and Bernays, Paul},
  date = {1934},
  publisher = {{Springer}},
  location = {{Berlin}},
  eprint = {5xkuAAAAMAAJ},
  eprinttype = {googlebooks},
  langid = {german},
  pagetotal = {498}
}

@incollection{hirstReverseMathematicsMatroids2017,
  title = {Reverse Mathematics of Matroids},
  booktitle = {Computability and {{Complexity}}: {{Essays Dedicated}} to {{Rodney G}}. {{Downey}} on the {{Occasion}} of {{His}} 60th {{Birthday}}},
  author = {Hirst, Jeffry L. and Mummert, Carl},
  editor = {Day, Adam and Fellows, Michael and Greenberg, Noam and Khoussainov, Bakhadyr and Melnikov, Alexander and Rosamond, Frances},
  date = {2017},
  pages = {143--159},
  publisher = {{Springer International Publishing}},
  location = {{Cham}},
  doi = {10.1007/978-3-319-50062-1_12},
  url = {https://doi.org/10.1007/978-3-319-50062-1_12},
  urldate = {2020-10-27},
  isbn = {978-3-319-50062-1},
  langid = {english},
  series = {Lecture {{Notes}} in {{Computer Science}}}
}

@article{jockuschDegreesMembersPi1972,
  title = {Degrees of members of {$\Pi^0_1$} classes},
  author = {Jockusch, Carl G. and Soare, Robert I.},
  date = {1972},
  journaltitle = {Pacific Journal of Mathematics},
  shortjournal = {Pacific J. Math.},
  volume = {40},
  pages = {605--616},
  publisher = {{Pacific Journal of Mathematics}},
  issn = {0030-8730},
  url = {https://projecteuclid.org/euclid.pjm/1102968559},
  urldate = {2020-10-27},
  mrnumber = {MR0309722},
  number = {3},
  zmnumber = {0209.02201}
}

@article{jockuschRamseyTheoremRecursion1972,
  title = {Ramsey's Theorem and Recursion Theory},
  author = {Jockusch, Carl G.},
  date = {1972-06},
  journaltitle = {The Journal of Symbolic Logic},
  volume = {37},
  pages = {268--280},
  publisher = {{Cambridge University Press}},
  issn = {0022-4812, 1943-5886},
  doi = {10.2307/2272972},
  url = {https://www.cambridge.org/core/journals/journal-of-symbolic-logic/article/ramseys-theorem-and-recursion-theory/8D9C0E8F7E0E78CF83651B1E622DAE84},
  urldate = {2020-10-27},
  langid = {english},
  number = {2}
}

@article{kechrisClassificationBaireClass1990,
  title = {A Classification of {{Baire}} Class 1 Functions},
  author = {Kechris, A. S. and Louveau, A.},
  date = {1990-03},
  journaltitle = {Transactions of the American Mathematical Society},
  shortjournal = {Transactions of the American Mathematical Society},
  volume = {318},
  pages = {209},
  issn = {00029947},
  doi = {10.2307/2001236},
  eprint = {2001236},
  eprinttype = {jstor},
  langid = {english},
  number = {1}
}

@article{kleeneRecursivePredicatesQuantifiers1943,
  title = {Recursive Predicates and Quantifiers},
  author = {Kleene, S. C.},
  date = {1943},
  journaltitle = {Transactions of the American Mathematical Society},
  shortjournal = {Trans. Amer. Math. Soc.},
  volume = {53},
  pages = {41--73},
  issn = {0002-9947, 1088-6850},
  doi = {10.1090/S0002-9947-1943-0007371-8},
  url = {https://www.ams.org/tran/1943-053-01/S0002-9947-1943-0007371-8/},
  urldate = {2020-10-27},
  langid = {english},
  number = {1}
}

@book{koComplexityTheoryReal1991,
  title = {Complexity {{Theory}} of {{Real Functions}}},
  author = {Ko, Ker-I.},
  date = {1991},
  publisher = {{Birkhäuser}},
  location = {{Boston}},
  doi = {10.1007/978-1-4684-6802-1},
  url = {https://www.springer.com/gp/book/9781468468045},
  urldate = {2020-08-11},
  isbn = {978-1-4684-6804-5},
  langid = {english},
  series = {Progress in {{Theoretical Computer Science}}}
}

@book{kurtzTheoriesIntegrationIntegrals2004,
  title = {Theories of {{Integration}}: {{The Integrals}} of {{Riemann}}, {{Lebesgue}}, {{Henstock}}-{{Kurzweil}}, and {{McShane}}},
  shorttitle = {Theories of {{Integration}}},
  author = {Kurtz, Douglas S. and Swartz, Charles W.},
  date = {2004},
  publisher = {{World Scientific}},
  eprint = {clhY9O__t6QC},
  eprinttype = {googlebooks},
  isbn = {978-981-238-843-8},
  langid = {english},
  pagetotal = {286}
}

@article{kurzweilGeneralizedOrdinaryDifferential1957,
  title = {Generalized Ordinary Differential Equations and Continuous Dependence on a Parameter},
  author = {Kurzweil, Jaroslav},
  date = {1957},
  journaltitle = {Czechoslovak Mathematical Journal},
  shortjournal = {Czech. Math. J.},
  volume = {7},
  pages = {418--449},
  issn = {0011-4642, 1572-9141},
  doi = {10.21136/CMJ.1957.100258},
  url = {https://dml.cz/handle/10338.dmlcz/100258},
  urldate = {2020-10-24},
  langid = {english},
  number = {3}
}

@article{kuyperEffectiveGenericityDifferentiability2014,
  title = {Effective Genericity and Differentiability},
  author = {Kuyper, Rutger and Terwijn, Sebastiaan},
  date = {2014-08-26},
  journaltitle = {Jou-rnal of Logic and Analysis},
  volume = {6},
  pages = {1--14},
  issn = {1759-9008},
  doi = {10.4115/jla.2014.6.4},
  url = {http://www.logicandanalysis.org/index.php/jla/article/view/215},
  urldate = {2020-10-27},
  langid = {english},
  number = {4}
}

@article{lebesgueIntegraleLongueurAire1902,
  title = {Intégrale, longueur, aire},
  author = {Lebesgue, Henri},
  date = {1902-12-01},
  journaltitle = {Annali di Matematica Pura ed Applicata},
  shortjournal = {Annali di Matematica, Serie III},
  volume = {7},
  pages = {231--359},
  issn = {0373-3114},
  doi = {10.1007/BF02420592},
  url = {https://doi.org/10.1007/BF02420592},
  urldate = {2020-10-26},
  langid = {french},
  number = {1}
}

@book{lebesgueLeconsIntegrationRecherche1904,
  title = {Leçons sur l'intégration et la recherche des fonctions primitives},
  shorttitle = {Leçons sur l'intégration et la recherche des fonctions primitives},
  author = {Lebesgue, Henri},
  date = {1904},
  publisher = {{Gaut-hier-Villars}},
  location = {{Paris}},
  eprint = {elBtAAAAMAAJ},
  eprinttype = {googlebooks},
  langid = {french},
  pagetotal = {156}
}

@book{leeHenstockKurzweilIntegrationEuclidean2011,
  title = {Henstock-{{Kurzweil Integration}} on {{Euclidean Spaces}}},
  author = {Lee, Tuo Yeong},
  date = {2011},
  publisher = {{World Scientific}},
  eprint = {cKPaMUmF4MAC},
  eprinttype = {googlebooks},
  isbn = {978-981-4324-58-8},
  langid = {english},
  pagetotal = {325}
}

@article{lobachevskyConciseOutlineFoundations1829,
  title = {A Concise Outline of the Foundations of Geometry},
  author = {Lobachevsky, Nikolai},
  date = {1829},
  journaltitle = {University of Kazan Messenger}
}

@article{luzinProprietesIntegraleDenjoy1912,
  title = {Sur Les Propriétés de l'intégrale de {{M}}. {{Denjoy}}},
  author = {Luzin, Nikolai},
  date = {1912},
  journaltitle = {Comptes rendus de l'Académie des Sciences},
  volume = {155},
  pages = {1475--1477}
}

@incollection{marconeLogicalStrengthNashWilliams1996,
  title = {On the Logical Strength of {{Nash}}-{{Williams}}' Theorem on Transfinite Sequences},
  booktitle = {Logic: From Foundations to Applications},
  author = {Marcone, Alberto},
  editor = {Hodges, Wilfrid and Hyland, Martin and Steinhorn, Charles and Truss, John},
  date = {1996-08-22},
  pages = {327--351},
  publisher = {{Clarendon Press}},
  location = {{New York}},
  isbn = {978-0-19-853862-2}
}

@book{markerModelTheoryIntroduction2002,
  title = {Model {{Theory}}: {{An Introduction}}},
  shorttitle = {Model {{Theory}}},
  author = {Marker, David},
  date = {2002},
  publisher = {{Springer-Verlag}},
  location = {{New York}},
  doi = {10.1007/b98860},
  url = {https://www.springer.com/gp/book/9780387987606},
  urldate = {2020-09-18},
  isbn = {978-0-387-98760-6},
  langid = {english},
  number = {217},
  series = {Graduate {{Texts}} in {{Mathematics}}}
}

@unpublished{mendiolaSignificanceAxiomChoice2016,
  title = {Significance of the Axiom of Choice in Mathematics},
  author = {Mendiola, José Eduardo},
  date = {2016-06-09},
  url = {https://www.researchgate.net/publication/303874180}
}

@inproceedings{metakidesIntroductionNonrecursiveMethods1982,
  title = {The Introduction of Non-Recursive Methods into Mathematics},
  booktitle = {The {{L}}. {{E}}. {{J}}. {{Brouwer Centenary Symposium}}, {{Noordwijkerhout}}},
  author = {Metakides, George and Nerode, Anil},
  editor = {Troelstra, A. S. and van Dalen, D.},
  date = {1982-01-01},
  volume = {110},
  pages = {319--335},
  publisher = {{Elsevier}},
  doi = {10.1016/S0049-237X(09)70135-1},
  url = {http://www.sciencedirect.com/science/article/pii/S0049237X09701351},
  urldate = {2020-10-25},
  langid = {english},
  options = {useprefix=true},
  series = {Studies in {{Logic}} and the {{Foundations}} of {{Mathematics}}}
}

@thesis{miletiPartitionTheoremsComputability2004,
  title = {Partition {{Theorems}} and {{Computability Theory}}},
  author = {Mileti, Joseph R.},
  date = {2004},
  institution = {{University of Illinois}},
  location = {{Urbana-Champaign}},
  url = {https://www.cambridge.org/core/product/identifier/S1079898600003140/type/journal_article},
  urldate = {2020-10-26},
  langid = {english},
  type = {PhD thesis}
}

@inproceedings{moriComputabilitiesFinecontinuousFunctions2001,
  title = {Computabilities of {{Fine}}-Continuous Functions},
  booktitle = {Computability and {{Complexity}} in {{Analysis}}},
  author = {Mori, Takakazu},
  editor = {Blanck, Jens and Brattka, Vasco and Hertling, Peter},
  date = {2001},
  pages = {200--221},
  publisher = {{Springer}},
  location = {{Berlin, Heidelberg}},
  doi = {10.1007/3-540-45335-0_13},
  isbn = {978-3-540-45335-2},
  langid = {english},
  series = {Lecture {{Notes}} in {{Computer Science}}}
}

@article{moriComputabilityWalshFunctions2002,
  title = {On the Computability of {{Walsh}} Functions},
  author = {Mori, Takakazu},
  date = {2002-07-28},
  journaltitle = {Theoretical Computer Science},
  shortjournal = {Theoretical Computer Science},
  volume = {284},
  pages = {419--436},
  issn = {0304-3975},
  doi = {10.1016/S0304-3975(01)00099-8},
  url = {http://www.sciencedirect.com/science/article/pii/S0304397501000998},
  urldate = {2020-10-20},
  langid = {english},
  number = {2}
}

@book{odifreddiClassicalRecursionTheory1999,
  title = {Classical {{Recursion Theory}}, {{Volume II}}},
  author = {Odifreddi, P.},
  date = {1999},
  publisher = {{North Holland}},
  eprint = {Z97uAAAAMAAJ},
  eprinttype = {googlebooks},
  isbn = {978-0-444-50205-6},
  langid = {english},
  number = {143},
  pagetotal = {970},
  series = {Studies in {{Logic}} and the {{Foundations}} of {{Mathematics}}}
}

@book{peanoArithmeticesPrincipiaNova1889,
  title = {Arithmetices principia: nova methodo exposita},
  shorttitle = {Arithmetices principia},
  author = {Peano, Giuseppe},
  date = {1889},
  publisher = {{Fratres Bocca}},
  location = {{Rome}},
  eprint = {UUFtAAAAMAAJ},
  eprinttype = {googlebooks},
  langid = {latin},
  pagetotal = {44}
}

@book{perronUberIntegralbegriff1914,
  title = {{\"U}ber den {I}ntegralbegriff},
  shorttitle = {Sitzungsberichte Der {{Heidelberger Akademie}} Der {{Wissenschaften}}, {{Mathematisch}}-{{Naturwissenschaftliche Klasse}}},
  author = {Perron, Oskar},
  date = {1914},
  publisher = {{Heidelberger Akademie der Wissenschaften}},
  doi = {10.11588/diglit.37437},
  url = {https://doi.org/10.11588/diglit.37437},
  series = {Sitzungsberichte}
}

@article{porterNotesComputableAnalysis2017,
  title = {Notes on Computable Analysis},
  author = {Porter, Michelle and Day, Adam and Downey, Rodney G.},
  date = {2017-01},
  journaltitle = {Theory of Computing Systems},
  shortjournal = {Theory Comput Syst},
  volume = {60},
  pages = {53--111},
  issn = {1432-4350, 1433-0490},
  doi = {10.1007/s00224-016-9732-y},
  url = {http://link.springer.com/10.1007/s00224-016-9732-y},
  urldate = {2020-10-27},
  langid = {english},
  number = {1}
}

@article{postRecursivelyEnumerableSets1944,
  title = {Recursively Enumerable Sets of Positive Integers and Their Decision Problems},
  author = {Post, Emil L.},
  date = {1944},
  journaltitle = {Bulletin of the American Mathematical Society},
  shortjournal = {Bull. Amer. Math. Soc.},
  volume = {50},
  pages = {284--316},
  issn = {0002-9904, 1936-881X},
  doi = {10.1090/S0002-9904-1944-08111-1},
  url = {https://www.ams.org/bull/1944-50-05/S0002-9904-1944-08111-1/},
  urldate = {2020-10-27},
  langid = {english},
  number = {5}
}

@book{pour-elComputabilityAnalysisPhysics1989,
  title = {Computability in {{Analysis}} and {{Physics}}},
  author = {Pour-El, Marian B. and Richards, J. Ian},
  date = {1989},
  publisher = {{Springer-Verlag}},
  location = {{Berlin}},
  number = {1},
  series = {Perspectives in {{Mathematical Logic}}}
}

@article{pour-elSimpleDefinitionComputable1975,
  title = {On a Simple Definition of Computable Function of a Real Variable—with Applications to Functions of a Complex Variable},
  author = {Pour-El, Marian B. and Caldwell, J.},
  date = {1975},
  journaltitle = {Zeitschrift für Mathematische Logik und Grundlagen der Mathematik},
  volume = {21},
  pages = {1--19},
  doi = {10.1002/malq.19750210102}
}

@thesis{riemannUberDarstellbarkeitFunction1854,
  title = {{\"U}ber die {D}arstellbarkeit einer {F}unction durch eine trigonometris-che {R}eihe},
  author = {Riemann, Bernhard},
  date = {1854},
  institution = {{{U}niversity of G{\"o}ttingen}},
  langid = {german},
  type = {Habilitation thesis}
}

@book{rogersTheoryRecursiveFunctions1967,
  title = {Theory of {{Recursive Functions}} and {{Effective Computability}}},
  author = {Rogers, Hartley},
  date = {1967},
  publisher = {{McGraw-Hill}},
  eprint = {nAkzAAAAMAAJ},
  eprinttype = {googlebooks},
  langid = {english},
  pagetotal = {526}
}

@unpublished{schwarzLecturesKarlWeierstrass1861,
  title = {Lectures of {{Karl Weierstrass}}},
  author = {Schwarz, H. A.},
  date = {1861},
  location = {{Mittag-Leffler Institute, Sweden}}
}

@article{shiojiFixedPointTheory1990,
  title = {Fixed Point Theory in Weak Second-Order Arithmetic},
  author = {Shioji, Naoki and Tanaka, Kazuyuki},
  date = {1990-05-22},
  journaltitle = {Annals of Pure and Applied Logic},
  shortjournal = {Annals of Pure and Applied Logic},
  volume = {47},
  pages = {167--188},
  issn = {0168-0072},
  doi = {10.1016/0168-0072(90)90068-D},
  url = {http://www.sciencedirect.com/science/article/pii/016800729090068D},
  urldate = {2020-10-27},
  langid = {english},
  number = {2}
}

@article{shoreInvariantsBooleanAlgebras2005,
  title = {Invariants, {B}oolean algebras and {$\ACA^+$}},
  author = {Shore, Richard A.},
  date = {2005-04-13},
  journaltitle = {Transactions of the American Mathematical Society},
  shortjournal = {Trans. Amer. Math. Soc.},
  volume = {358},
  pages = {989--1014},
  issn = {0002-9947},
  doi = {10.1090/S0002-9947-05-03802-X},
  url = {http://www.ams.org/journal-getitem?pii=S0002-9947-05-03802-X},
  urldate = {2020-09-25},
  langid = {english},
  number = {03}
}

@inbook{simpsonSubsystemsReverseMathematics1987,
  title = {Subsystems of {$Z_2$} and reverse mathematics},
  booktitle = {Proof {{Theory}}},
  author = {Simpson, Stephen G.},
  date = {1987},
  edition = {2},
  publisher = {{North Holland}},
  bookauthor = {Takeuti, Gaisi},
  number = {81},
  series = {Studies in {{Logic}} and {{Foundations}} of {{Mathematics}}}
}

@book{simpsonSubsystemsSecondOrder2009,
  title = {Subsystems of {{Second Order Arithmetic}}},
  author = {Simpson, Stephen G.},
  date = {2009},
  edition = {2},
  publisher = {{Cambridge University Press}},
  location = {{Cambridge}},
  doi = {10.1017/CBO9780511581007},
  url = {https://www.cambridge.org/core/books/subsystems-of-second-order-arithmetic/EA16CB4305831530B7015D6BC46B7424},
  urldate = {2020-08-11},
  isbn = {978-0-521-88439-6},
  series = {Perspectives in {{Logic}}}
}

@book{soareRecursivelyEnumerableSets1987,
  title = {Recursively {{Enumerable Sets}} and {{Degrees}}},
  shorttitle = {Recursively {{Enumerable Sets}} and {{Degrees}}},
  author = {Soare, Robert I.},
  date = {1987},
  publisher = {{Springer-Verlag}},
  location = {{Berlin, Heidelberg}},
  url = {https://www.springer.com/gp/book/9783540666813},
  urldate = {2020-10-18},
  isbn = {978-3-540-66681-3},
  langid = {english},
  series = {Perspectives in {{Mathematical Logic}}}
}

@book{soareTuringComputability2016,
  title = {Turing {{Computability}}},
  shorttitle = {Turing {{Computability}}},
  author = {Soare, Robert I.},
  date = {2016},
  publisher = {{Springer-Verlag}},
  location = {{Berlin, Heidelberg}},
  doi = {10.1007/978-3-642-31933-4},
  url = {https://www.springer.com/gp/book/9783642319327},
  urldate = {2020-10-18},
  isbn = {978-3-642-31932-7},
  langid = {english},
  series = {Theory and {{Applications}} of {{Computability}}}
}

@article{solomonPiCAOrderTypes2001,
  title = {{$\PiCA$} and order types of countable ordered groups},
  author = {Solomon, Reed},
  date = {2001},
  journaltitle = {The Journal of Symbolic Logic},
  volume = {66},
  pages = {192--206},
  publisher = {{Association for Symbolic Logic}},
  issn = {0022-4812},
  doi = {10.2307/2694917},
  eprint = {2694917},
  eprinttype = {jstor},
  number = {1}
}

@article{solovayHyperarithmeticallyEncodableSets1978,
  title = {Hyperarithmetically Encodable Sets},
  author = {Solovay, Robert M.},
  date = {1978},
  journaltitle = {Transactions of the American Mathematical Society},
  shortjournal = {Trans. Amer. Math. Soc.},
  volume = {239},
  pages = {99--122},
  issn = {0002-9947, 1088-6850},
  doi = {10.1090/S0002-9947-1978-0491103-7},
  url = {https://www.ams.org/tran/1978-239-00/S0002-9947-1978-0491103-7/},
  urldate = {2020-10-27},
  langid = {english}
}

@article{speckerNichtKonstruktivBeweisbare1949,
  title = {Nicht Konstruktiv Beweisbare {{Sätze}} Der {{Analysis}}},
  author = {Specker, Ernst},
  date = {1949},
  journaltitle = {The Journal of Symbolic Logic},
  volume = {14},
  pages = {145--158},
  publisher = {{[Association for Symbolic Logic, Cambridge University Press]}},
  issn = {0022-4812},
  doi = {10.2307/2267043},
  eprint = {2267043},
  eprinttype = {jstor},
  number = {3}
}

@thesis{steelDeterminatenessSubsystemsAnalysis1977,
  title = {Determinateness and Subsystems of Analysis},
  author = {Steel, John Robert},
  date = {1977},
  institution = {{University of California}},
  location = {{Berkeley}},
  url = {https://catalog.hathitrust.org/Record/101649807},
  urldate = {2020-10-27},
  pagetotal = {iii, 107 l.},
  type = {PhD thesis}
}

@article{tanakaWeakAxiomsDeterminacy1991,
  title = {Weak axioms of determinacy and subsystems of analysis {II}: {$\Delta^0_2$} games},
  author = {Tanaka, Kazuyuki},
  date = {1991-04},
  journaltitle = {Annals of Pure and Applied Logic},
  shortjournal = {Annals of Pure and Applied Logic},
  volume = {52},
  pages = {181--193},
  issn = {01680072},
  doi = {10.1016/0168-0072(91)90045-N},
  url = {https://linkinghub.elsevier.com/retrieve/pii/016800729190045N},
  urldate = {2020-10-27},
  langid = {english},
  number = {1-2}
}

@thesis{traylorEquivalenceAxiomChoice1962,
  title = {On the Equivalence of the Axiom of Choice, {{Zorn}}'s Lemma, and the Well-Ordering Theorem},
  author = {Traylor, Grace Joy},
  date = {1962},
  institution = {{Atlanta University}},
  langid = {english},
  type = {Masters thesis}
}

@article{turingComputableNumbersApplication1937,
  title = {On Computable Numbers, with an Application to the {{Entscheidungsproblem}}},
  author = {Turing, Alan},
  date = {1937},
  journaltitle = {Proceedings of the London Mathematical Society},
  volume = {s2-42},
  pages = {230--265},
  url = {https://londmathsoc.onlinelibrary.wiley.com/doi/abs/10.1112/plms/s2-42.1.230},
  urldate = {2020-10-25},
  number = {1}
}

@article{vaughtAlfredTarskiWork1986,
  title = {Alfred {{Tarski}}'s Work in Model Theory},
  author = {Vaught, Robert L.},
  date = {1986},
  journaltitle = {The Journal of Symbolic Logic},
  volume = {51},
  pages = {869--882},
  publisher = {{[Association for Symbolic Logic, Cambridge University Press]}},
  issn = {0022-4812},
  doi = {10.2307/2273900},
  eprint = {2273900},
  eprinttype = {jstor},
  number = {4}
}

@article{walshClosedSetNormal1923,
  title = {A Closed Set of Normal Orthogonal Functions},
  author = {Walsh, J. L.},
  date = {1923},
  journaltitle = {American Journal of Mathematics},
  volume = {45},
  pages = {5--24},
  publisher = {{Johns Hopkins University Press}},
  issn = {0002-9327},
  doi = {10.2307/2387224},
  eprint = {2387224},
  eprinttype = {jstor},
  number = {1}
}

@article{zermeloUberGrenzzahlenUnd1930a,
  title = {{\"U}ber {G}renzzahlen und {M}engenbereiche: neue {U}ntersuchungen {\"u}ber die {G}rundlagen der {M}engenlehre},
  author = {Zermelo, Ernst},
  date = {1930},
  journaltitle = {Fundamenta Mathematicae},
  volume = {16},
  pages = {29--47}
}

\end{document}